\documentclass[12pt]{amsart}
\usepackage{amssymb}  
\usepackage{latexsym} 
\usepackage[all]{xy}
\usepackage{comment}
\usepackage{url}

\newcommand{\Z}{{\mathbb Z}}

\newcommand{\C}{{\mathbb C}}
\newcommand{\Q}{{\mathbb Q}}
\newcommand{\Qbar}{{\overline{\Q}}}
\newcommand{\rhobar}{{\overline{\rho}}}

\newcommand{\PP}{{\mathbb P}}
\newcommand{\FF}{{\mathcal F}}
\newcommand{\GG}{{\mathcal G}}
\newcommand{\CL}{{\mathcal L}}
\newcommand{\OO}{{\mathcal O}}

\newcommand{\HH}{{\mathcal H}}

\newcommand{\Kbar}{\overline{K}}
\newcommand{\Lbar}{\overline{L}}

\newcommand{\BF}{{\mathbf{F}}}
\newcommand{\BG}{{\mathbf{G}}}

\newcommand{\GL}{\operatorname{GL}}

\newcommand{\SL}{\operatorname{SL}}
\newcommand{\PGL}{\operatorname{PGL}}
\newcommand{\Gal}{\operatorname{Gal}}

\newcommand{\Aut}{\operatorname{Aut}}

\newcommand{\eps}{\varepsilon}
\newcommand{\Gm}{\mathbb{G}_{\text{\rm m}}}
\newcommand{\isom}{\cong}
\newcommand{\inj}{\hookrightarrow}
\def\x{{\mathbf x}}

\newcommand{\rank}{{\operatorname{rank}}}

\newcommand{\qq}{\mathfrak{q}}
\newcommand{\rr}{\mathfrak{r}}

\newcommand{\e}{{\mathbf e}}
\newcommand{\f}{{\mathbf f}}
\newcommand{\uu}{{\mathbf u}}
\newcommand{\vv}{{\mathbf v}}
\newcommand{\zz}{\zeta}
\newcommand{\MMPQ}{\langle M_P,M_Q \rangle}
\def\I{{\mathcal I}}

\def\MatP{{\begin{pmatrix}
    1   &    0   &    0    & \cdots &      0     \\
    0   & \zeta_n  &    0    & \cdots &      0     \\
    0   &    0   & \zeta_n^2 & \cdots &      0     \\ 
 \vdots & \vdots & \vdots  &        &   \vdots   \\
    0   &   0    &    0    & \cdots & \zeta_n^{n-1}  
\end{pmatrix}}}

\def\MatQ{{\begin{pmatrix}
    0   &    0   & \cdots &   0    &      1     \\ 
    1   &    0   & \cdots &   0    &      0     \\
    0   &    1   & \cdots &   0    &      0     \\ 
 \vdots & \vdots &        & \vdots &   \vdots   \\
    0   &   0    & \cdots &   1    &      0       
\end{pmatrix}}}
\def\MatR{{\begin{pmatrix}
    1   &    0   & \cdots  &   0    &      0     \\
    0   &    0   & \cdots  &   0    &      1     \\
    0   &    0   & \cdots  &   1    &      0     \\ 
 \vdots & \vdots &         & \vdots &   \vdots   \\
    0   &    1   & \cdots  &   0    &      0       
\end{pmatrix}}}
\def\PSL{{\operatorname{PSL}}}
\def\Diag{\operatorname{Diag}}

\def\la{{\lambda}}
\def\lala{{{\boldsymbol \la}}}
\def\dv{{\mid}}

\def\Im{{\rm Im}}

\newcommand{\DD}{{\mathbf D}}
\newcommand{\cc}{{\mathbf c}}

\newfont{\wncyr}{wncyr10 at 12pt}
\newfont{\wncyrten}{wncyr10 at 10pt}

\newenvironment{Proof}{\par\noindent{\sc Proof:}}%
                      {\hspace*{\fill}\nobreak$\Box$\par\medskip}
\newenvironment{ProofOf}[1]{\par\noindent{\sc Proof of #1:}}%
                       {\hspace*{\fill}\nobreak$\Box$\par\medskip}

\newtheorem{Proposition}{Proposition}[section]
\newtheorem{Theorem}[Proposition]{Theorem}
\newtheorem{Lemma}[Proposition]{Lemma}
\newtheorem{Corollary}[Proposition]{Corollary}

\theoremstyle{definition}

\newtheorem{Definition}[Proposition]{Definition}
\newtheorem{Remark}[Proposition]{Remark}

\newtheorem{Example}[Proposition]{Example}

\begin{document}
\date{9th May 2011}
\title[On families of $n$-congruent elliptic curves]
{On families of $n$-congruent elliptic curves}

\author{T.A.~Fisher}
\address{University of Cambridge,
          DPMMS, Centre for Mathematical Sciences,
          Wilberforce Road, Cambridge CB3 0WB, UK}
\email{T.A.Fisher@dpmms.cam.ac.uk}

\renewcommand{\baselinestretch}{1.1}
\renewcommand{\arraystretch}{1.3}

\renewcommand{\theenumi}{\roman{enumi}}

\begin{abstract}
We use an invariant-theoretic method to compute certain twists
of the modular curves $X(n)$ for $n=7,9,11$. Searching for rational
points on these twists enables us to find non-trivial pairs of 
$n$-congruent elliptic curves over $\Q$, i.e. pairs of non-isogenous
elliptic curves over $\Q$ whose $n$-torsion subgroups are isomorphic as
Galois modules. We also show by giving explicit non-trivial examples 
over $\Q(T)$ that there are infinitely many examples over $\Q$ 
in the cases $n=9$ and $n=11$.
\end{abstract}

\maketitle

\section{Introduction}

Elliptic curves $E_1$ and $E_2$ over a field $K$ are 
{\em $n$-congruent} if their $n$-torsion subgroups $E_1[n]$ and 
$E_2[n]$ are isomorphic as Galois modules. They are 
{\em directly $n$-congruent}
if the isomorphism $\phi : E_1[n] \isom E_2[n]$ respects the Weil pairing
and {\em reverse $n$-congruent} if $e_n(\phi P,\phi Q) = e_n(P,Q)^{-1}$
for all $P,Q \in E_1[n]$. The elliptic curves directly $n$-congruent 
to a given elliptic curve $E$ are parametrised by the modular curve 
$Y_E(n) = X_E(n) \setminus \{ \text{cusps} \}$. 

For $n \le 5$ we have $X_E(n) \isom \PP^1$ and the corresponding 
families of elliptic curves were computed by Rubin and Silverberg 
\cite{RubinSilverberg}, \cite{RubinSilverberg2}, \cite{Silverberg}. 
However for $n \ge 7$ the genus is greater than $1$.
This prompted Mazur \cite{MazurQ} to ask whether there are any pairs of 
non-isogenous elliptic curves over $\Q$ that are directly $n$-congruent for any
$n \ge 7$. This was answered by Kraus and Oesterl\'e \cite{KO} 
who gave the example of the directly $7$-congruent elliptic curves 
$152a1$ and $7448e1$. The labels here are those 
in Cremona's tables \cite{Cr}. Nowadays
it is easy to find further examples by searching in 
Cremona's tables, for example
\[ \begin{array}{l@{\qquad}rcl} 
n=11 & 190b1 &+& 2470a1, \\
n=13 & 52a2 &+& 988b1, \\
n=17 & 3675b1 &-& 47775b1. 
\end{array} \]
In each case the $n$-congruence is proved by computing sufficiently many
traces of Frobenius; see \cite[Proposition 4]{KO}. Then 
\cite[Proposition 2]{KO} shows that the congruences are direct
or reverse as indicated by the $\pm$.

Motivated by Mazur's question, Kani and Schanz~\cite{KS} studied the
geometry of the surfaces that parametrise pairs of $n$-congruent of 
elliptic curves. This prompted them to conjecture that for any $n \le 12$ 
there are infinitely many pairs of $n$-congruent non-isogenous elliptic curves 
over $\Q$. It is understood that we are looking for examples 
with distinct $j$-invariants, since otherwise from any single example 
we could construct infinitely many
by taking quadratic twists. The conjecture was proved in the case $n=7$ 
by Halberstadt and Kraus \cite{HK1}, who subsequently \cite{HK} gave
an explicit formula for $X_E(7)$ and used it to show that 
there are infinitely many $6$-tuples of 
directly $7$-congruent non-isogenous elliptic curves over $\Q$. 
We find the corresponding formulae for $X_E(9)$ and $X_{E}(11)$ and
use them to construct explicit infinite families of pairs of elliptic curves
that prove the conjecture 
for $n=9$ and $n=11$. In contrast the proof of the conjecture for 
$n=11$ in \cite{KR} does not yield a single explicit example.

We briefly mention three further motivations for studying 
$n$-congruence of elliptic curves.
\begin{itemize}
\item The modular approach to solving Diophantine equations sometimes
requires us to find all elliptic curves $n$-congruent to a given 
elliptic curve. For example the paper of Poonen, Schaefer and Stoll 
\cite{PSS} makes essential use of the formula for $X_E(7)$ due to 
Halberstadt and Kraus.
\item There is a correspondence between pairs of reverse $n$-congruent
elliptic curves and curves of genus $2$ that admit a degree $n$ morphism
to an elliptic curve. See for example \cite{Frey}.
\item It was observed by Cremona and Mazur \cite{CM} that if elliptic
curves $E$ and $F$ are $n$-congruent then the Mordell-Weil group of
$F$ can sometimes be used to explain elements of the Tate-Shafarevich
group of $E$. 
\end{itemize}

As each of these motivations makes clear, we should also be interested
in congruences that do not respect the Weil pairing. The elliptic curves 
reverse $n$-congruent to $E$ are parametrised by the modular curve 
$Y^-_E(n) = X^-_E(n) \setminus \{ \text{cusps} \}$. The families of 
elliptic curves parametrised by $Y_E^-(3)$ and $Y_E^-(4)$ were 
computed in \cite{g1hess}, and the case of $Y_E^-(5)$ will be treated 
in \cite{enqI}. An equation for
$X_E^-(7)$ was given in \cite[Section 7.2]{PSS}. We find corresponding
formulae for $X_E^-(9)$ and $X_E^-(11)$. In the cases $n=7$ and $n=9$ 
we use these formulae to construct explicit infinite families of 
pairs of reverse $n$-congruent non-isogenous elliptic curves over $\Q$. 
We do not know if any such families exist in the case $n=11$. 

In Section~\ref{sec:defcrvs} we recall the definitions of 
$X(n)$ and its twists. We then state the formulae for $X_E(n)$
and $X_E^-(n)$ for $n=3,7,9,11$ in Section~\ref{sec:statres}. We include
the cases $n=3$ and $n=7$ since our methods are quite different from 
those in \cite{g1hess} and \cite{HK}. The results in the case $n=3$ 
are also needed in Section~\ref{modint} to treat the case $n=9$. In the next
two subsections we explain our basic strategy for computing twists and
illustrate how it works in the case $n=3$.

In Section~\ref{modcurve} we recall Klein's equations for $X(n)$ for
$n \ge 5$ an odd integer. The original approach of Klein was via theta
functions, but our treatment is purely algebraic. 
We also recall explicit formulae for the action of $\SL_2(\Z/n\Z)$ 
on $X(n)$. Then in Section~\ref{invthy} we use invariant theory for 
$\SL_2(\Z/n\Z)$ to compute the twists $X_E(n)$ and $X^-_E(n)$ for $n=7,9,11$.

In Section~\ref{sec:diag} we show that in the special case of an elliptic
curve $E$ whose $n$-torsion contains a copy of the Galois module $\mu_n$,
there are particularly simple formulae for $X_E(n)$ and $X^-_E(n)$. 
This section may be read independently of Section~\ref{invthy}.
Again the result for $X_E(7)$ is already given in \cite{HK}.

Our formulae reduce the problem of finding elliptic curves
$n$-congruent to $E$ to that of finding rational
points on $X_E(n)$ and $X^-_E(n)$. However before 
searching for rational points it helps to simplify the equations by
making a change of co-ordinates. In Section~\ref{sec:minred} 
we describe how to do this in the case $K = \Q$. The problem 
naturally splits into two parts called minimisation and reduction.

In  Section~\ref{modint} we give formulae for the families of elliptic 
curves parametrised by $Y_E(n)$ and $Y_E^-(n)$ for $n=7,9,11$. 
This is easiest in the case $n=9$ since the natural maps $X_E(9) \to X_E(3)$
and $X^-_E(9) \to X^-_E(3)$ have a simple geometric description. 
In the cases $n=7$ and $n=11$ it is easy to 
compute the maps $j : X_E(n) \to \PP^1$ and $j : X^-_E(n) \to \PP^1$. 
However determination of the right quadratic twists takes 
considerably more work. (In specific numerical examples 
one can always fall back on the method in \cite{KO}, \cite{HK1}.)
In the case of $Y_E(7)$ a formula is given in \cite{HK}, although
this formula does not quite cover all cases. We give a new proof 
leading to formulae that work in all cases. We then generalise to 
the families of elliptic curves 
parametrised by $Y_E^-(7)$, $Y_E(11)$ and $Y^-_E(11)$.

Finally in Section~\ref{sec:ex} we give examples of two different sorts.
First we have written a program in Magma \cite{Magma} 
that given an elliptic curve 
$E/\Q$ and $n \in \{7,9,11\}$ searches for rational points 
(up to a specified height bound) 
on minimised and reduced models for $X_E(n)$ and $X_E^-(n)$  and 
returns the corresponding list of elliptic curves $n$-congruent to $E$. 
For $n=9$ and $n=11$ we have run this program on every elliptic curve in 
the Cremona database. In particular we found three triples of directly 
$9$-congruent non-isogenous elliptic curves over $\Q$.

Secondly we have found some non-trivial pairs of $n$-congruent
elliptic curves over $\Q(T)$ for $n=7,9,11$. The infinite families
mentioned above are obtained by specialising $T$.

All computer calculations in support of this work were preformed
using Magma \cite{Magma}. A Magma file checking all our formulae, 
together with extended versions of the tables in Section~\ref{sec:ex},
is available from the author's website \cite{mywebsite}.

\subsection{Some modular curves}
\label{sec:defcrvs}

We work over a field $K$ of characteristic $0$ and write $\Kbar$ for
the algebraic closure. Let $n \ge 3$ be an integer, and $M$ 
a Galois module isomorphic to $(\Z/n\Z)^2$ as an abelian group
and equipped with a non-degenerate alternating Galois equivariant
pairing $M \times M \to \mu_n$. We temporarily write $Y_M$ for the
algebraic curve defined over $K$ whose $L$-rational points ($L$ a
field extension of $K$) parametrise the isomorphism classes of pairs
$(E,\phi)$, where $E$ is an elliptic curve defined over $L$ and $\phi : 
E[n] \isom M$ is a symplectic isomorphism (i.e. one that matches up the
given pairing on $M$ with the Weil pairing on $E[n]$) commuting with
the action of $\Gal(\Lbar/L)$.
Two such pairs $(E_1,\phi_1)$ and $(E_2,\phi_2)$ are 
isomorphic if there is an $L$-isomorphism $\alpha : E_1 \to E_2$ such
that $\phi_1 = \phi_2 \circ (\alpha|_{E_1[n]})$.

Let $X_M$ be the smooth projective model of $Y_M$. We define $X(n)$
to be $X_M$ in the case $M = \mu_n \times \Z/n\Z$ with pairing
\[  \langle (\zeta,a) , (\xi,b) \rangle = \zeta^{b} \xi^{-a}. \]
Given an elliptic curve $E$ over $K$ we define $X_E(n)$
to be $X_M$ in the case $M$ is $E[n]$ equipped with the Weil pairing.
More generally we may take $M$ to be $E[n]$ equipped with the $r$th
power of the Weil pairing for any $r \in (\Z/n\Z)^\times$. However
the curve obtained only depends on the class of $r$ mod squares.
Since we will later be specialising to $n \in \{3,7,9,11\}$ it 
suffices to take $r = \pm 1$. Taking $r = 1$
gives the curve $X_E(n)$ defined above. Taking $r=-1$ gives the 
curve $X^-_E(n)$.

We identify $\SL_2(\Z/n\Z)$ with the group of symplectic
automorphisms of $\mu_n \times \Z/n\Z$. There is then a natural
action of $\PSL_2(\Z/n\Z) := \SL_2(\Z/n\Z)/ \{\pm I_2\}$ 
on $X(n)$ with quotient map
$j : X(n) \to \PP^1$. 
From the analytic theory we know  
that the $j$-map is ramified above 0, 1728 and $\infty$ with 
ramification indexes 3, 2 and $n$. Hence by the Riemann-Hurwitz
formula the genus of $X(n)$ is 
\[ g(n) = \frac{n-6}{12n} \, \# \PSL_2(\Z/n\Z) +1 \]
where for $n \ge 3$ we have 
$\# \PSL_2(\Z/n\Z) = (n^3/2) \prod_{p \dv n} (1-1/p^2)$. 
For some small values of $n$ the genus is as follows.
\[ \begin{array}{c|ccccccccccccccccccccc}
n & 2 & 3 & 4 & 5 & 6 & 7 & 8 & 9 & 10 & 11 & 12 & 13 & 14 & 15 
& 16 & 17 
\\ \hline
g(n) & 0 & 0 & 0 & 0 & 1 & 3 & 5 & 10 & 13 & 26 & 25 & 50 & 49 & 73 
& 81 & 133 
\end{array} \]

\subsection{Statement of results}
\label{sec:statres}

The family of elliptic curves parametrised by $Y_E(n)$ 
for $n=2,3,4,5$ is given by formulae of Rubin and Silverberg 
\cite{RubinSilverberg}, \cite{RubinSilverberg2}, \cite{Silverberg}.
In \cite{g1hess} we developed an alternative invariant-theoretic
approach that also gives formulae for $Y_E^-(n)$. In the case
$n=3$ we have

\begin{Theorem}
\label{MainThm3}
Let $E$ be an elliptic curve 
with Weierstrass equation $y^2 = x^3 - 27 c_4  x - 54 c_6$.
Then the families of elliptic curves parametrised by $Y_E(3)$ and $Y_E^-(3)$ are 
\[ y^2 =  x^3 - 27 \cc_4(\la,\mu)  x - 54 \cc_6(\la,\mu) \]
and
\[ y^2 =  x^3 - 27 \cc_4^{*}(\la,\mu)  x - 54 \cc_6^{*}(\la,\mu) \]
where
\begin{align*}
\cc_4(\la,\mu) & =
c_4 \la^4 + 4 c_6 \la^3 \mu + 6 c_4^2 \la^2 \mu^2 
+ 4 c_4 c_6 \la \mu^3 - (3 c_4^3 - 
    4 c_6^2) \mu^4, \\
 \cc_6(\la,\mu) & =
c_6 \la^6 + 6 c_4^2 \la^5 \mu + 15 c_4 c_6 \la^4 \mu^2 
  + 20 c_6^2 \la^3 \mu^3 \\ &~\qquad + 15 c_4^2 c_6 \la^2 \mu^4 
+ 6(3 c_4^4 - 2 c_4 c_6^2) \la \mu^5 + (9 c_4^3 c_6 - 8 c_6^3) \mu^6,
\end{align*}
and 
\begin{align*}
\cc_4^{*}(\la,\mu) & = -4 (\la^4 - 6 c_4 \la^2 \mu^2 - 8 c_6 \la \mu^3 
 - 3 c_4^2 \mu^4 ) /(c_4^3-c_6^2), \\
\cc_6^{*}(\la,\mu) & = -8 \cc_6(\la,\mu)/(c_4^3-c_6^2)^2.
\end{align*}
\end{Theorem}
\begin{Proof} These formulae (and their relation to those in 
\cite{RubinSilverberg}) may be found in \cite{g1hess}.
We recall the idea of the proof. Let $F$ be the ternary cubic obtained 
by homogenising the Weierstrass equation for $E$. Then the family 
of curves parametrised by $X_E(3)$ is the pencil of plane 
cubics spanned by $F$ and its Hessian, each with base point the same 
as that for $E$. It only remains to put these curves in 
Weierstrass form. One way of carrying out this 
final calculation is by using the formula for the Jacobian 
(due to Weil) coming from classical invariant theory. For this 
reason the polynomials $\cc_4(\la,\mu)$ and $\cc_6(\la,\mu)$ may 
already be found in \cite[Art.~230]{Salmon}, although not with the 
interpretation given here.  

The corresponding formula for $Y^-_E(3)$ 
is obtained by replacing the Hessian (which is a covariant) by suitable
contravariants.
\end{Proof}

A formula for $X_E(7)$ was obtained by Halberstadt and Kraus
\cite{HK}. Their method relies on studying the points on the Klein
quartic \[X(7) = \{x^3 y + y^3 z + z^3 x = 0\} \subset \PP^2\]
corresponding to an elliptic curve $E$ and the elliptic curves $E_a$,
$E_b$, $E_c$ that are $2$-isogenous to $E$. By combining this 
result with some classical invariant theory Poonen, Schaefer and Stoll
\cite[Section 7.2]{PSS} then gave a formula for $X^-_E(7)$.
\begin{Theorem}[Halberstadt, Kraus, Poonen, Schaefer, Stoll]
\label{MainThm7}
Let $E$ be an elliptic curve 
with Weierstrass equation $y^2 = x^3 + a x + b$.
Then $X_E(7) \subset \PP^2$ has equation $\mathcal F = 0$ where
\begin{align*}
\FF =  a x^4 &+ 7 b x^3 z + 3 x^2 y^2 - 3 a^2 x^2 z^2 - 6 b x y z^2 \\
& ~- 5 a b x z^3 + 2 y^3 z + 3 a y^2 z^2 + 2 a^2 y z^3 - 4 b^2 z^4 ,
\end{align*}
and $X^-_E(7)\subset \PP^2 $ has equation $\mathcal G = 0$ where
\begin{align*}
\GG = -a^2 & x^4 + 2 a b x^3 y - 12 b x^3 z - (6 a^3 + 36 b^2) x^2 y^2 
    + 6 a x^2 z^2 \\ &~+ 2 a^2 b x y^3 - 12 a b x y^2 z + 18 b x y z^2 
    + (3 a^4 + 19 a b^2) y^4 \\ &~- (8 a^3 + 42 b^2) y^3 z 
      + 6 a^2 y^2 z^2 - 8 a y z^3 + 3 z^4.
\end{align*}
\end{Theorem}

We give new proofs of Theorems~\ref{MainThm3} and~\ref{MainThm7}. 
We then extend to the cases $n=9$ and $n=11$.
Although we believe that the formulae in Theorems~\ref{MainThm9} 
and~\ref{MainThm11} below are correct for all elliptic curves
$E$ we assume for simplicity that $j(E) \not= 0,1728$.
In Section~\ref{modcurve} we recall 
that $X(9)$ may be embedded in $\PP^3$ as the
complete intersection of two cubics: 
\[ X(9) = \left\{ 
\begin{aligned}
x^2 y +y^2 z +z^2 x & =  0 \\
x y^2 + y z^2 + z x^2 & =  t^3 
\end{aligned} \right\} \subset \PP^3. \]
\begin{Theorem}
\label{MainThm9}
Let $E$ be an elliptic curve with Weierstrass equation 
$y^2 = x^3 + a x + b$.
If $j(E) \not= 0,1728$ then $X_E(9) \subset \PP^3$ has equations
${\mathcal F}_1 = {\mathcal F}_2 = 0$ where
\begin{align*}
 {\mathcal F}_1  & =  
    x^2 t +  6 x y z + 6 b x t^2 + 6 y^3 - 9 a y^2 t + 6 a^2 y t^2 - 3 b z^3 + 
        3 a^2 z^2 t \\ & \qquad  + 9 a b z t^2  - (a^3 - 12 b^2) t^3, \\
 {\mathcal F}_2  & = 
    x^2 z  +  6 x y^2 - 6 a x y t + 2 a^2 x t^2 - 9 a y^2 z - 18 b y z^2 + 
        12 a^2 y z t  + a^2 z^3 \\ & \qquad 
+ 9 a b z^2 t - 3 a^3 z t^2 + a^2 b t^3,
\end{align*}
and $X^-_E(9) \subset \PP^3$ has equations $ {\mathcal G}_1 = {\mathcal G}_2 = 0$ where
\begin{align*}
 {\mathcal G}_1 & =  
    9 x^2 y + 3 x^2 z - 6 a x y t + 6 b x t^2 - 6 a y^3 
      + 27 b y^2 t + 3 a y z^2 \\ & \qquad + 
        18 b y z t + 3 a^2 y t^2 + a z^3 + 3 b z^2 t + a^2 z t^2 - a b t^3, \\
  {\mathcal G}_2 & =    x^3 + 6 a x y z + 18 b x y t + 3 a x z^2 + 6 b x z t 
  + a^2 x t^2 + 9 b y^3 + 6 a^2 y^2 t \\ & \qquad - 9 b y z^2 + 6 a^2 y z t - 3 a b y t^2 
  - 4 b z^3 + 2 a^2 z^2 t + 2 b^2 t^3.
\end{align*}
\end{Theorem}
It was observed by Klein that $X(11)$ may be embedded 
in $\PP^4$ as the singular locus of the Hessian of the cubic threefold
\[ \{ v^2 w + w^2 x + x^2 y + y^2 z + z^2 v = 0 \} \subset \PP^4.  \]

\begin{Theorem}
\label{MainThm11}
Let $E$ be an elliptic curve with Weierstrass equation $y^2 = x^3 + a x + b$.
If $j(E) \not= 0,1728$ then $X_E(11) \subset \PP^4$ is the singular locus of
the Hessian of 
\begin{align*}
{\mathcal F} &=  v^3 + a v^2 z - 2 a v x^2 + 4 a v x y - 6 b v x z + a v y^2 + 6 b v y z 
     + a^2 v z^2 - w^3 \\ & + a w^2 z - 4 a w x^2 - 12 b w x z + a^2 w z^2 - 2 b x^3 
     + 3 b x^2 y + 2 a^2 x^2 z + 6 b x y^2 \\ & + 4 a b x z^2 + b y^3 - a^2 y^2 z 
     + a b y z^2 + 2 b^2 z^3,
\end{align*}
and $X^-_E(11) \subset \PP^4$ is the singular locus of the Hessian of 
\begin{align*}
{\mathcal G} = v^2 & z + 2 v w y + 4 v x y + 2 w^2 x - a w^2 z + 2 w x^2 - 2 a w y^2 
      - 6 b w y z \\ & + 6 x^3 - 6 a x^2 z + 2 a^2 x z^2 + b y^3 - 2 a^2 y^2 z 
      - 5 a b y z^2 - b^2 z^3.
\end{align*}
\end{Theorem}

\subsection{Preliminaries on twisting}
We identify $\SL_2(\Z/n\Z)$ with the group of 
symplectic automorphisms of $\mu_n \times \Z/n\Z$ via 
$(\begin{smallmatrix} a & b \\ c & d \end{smallmatrix}) :
(\zeta_n^x, y) \mapsto (\zeta_n^{a x + by},c x + dy)$,
where $\zeta_n$ is a fixed primitive $n$th root of unity.
The Galois action on $\SL_2(\Z/n\Z)$ will be that arising from
this identification. There is then a Galois equivariant 
group homomorphism
\begin{equation}
\label{rhobar-v1}
 \rhobar : \SL_2(\Z/n\Z) \to \Aut ( X(n) ). 
\end{equation}
where the action of $\rhobar(\gamma)$ on $Y(n)$ is given by
$(E,\phi) \mapsto (E, \gamma \circ \phi)$.

\begin{Lemma}
\label{firstlem}
Let $E/K$ be an elliptic curve and $\phi : E[n] \isom \mu_n \times \Z/n\Z$
a symplectic, respectively anti-symplectic, isomorphism 
over $\Kbar$. Then there is an isomorphism $\alpha : X_E(n) \to X(n)$,
respectively $\alpha : X^-_E(n) \to X(n)$, over $\Kbar$ such that
\[ \sigma(\alpha) \alpha^{-1} = \rhobar ( \sigma(\phi) \phi^{-1}) \]
for all $\sigma \in \Gal(\Kbar/K)$.
\end{Lemma}
\begin{Proof}
The points on $Y_E(n)$, respectively $Y^-_E(n)$, correspond to 
pairs $(F,\psi)$ where $F$ is an elliptic curve
and $\psi : F[n] \isom E[n]$ is a symplectic, respectively 
anti-symplectic, isomorphism. The modular interpretation of
$\alpha$ is that it sends $(F,\psi)$ to $(F, \phi \circ \psi)$.
\end{Proof}

We are interested in the case $X(n)$ is embedded in $\PP^{m-1}$
for some $m$, and $\rhobar$ is realised as a projective representation 
(also denoted $\rhobar$ by abuse of notation):
\begin{equation*}
\label{rhobar-v2}
 \rhobar : \SL_2(\Z/n\Z) \to \PGL_m(\Kbar). 
\end{equation*}
Writing $\propto$ for equality in $\PGL_m(\Kbar)$, 
we further suppose there 
exists\footnote{By Remark~\ref{rem:rho}(iii) below we can usually 
take $\eps = I_m$.} $\eps \in \GL_m(K)$  such that
\begin{equation}
\label{iota}
\rhobar (\iota \gamma \iota) \propto 
\eps \rhobar(\gamma)^{-T} \eps^{-1}
\end{equation}
for all $\gamma \in \SL_2(\Z/n\Z)$, where
$\iota = (\begin{smallmatrix} 1 & 0 \\ 0 & -1 \end{smallmatrix})$.
Our strategy for computing $X_E(n)$ and $X_E^-(n)$ as twists
of $X(n)$ is explained by the following lemma.
 
\begin{Lemma}
\label{TwistLem}
Let $E/K$ be an elliptic curve and $\phi : E[n] \isom \mu_n \times \Z/n\Z$
a symplectic isomorphism over $\Kbar$.
Suppose $h_1,h_2 \in \GL_m(\Kbar)$ satisfy
\[\sigma(h_1) h_1^{-1} \propto \rhobar ( \sigma (\phi) \phi^{-1})
\qquad \quad 
\sigma(h_2) h_2^{-1} \propto \rhobar ( \sigma (\phi) \phi^{-1})^{-T}
\] for 
all $\sigma \in \Gal(\Kbar/K)$. Then $X_E(n) \subset \PP^{m-1}$ and
$X^-_E(n) \subset \PP^{m-1}$ are the twists of $X(n) \subset \PP^{m-1}$
given by $ X_E(n) \isom X(n) ;\, \x \mapsto h_1 \x$
and $ X^-_E(n) \isom X(n) ;\, \x \mapsto \eps h_2 \x.$
\end{Lemma}
\begin{Proof}
We apply Lemma~\ref{firstlem} to the pairs $(E,\phi)$ and 
$(E, \iota \circ \phi)$.
\end{Proof}

\begin{Remark}
If the projective representation $\rhobar$ lifts to a 
representation 
\[ \rho : \SL_2(\Z/n\Z) \to \GL_m(\Kbar) \]
then the existence of the matrix $h_1$ in Lemma~\ref{TwistLem} is 
clear from the generalised form of Hilbert's Theorem 90 which 
states that $H^1(\Gal(\Kbar/K), \GL_m(\Kbar))=0$. We could then take
$h_2 = h_1^{-T}$. We will instead use invariant theory for the group 
$\SL_2(\Z/n\Z)$ to give explicit formulae for $h_1$ and $h_2$.
\end{Remark}

\subsection{Introductory example: the case $n=3$}
\label{sec:n=3}
The family of curves parametrised 
by $X(3) \isom \PP^1$ is the Hesse pencil of plane cubics
\begin{equation*}
   \{ A(x^3 + y^3 + z^3) + B x y z = 0 \} \subset \PP^2 
\end{equation*}
with base point $(x:y:z) = (0:1:-1)$. 
An equation for this family in Weierstrass form is
$y^2 = x^3 - 27 c_4(A,B) x - 54 c_6(A,B)$ with discriminant 
$\Delta = (c_4^3 - c_6^2)/1728 = D^3$ where
\begin{align*}
D(A,B) &= -27 A^4 - A B^3 \\
c_4(A,B) &= -216 A^3 B + B^4 \\
c_6(A,B) &= 5832 A^6 - 540 A^3 B^3 - B^6.
\end{align*}

The action of $\SL_2(\Z/3\Z)$ on $X(3) \isom \PP^1$ 
lifts to a representation
\[ \rho : \SL_2(\Z/3\Z) \to \SL_2(\Kbar) \]  
given on the generators 
$S = (\begin{smallmatrix} 0 & 1 \\ -1 & 0 \end{smallmatrix})$
and $T = (\begin{smallmatrix} 1 & 1 \\ 0 & 1 \end{smallmatrix})$ by
\[  \rho(S) =
\frac{-1}{\zeta_3 - \zeta_3^2}
\begin{pmatrix} 1 & 1/3 \\ 6 & -1 \end{pmatrix}, \qquad
\rho(T) = 
\begin{pmatrix} \zeta_3 & 0 \\ 0 & \zeta_3^2 \end{pmatrix}. \]
The hypothesis~(\ref{iota}) is satisfied with $\eps = 
(\begin{smallmatrix} 1 & 0 \\ 0 & 18 \end{smallmatrix})$.

\begin{Lemma}
\label{prop3}
Let $E/K$ be an elliptic curve and 
$\phi : E[3] \isom \mu_3 \times \Z/3\Z$ a symplectic isomorphism
over $\Kbar$.
Let $(A:B)$ be the corresponding $\Kbar$-point on $X(3) \isom \PP^1$
with co-ordinates $(A,B)$ scaled so that
\begin{equation}
\label{c4c6:three}
c_4(A,B) = c_4(E) \quad \text{ and } \quad
c_6(A,B) = c_6(E) 
\end{equation}
where $E$ has Weierstrass equation $y^2 = x^3 - 27 c_4(E) x - 54 c_6(E)$.
If $j(E) \not= 0,1728$ then  
\[ h_1 = \begin{pmatrix} A & - \frac{\partial D}{\partial B} \\
B & \frac{\partial D}{\partial A} \end{pmatrix} \]
satisfies the hypotheses of Lemma~\ref{TwistLem}.
\end{Lemma}

\begin{Proof}
Since $(A:B) \in X(3)$ is not a cusp we have $\det h_1 = 4 D \not=0$
and so $h_1 \in \GL_2(\Kbar)$.
Let $G = \Im(\rho)$. There is a unique character 
$\chi: G \to \mu_3$ such that
\begin{equation}
\label{inv3}
\begin{aligned}
D \circ g & = \chi(g) D \\
c_4 \circ g & = \chi(g)^2 c_4 \\
c_6 \circ g & = \chi(g)^3  c_6 
\end{aligned}
\end{equation}
for all $g \in G$.
It follows by the chain rule that 
\begin{equation}
\label{col3}
   \left( \begin{array}{r} \!\!\!  - \frac{\partial D}{\partial B} \\
 \frac{\partial D}{\partial A}  \end{array} \right) \circ g = 
\chi(g) g
\left( \begin{array}{r} \!\!\! - \frac{\partial D}{\partial B} \\
 \frac{\partial D}{\partial A}  \end{array} \right)
\end{equation}
for all $g \in G$.

Let $\xi_\sigma = \sigma(\phi) \phi^{-1} \in \SL_2(\Z/3\Z)$. 
Since $\rho$ describes the action of $\SL_2(\Z/3\Z)$ on $X(3) \isom \PP^1$
we have
\begin{equation}
\label{ddaggerthree}
\sigma \begin{pmatrix} A \\ B \end{pmatrix} 
 = \lambda_\sigma \rho( \xi_\sigma) \begin{pmatrix} A \\ B \end{pmatrix} 
\end{equation}
for some $\lambda_\sigma \in \Kbar^\times$. Then by~(\ref{inv3}) we have
\[ \sigma (c_4(A,B)) = \lambda_\sigma^{4} \chi_\sigma^{2} c_4(A,B) \quad
\text{ and } \quad
\sigma (c_6(A,B)) = \lambda_\sigma^{6} \chi_\sigma^{3} c_6(A,B) \]
where $\chi_\sigma = \chi(\rho(\xi_\sigma))$.
Since $c_4(E), c_6(E) \in K$ it follows by~(\ref{c4c6:three}) 
and our assumption 
$j(E) \not= 0,1728$ that $\lambda_\sigma^{2} \chi_\sigma = 1$. 
Using~(\ref{col3}) and~(\ref{ddaggerthree}) we compute
\[ \sigma(h_1) = h_1 \circ (\lambda_\sigma \rho(\xi_\sigma)) 
= \rho(\xi_\sigma) h_1 \begin{pmatrix} \lambda_\sigma & 0 \\
0 & \lambda_\sigma^3 \chi_\sigma  \end{pmatrix}  
\propto \rho(\xi_\sigma) h_1. 
\] 
Hence $\sigma(h_1) h_1^{-1} \propto \rho(\xi_\sigma)$ 
as required.
\end{Proof}

Lemma~\ref{TwistLem} shows that $X_E(3)$ and $X_E^-(3)$ are the twists 
of $X(3) \isom \PP^1$ by $h_1$ and $\eps h_1^{-T}$. 
To compute the families of curves parametrised by $X_E(3)$ 
and $X_E^-(3)$ we put
\begin{align*} 
\cc_4(\lambda,\mu) & =  c_4(\lambda A - 
\mu \tfrac{\partial D}{\partial B},\lambda B + 
\mu \tfrac{\partial D}{\partial A}), \\
\cc_6(\lambda,\mu) & =  c_6(\lambda A - 
\mu \tfrac{\partial D}{\partial B},\lambda B + 
\mu \tfrac{\partial D}{\partial A}), 
\end{align*}
and
\begin{align*} 
\cc_4^*(\la,\mu) & =  \frac{1}{(36D)^4} c_4(
\lambda B + \mu \tfrac{\partial D}{\partial A}, 
-18 (\lambda A - \mu \tfrac{\partial D}{\partial B})), \\
\cc_6^*(\la,\mu) & =  \frac{1}{(36D)^6} c_6(
\lambda B + \mu \tfrac{\partial D}{\partial A}, 
-18 (\lambda A - \mu \tfrac{\partial D}{\partial B})).
\end{align*}
The coefficients of $\cc_4(\lambda,\mu)$, $\cc_6(\lambda,\mu)$,
$\cc_4^*(\la,\mu)$, $\cc_6^*(\la,\mu)$ 
may be written as rational functions in $c_4(A,B)$ and $c_6(A,B)$, 
thereby giving the formulae in Theorem~\ref{MainThm3}. 
The following alternative proof of Theorem~\ref{MainThm3} has the
advantage over that in \cite{g1hess} of not requiring any 
knowledge of the Hessian or the contravariants. We include it
to illustrate the approach we take for $n=7,9,11$.

\medskip

\begin{ProofOf}{Theorem~\ref{MainThm3}}
We assume $j(E) \not= 0,1728$ so that Lemma~\ref{prop3} applies. 
By \cite[Proposition 2.1]{RubinSilverberg} the formulae in the 
statement of the theorem define a family of elliptic curves
with constant mod $3$ Galois representation. In the case of 
$Y_E(3)$ the proof is completed by specialising to $(\la:\mu) = (1:0)$.
In general, that is, to give an argument that also applies to $Y^-_E(3)$,
we note that since we have specified a family of elliptic curves with
the right $j$-invariant, it must be correct up to quadratic twist, say by 
$\delta \in K^\times$. It remains to show that $\delta$ is a square.
As noted in \cite[Section 7.1]{HK} it suffices to prove this in the case 
$\phi : E[3] \isom \mu_3 \times \Z/3\Z$ is defined over $K$. But in
that case the families of curves parametrised by $X(3)$, $X_E(3)$
and $X_E^-(3)$ are the same, and this is born out by our construction.
\end{ProofOf}

\section{Equations for $X(\lowercase{n})$}
\label{modcurve}

We recall equations of Klein \cite{K1}, \cite{K2}, \cite{K3} 
for the modular curves $X(n)$. Our treatment follows the survey in 
\cite[Chapter 4]{mythesis}, but see also \cite{AR}, \cite{V2}.

\subsection{Klein's equations}
\label{klein-xn}

Let $E$ be an elliptic curve and $P,Q$ a basis for $E[n]$ 
with $e_n(P,Q) = \zeta_n$. If we embed $E \subset \PP^{n-1}$ 
by a complete linear system $|D|$ of degree~$n$ then the translation
maps $\tau_P$ and $\tau_Q$ extend to automorphisms of
$\PP^{n-1}$. The following lemma is recalled from \cite{JEMS}.

\begin{Lemma}
\label{tlab}
(i) We may change co-ordinates on $\PP^{n-1}$ (over $\Kbar$) so
that $\tau_P$ and $\tau_Q$ are given by
\[ M_P =\MatP \quad \text{ and } \quad M_Q =\MatQ.  \]
(ii) If $n$ is odd and $[-1]^*D \sim D$  then there is a unique choice of 
co-ordinates (over $\Kbar$) such that $\tau_P$, $\tau_Q$ and multiplication
by $-1$ are given by $M_P$, $M_Q$ and
\begin{eqnarray*}
[-1] = \MatR.   
\end{eqnarray*}
\end{Lemma}

We restrict to $n \ge 5$ an odd integer.
Given a symplectic isomorphism $\phi : E[n] \isom \mu_n \times \Z/n\Z$
we let $P,Q$ be the basis for $E[n]$ with $\phi(iP+jQ) = (\zeta_n^i,j)$
for all $i,j \in \Z/n\Z$.
If we embed $E \subset \PP^{n-1}$ via the complete 
linear system $|n.0_E|$ and 
choose co-ordinates as in Lemma~\ref{tlab} then  
sending $(E,\phi)$ to the image of $0_E$ defines
an embedding $Y(n) \subset \PP^{n-1}$. Indeed if we know the 
co-ordinates of $0_E \in \PP^{n-1}$ then 
$M_P$ and $M_Q$ allow us to write down $n^2$ points on $E$.
Since $E$ is defined by quadrics, it follows by B\'ezout's theorem
that these $n^2$ points suffice to determine $E$.
The embedding $Y(n) \subset \PP^{n-1}$
is defined over $K$ 
since ${\MMPQ} \subset \PGL_n(\Kbar)$ is isomorphic 
to $\mu_n \times \Z/n\Z$ as a Galois module.

We write $(x_0: x_1 : \ldots: x_{n-1})$ for our co-ordinates on 
$\PP^{n-1}$ and agree to read all subscripts mod $n$.
Since $n$ is odd we have 
\[ n.0_E \sim 0_E + P + 2P + \ldots + (n-1)P. \]
Therefore $0_E$ belongs to exactly one of the hyperplanes fixed 
by $M_P$. But $0_E$ is fixed by $[-1]$ so we have either
\[ \begin{array}{rrclr}
            & \quad 0 & = & (0:a_1:a_2: \ldots :a_2:a_1) & (+) \\
\text{ or } & 0 & = & (0:a_1:a_2: \ldots :-a_2:-a_1) & \quad \quad \quad (-) 
\end{array} \]
where $a_1,a_2,\ldots $ are non-zero. 
 
Let $W$ be the vector space of quadrics on $\PP^{n-1}$ and $V$ the
subspace of quadrics vanishing on $E$. Then $\dim W = n(n+1)/2$
and $\dim V = n(n-3)/2$. The action of $M_P$ allows us to write these
as direct sums $V =\oplus V_i$ and $W= \oplus W_i$ with
 \[ V_i \subset W_i = \langle x_i^2, x_{i-1} x_{i+1}, \ldots \rangle.\]
Since $n$ is odd it follows by the action of $M_Q$ that 
$\dim V_i = (n-3)/2$ and $\dim W_i = (n+1)/2$. The requirement
that the quadrics in $V_0$ vanish at $0_E$ and its translates under $M_Q$
imposes some linear conditions on the coefficients of these quadrics.
This leads us to rule out the case $(+)$ and to make the following definition.
\begin{Definition}
\label{def:zn}
For $n \ge 5$ an odd integer let $Z(n) \subset \PP^{n-1}$ be the subvariety
defined by $a_0=0$, $a_{n-i}=-a_i$ and
\begin{equation}
\label{K1}
 \rank ( a_{i-j} a_{i+j} )_{i,j=0}^{n-1} \le 2.
\end{equation}
\end{Definition}

We note that~(\ref{K1}) is equivalent to the vanishing of the 
$4 \times 4$ Pfaffians of this skew-symmetric matrix.
The above construction shows that $Y(n) \subset Z(n)$. It is 
natural to ask whether $Z(n) = X(n)$. V\'elu \cite{V2} proved this 
in the case $n=p$ is a prime. However when $n$ is composite
$Z(n)$ has extra components.

When $n=7$ we put $0_E=(0: a: b: -c : c : -b: -a)$. Then $Z(7)$ 
is defined by  
\[ \rank \begin{pmatrix}
  0   & -a^2 & -b^2 & -c^2 \\
  a^2 &    0 &   ac &  -bc \\
  b^2 &  -ac &    0 &   ab \\
  c^2 &   bc &  -ab &    0 
\end{pmatrix} \le 2. \]
Thus $X(7) = Z(7)$ is the Klein quartic 
$ \{ a^3b+b^3c+c^3a=0 \} \subset \PP^2$.

When $n=9$ we put
$0_E= (0: a: -b: d : c : -c : -d : b : -a)$. Then
$Z(9)$ is defined by 
\[ \rank \begin{pmatrix}
  0   & -a^2 & -b^2 & -d^2 &-c^2 \\
  a^2 &    0 &  -ad &   bc &  cd \\
  b^2 &   ad &    0 &   ac & -bd \\
  d^2 &  -bc &  -ac &    0 & -ab \\
  c^2 &  -cd &   bd &   ab &   0 \\
\end{pmatrix} \le 2, \]
equivalently
\begin{align*}
(a^2 b+b^2 c+c^2 a) d & = 0 & b c^3 - b a^3 - c d^3 & = 0 \\
a b^3 - a c^3 - b d^3 & = 0 & c a^3 - c b^3 - a d^3 & = 0.
\end{align*}
Adding together the last three equations 
and factoring shows that $Z(9)$ is the union of the curve
\[ X(9) = \left\{ 
\begin{aligned}
a^2 b +b^2 c +c^2 a & = 0 \\
a b^2 + b c^2 + c a^2 & = d^3 
\end{aligned} \right\} \subset \PP^3 \]
and four isolated points
\[ (0:0:0:1),  \, \, (1:1:1:0), \, \, 
(1:\zeta_3:\zeta_3^2:0), \, \, (1:\zeta_3^2:\zeta_3:0). \]

When $n=11$ we put $0_E= (0: a: -c: b: e: d: -d: -e: -b: c: -a)$. 
Then $X(11) \subset \PP^4$ 
is the singular locus of the Hessian of the cubic threefold
\[ \{ a^2 b + b^2 c + c^2 d + d^2 e + e^2 a = 0 \} \subset \PP^4. \]
We refer to \cite{AR} for further details. 
We checked using Magma that the homogeneous 
ideal of $X(11)$ is generated by the $4 \times 4$ minors of the 
Hessian matrix of the cubic form.

\subsection{The action of $\SL_2(\Z/n\Z)$} 
\label{sec:gpaction}

We continue to take $n \ge 5$ odd. In Section~\ref{klein-xn} we defined
an embedding $X(n) \subset \PP^{m-1}$ where $m = (n-1)/2$. The action
of $\SL_2(\Z/n\Z)$ on $X(n)$ extends to a projective representation
$\overline{\rho}: \SL_2(\Z/n\Z) \to \PGL_m(\Kbar)$. We show
that it lifts to a representation. See \cite[Appendix I]{AR}
for a discussion of how this relates to work of Weil.

\begin{Proposition}
\label{rho-lifts}
The projective representation
$\overline{\rho}: \SL_2(\Z/n\Z) \to \PGL_m(\Kbar)$ lifts to a 
representation $\rho: \SL_2(\Z/n\Z) \to \GL_m(\Kbar)$.
\end{Proposition}
\begin{Proof}
The action of $\gamma = (\begin{smallmatrix} a & b \\ c & d 
\end{smallmatrix}) \in \SL_2(\Z/n\Z)$ on $Y(n)$ is given by
\begin{equation}
\label{G-act}
\gamma : (E,P,Q) \mapsto (E,d P - c Q, -b P + aQ).
\end{equation}  

If we embed $X(n) \subset \PP^{n-1}$ as 
described in Section~\ref{klein-xn} then 
the action~(\ref{G-act}) extends to 
a projective representation $\overline{\pi} :\SL_2(\Z/n\Z) 
\to \PGL_n(\Kbar)$ where the image of 
$\gamma = (\begin{smallmatrix} a & b \\ c & d \end{smallmatrix})$
is uniquely determined by the properties that 
\begin{equation}
\label{comm-rels}
   \overline{\pi}(\gamma)^{-1} M_P^u M_Q^v \overline{\pi}(\gamma) 
   \propto M_P^{du - bv} M_Q^{-cu+av} 
\end{equation}
for all $u,v \in \Z/n\Z$, and $\overline{\pi}(\gamma)$ 
commutes with $[-1]$. We regard $\overline{\pi}$
as describing an action 
on $\PP^{n-1} = \PP(W)$ where $W$ is an 
$n$-dimensional vector space. The action of $[-1]$ gives an eigenspace
decomposition $W = W_+ \oplus W_-$ with $\dim W_{\pm} = (n \pm 1)/2$.
We may then identify $\overline{\rho}$ with the restriction of 
$\overline{\pi}$ to $\PP(W_-) = \PP^{m-1}$. To prove the proposition
we show more generally that $\overline{\pi}$ lifts to a representation
$\pi :\SL_2(\Z/n\Z) \to \GL_n(\Kbar)$. 

Let $S = (\begin{smallmatrix} 0 &  1 \\ -1 & 0 \end{smallmatrix})$
and $T = (\begin{smallmatrix} 1 &  1 \\ 0 & 1 \end{smallmatrix})$
be the usual generators of $\SL_2(\Z/n\Z)$. In view of the
relations $(ST)^3 = S^4 = T^n = I_2$, the only $1$-dimensional characters
of $\SL_2(\Z/n\Z)$ are the ones in the case $n$ is a multiple of $3$
that factor via $\PSL_2(\Z/3\Z) \isom A_4$. Using~(\ref{G-act}) and 
(\ref{comm-rels}) we compute
\begin{equation}
\label{get_SandT}
\overline{\pi}(S) \propto (\zeta_n^{ij})_{i,j=0}^{n-1} \quad \qquad
\overline{\pi}(T) \propto \Diag (\zeta_n^{i^2/2})_{i=0}^{n-1}
\end{equation}
where the exponents are read as elements of $\Z/n\Z$.

If $M \in \GL_n(\Kbar)$ acts on $W_{\pm}$ then 
we write $M_\pm$ for its restriction to $W_\pm$.
Since \[3(1^2 + 2^2 + \ldots + m^2) \equiv 0 \pmod{n}\] it is clear
by~(\ref{get_SandT}) that there is a lift $\pi(T)$ of $\overline{\pi}(T)$ 
and $1$-dimensional characters $\chi_{\pm}$ of $\SL_2(\Z/n\Z)$
such that $\det(\pi(T)_\pm) = \chi_{\pm}(T)$.
Next we lift $\overline{\pi}(S)$ to a matrix $\pi(S)$ such that
\begin{equation}
\label{lift-S}
  \pi(S) \pi(T)^{-1} \pi(S) = \pi(T) \pi(S) \pi(T).  
\end{equation}
Restricting to $W_{\pm}$ and taking determinants it follows that
$\det(\pi(S)_\pm) = 1 = \chi_{\pm}(S)$. 
For each $\gamma \in \SL_2(\Z/n\Z)$ we now let $\pi(\gamma)$ be the
unique lift of $\overline{\pi}(\gamma)$ such that
$\det(\pi(\gamma)_\pm) = \chi_{\pm}(\gamma)$. These lifts exist since
$S$ and $T$ generate $\SL_2(\Z/n\Z)$ and are unique since $\dim W_-$
and $\dim W_+$ are coprime. It is evident that the map $\pi$ so 
defined is a group homomorphism.
\end{Proof}

\begin{Remark}
\label{rem:rho}
(i) A calculation using~(\ref{lift-S}) shows that 
$\pi(S) = g_n^{-1} (\zeta^{ij})_{i,j=0}^{n-1}$
where the Gauss sum $g_n = \sum_{i=0}^{n-1} \zeta^{-i^2/2}$ 
satisfies $g_n^2 = (-1)^{(n-1)/2} n$. \\ 
(ii) If we take $0_E = (0:a_1:a_2: \ldots :-a_2:-a_1)$ then with
respect to co-ordinates $(a_1: \ldots: a_m)$ we may take 
\begin{equation*}
\rho(S) = g_n^{-1} (\zeta^{ij}-\zeta^{-ij})_{i,j=1}^{m} \quad \qquad
\rho(T) = \Diag (\zeta^{i^2/2})_{i=1}^{m}.
\end{equation*}
In particular $\rho(-I_2) = (-1)^{(n+1)/2} I_m$. \\
(iii) Since $\rho(S)$ and $\rho(T)$ are symmetric 
we see that~(\ref{iota}) holds with $\eps = I_m$. \\
(iv) If $3 \nmid n$ then $\SL_2(\Z/n\Z)$ has no $1$-dimensional characters
and so the lifts we have constructed are unique. If $3 \mid n$ 
then $3 \nmid m$ and we can make $\rho$ unique by demanding that
$\det \rho(T) = 1$, equivalently $\rho$ takes values in $\SL_m(\Kbar)$.
\end{Remark}

\section{Equations for $X_E(n)$ and $X^-_E(n)$}
\label{invthy}

We derive our equations for $X_E(n)$ and $X^-_E(n)$ by using  
invariant theory for the group $\SL_2(\Z/n\Z)$ to twist the equations
for $X(n)$ in Section~\ref{klein-xn}. We split into the
cases $n=7,9,11$.

\subsection{Formulae in the case $n=7$}
\label{sec:formulae7}
We recall that $X(7)$ is the Klein quartic $\{ F =0 \} \subset \PP^2$ 
where $F = a^3 b + b^3 c + c^3 a$. 
Let $G \isom \PSL_2(\Z/7\Z)$ be the image of 
$\rho : \SL_2(\Z/7\Z) \to \GL_3(\Kbar)$. It is generated by
\[ \frac{1}{g_7} \begin{pmatrix} 
\zz_7     - \zz_7^{6} & \zz_7^{2} - \zz_7^{5} &  \zz_7^{4} - \zz_7^{3} \\
\zz_7^{2} - \zz_7^{5} & \zz_7^{4} - \zz_7^{3} &  \zz_7 - \zz_7^{6}   \\
\zz_7^{4} - \zz_7^{3} & \zz_7 - \zz_7^{6} & \zz_7^{2} - \zz_7^{5} 
\end{pmatrix}
\quad \text{ and } \quad 
\begin{pmatrix}
\zz_7 & 0 & 0 \\
0 & \zz_7^4 & 0  \\
0 & 0 & \zz_7^2 
\end{pmatrix} \]
where 
$g_7 = 1 + 2(\zeta_7^3 + \zeta_7^5 + \zeta_7^6)$.

\begin{Definition}
An invariant of degree $m$ is a homogeneous polynomial $I = I(a,b,c)$ of degree
$m$ such that $I \circ g = I$ for all $g \in G$. 
\end{Definition}
Recalling formulae of Klein we put
\[  H = (-1/54) \times \left| \begin{array}{ccc}  \vspace{1ex}
\frac{\partial^2 F}{\partial a^2} & 
\frac{\partial^2 F}{\partial a \partial b} &
\frac{\partial^2 F}{\partial a \partial c} \\ \vspace{1ex}
\frac{\partial^2 F}{\partial a \partial b} &
\frac{\partial^2 F}{\partial b^2} & 
\frac{\partial^2 F}{\partial b \partial c} \\
\frac{\partial^2 F}{\partial a \partial c} &
\frac{\partial^2 F}{\partial b \partial c} &
\frac{\partial^2 F}{\partial c^2} 
\end{array} \right|,  \]
\[  c_4 = (1/9) \times \left| \begin{array}{cccc}  \vspace{1ex}
\frac{\partial^2 F}{\partial a^2} & 
\frac{\partial^2 F}{\partial a \partial b} &
\frac{\partial^2 F}{\partial a \partial c} & 
\frac{\partial H}{\partial a } \\ \vspace{1ex}
\frac{\partial^2 F}{\partial a \partial b} &
\frac{\partial^2 F}{\partial b^2} & 
\frac{\partial^2 F}{\partial b \partial c} &
\frac{\partial H}{\partial b}  \\ \vspace{1ex}
\frac{\partial^2 F}{\partial a \partial c} &
\frac{\partial^2 F}{\partial b \partial c} &
\frac{\partial^2 F}{\partial c^2} &
\frac{\partial H}{\partial c}  \\
\frac{\partial H}{\partial a} &
\frac{\partial H}{\partial b} &
\frac{\partial H}{\partial c} &
0
\end{array} \right|,  
\qquad
  c_6 = (1/14) \times \left| \begin{array}{ccc}  \vspace{1ex}
\frac{\partial F}{\partial a} & 
\frac{\partial F}{\partial b} &
\frac{\partial F}{\partial c} \\ \vspace{1ex}
\frac{\partial H}{\partial a} &
\frac{\partial H}{\partial b} & 
\frac{\partial H}{\partial c} \\
\frac{\partial c_4}{\partial a} &
\frac{\partial c_4}{\partial b} &
\frac{\partial c_4}{\partial c} 
\end{array} \right|.  \]
The ring of invariants $K[a,b,c]^G$ is generated by 
$F, H, c_4$ and $c_6$ subject to a single relation
which reduces when we set $F=0$ to
\[ c_4^3 - c_6^2 \equiv 1728 H^7 \pmod{F}. \] 
Since $F, H, c_4$ and $c_6$ have degrees $4, 6, 14$ and $21$ it is clear
that every invariant of odd degree is divisible by $c_6$.

\begin{Lemma}
\label{lem:j7}
The $j$-invariant $X(7) \to \PP^1$ is given by $j = c_4^3/H^7$.
\end{Lemma}
\begin{Proof}
Both $j$ and $j_0 = c_4^3/H^7$ define maps $X(7) \to \PP^1$ 
that quotient out by the action of $G \isom \PSL_2(\Z/7\Z)$. 
So they can differ by at most a M\"obius map. We recall that $j$
is ramified above $0$, $1728$ and $\infty$ with ramification
indices $3$, $2$ and $7$. Since
\begin{align*}
\# \{ F = c_4 = 0 \} &\le 4 \deg (c_4) = \tfrac{1}{3} |G| \\
\# \{ F = c_6 = 0 \} &\le 4 \deg (c_6) = \tfrac{1}{2} |G| \\
\# \{ F = H = 0 \} &\le 4 \deg (H) = \tfrac{1}{7} |G| 
\end{align*}
and $j_0 - 1728 = c_6^2/H^7$ it follows that $j=j_0$ as required.
\end{Proof}

\begin{Definition}
\label{defc7}
A covariant column, respectively contravariant column, of degree $m$ 
is a column vector $\vv=(v_1,v_2,v_3)^T$ of homogeneous polynomials
of degree $m$ in variables $a,b,c$ such that $\vv \circ g = g \vv$, 
respectively $\vv \circ g = g^{-T} v$, for all $g \in G$.
\end{Definition}

We note that $\x = (a,b,c)^T$ is a covariant column of degree $1$,
whereas if $I$ is an invariant of degree $m$ then 
$\nabla I = (\tfrac{\partial I}{\partial a},\tfrac{\partial I}{\partial b},
\tfrac{\partial I}{\partial c})^T$ is a contravariant column of
degree $m-1$.

\begin{Lemma}
\label{Lem7}
Let $E/K$ be an elliptic curve and $\phi : E[7] \isom \mu_7 \times \Z/7\Z$ 
a symplectic isomorphism over $\Kbar$. Let $(a:b:c)$ be the
corresponding $\Kbar$-point on $X(7) \subset \PP^2$ with co-ordinates
$(a,b,c)$ scaled so that 
\begin{equation}
\label{c4c6}
c_4(a,b,c) = c_4(E) \quad \text{ and } \quad c_6(a,b,c) = c_6(E) 
\end{equation}
where $E$ has Weierstrass equation $y^2 = x^3 - 27 c_4(E) x - 54 c_6(E)$.
If $j(E) \not= 0,1728$ and $h \in \GL_3(\Kbar)$ is a matrix whose columns 
are covariant columns, respectively contravariant columns, 
of the same degree mod $7$ evaluated at $(a,b,c)$ then 
\[ \sigma(h) h^{-1} \propto \rho ( \sigma(\phi) \phi^{-1} ), \]
respectively
\[ \sigma(h) h^{-1} \propto \rho ( \sigma(\phi) \phi^{-1} )^{-T}, \]
for all $\sigma \in \Gal(\Kbar/K)$. 
\end{Lemma}
\begin{Proof}
Let $\xi_\sigma = \sigma(\phi) \phi^{-1} \in \SL_2(\Z/7\Z)$. 
Recalling that $\rho$ describes the action of $\SL_2(\Z/7\Z)$ 
on $X(7) \subset \PP^2$ we have
\begin{equation}
\label{ddagger}
\sigma ((a,b,c)^T) = \lambda_\sigma \rho( \xi_\sigma)(a,b,c)^T
\end{equation}
for some $\lambda_\sigma \in \Kbar^\times$. Since $c_4$ and $c_6$ are invariants
of degrees $14$ and $21$ we deduce
\[ \sigma (c_4(a,b,c)) = \lambda_\sigma^{14} c_4(a,b,c) 
\quad \text{ and } \quad
 \sigma (c_6(a,b,c)) = \lambda_\sigma^{21} c_6(a,b,c) \]
for all $\sigma \in \Gal(\Kbar/K)$.
Since $c_4(E), c_6(E) \in K$ it follows by~(\ref{c4c6}) and our assumption 
$j(E) \not= 0,1728$ that $\lambda_\sigma^{14} = \lambda_\sigma^{21} = 1$.
Hence $\lambda_\sigma$ is a $7$th root of unity. Now suppose the columns of 
$h$ are obtain by specialising polynomials whose degrees are all congruent
to $r$ mod $7$. Then by~(\ref{ddagger}) and Definition~\ref{defc7} we have
\[ \sigma(h) = h \circ (\lambda_\sigma \rho(\xi_\sigma)) 
= \lambda_\sigma^r \rho(\xi_\sigma) h,\]
respectively
\[ \sigma(h) = h \circ (\lambda_\sigma \rho(\xi_\sigma)) 
= \lambda_\sigma^r \rho(\xi_\sigma)^{-T} h.\]
Hence $\sigma(h) h^{-1} \propto \rho ( \xi_\sigma )$, 
respectively
$\sigma(h) h^{-1} \propto \rho ( \xi_\sigma )^{-T}$, as required.
\end{Proof}

We use Lemmas~\ref{TwistLem} and~\ref{Lem7} to compute equations
for $X_E(7)$ and $X^-_E(7)$. First we classify the covariant and contravariant
columns. It is evident that
\begin{itemize}
\item
The dot product of a covariant column and a contravariant column is 
an invariant. 
\item
The cross product of two covariant columns is a contravariant column.
\item
The cross product of two contravariant columns is a covariant column.
\end{itemize}
We also write $[\vv_1,\vv_2,\vv_3] = (\vv_1 \times \vv_2) \cdot \vv_3$ for
the scalar triple product.
It is easy to solve for the covariant and contravariant columns
of any given degree by linear algebra.
Let $\e$ and $\f$ be the covariant columns of
degrees $9$ and $11$ given by
\small
\begin{align*}
 \e &= \frac{ I_{22} (\nabla F \times \nabla H)
  - c_4 (\nabla F \times \nabla c_4)
  + 12 H^2 (\nabla H \times \nabla c_4) }{14 c_6} \\
 \f & =  \frac{ I_{24} (\nabla F \times \nabla H)
  - (16 F^4 - 104 FH^2) (\nabla F \times \nabla c_4)
  + c_4 (\nabla H \times \nabla c_4) }{14 c_6} 
\end{align*}
\normalsize
where 
\small
\begin{align*}
I_{22} &= 448 F^4 H - 48 F^2 c_4 - 2048 F H^3 \\
I_{24} &= 128 F^6 - 160 F^3 H^2 - 236 F H c_4 - 336 H^4.
\end{align*}
\normalsize
\begin{Lemma} 
\begin{enumerate}
\item The covariant columns of odd degree, respectively even degree, 
form a free $k[F,H,c_4]$-module of rank $3$ generated by $\x$, $\e$, $\f$,
respectively $\nabla F \times \nabla H, \nabla F \times \nabla c_4,
\nabla H \times \nabla c_4$.
\item The contravariant columns of odd degree, respectively even degree, 
form a free $k[F,H,c_4]$-module of rank $3$ 
generated by $\nabla F$, $\nabla H$, 
$\nabla c_4$, respectively $\x \times \e, \x \times \f, \e \times \f$.
\end{enumerate}
\end{Lemma}
\begin{Proof}
By direct calculation we have $[\x,\e,\f] = -c_6$, whereas the 
definition of $c_6$ may be rewritten as 
$[\nabla F,\nabla H,\nabla c_4] = 14 c_6$. 

Let $\vv$ be a covariant column of odd degree. We write $\vv = I_1 \x + I_2 \e + I_3 \f$
where $I_1,I_2,I_3$ are rational functions in $a,b,c$. Taking the
dot product with $\e \times \f$ shows that $[\vv, \e, \f] = I_1 [ \x,\e,\f]$.
But $[\vv,\e,\f]$ is an invariant of odd degree and therefore divisible
by $c_6$. It follows that $I_1$ is an invariant and likewise for $I_2$ and $I_3$.

The other cases are similar.
\end{Proof}

\begin{Theorem} 
\label{Thm7}
Let $E/K$ be an elliptic curve with Weierstrass equation 
$y^2 = x^3 - 27 c_4 x - 54 c_6$ and let $\Delta = (c_4^3- c_6^2)/1728$. 
If $j(E) \not=0,1728$ then $X_E(7) \subset \PP^2$ has equation
$\BF = 0$ where
\begin{align*}
 \BF =  12 x^3 z + 108 x^2 & y^2 + 3 c_4 x^2 z^2  + 72  c_4 x y^2 z 
   - 108 c_4 y^4  - 12 c_6 x y z^2 \\ & + 84 c_6 y^3 z +  c_4^2 x z^3 
    - 15 c_4^2 y^2 z^2 + c_4 c_6 y z^3 + 768 \Delta z^4 ,
\end{align*}
and $X_E^-(7) \subset \PP^2$ has equation
$\BG = 0$ where
\begin{align*}
\BG = 3 x^4  + c_4 & x^3 z - 18 c_4 x^2 y^2 - 3 c_6 x^2 y z + 24 c_6 x y^3 
    + 3 c_4^2 x y^2 z \\ &~- 9 c_4^2 y^4 - c_4 c_6 y^3 z + 168 \Delta x z^3 
    + 1728 \Delta y^2 z^2 + 5 c_4 \Delta z^4.
\end{align*}
\end{Theorem}

\begin{Proof}
The covariant columns $\x$, $\nabla F \times \nabla H$, $H \e$ have
degrees 1, 8, 15, and the contravariant columns 
$\nabla F$, $\x \times \e$, $H^2 \nabla H$ 
have degrees 3, 10, 17. The determinants of the matrices formed from
these columns are
\begin{equation}
\label{det7}
\begin{aligned}
 \det(\x,(\nabla F \times \nabla H) , H \e ) & =  72 H^4 - 4 c_4 F H \\
 \det(\nabla F, (\x \times \e), H^2 \nabla H) & =  72 H^5 - 4 c_4 F H^2.
\end{aligned}
\end{equation} 

The coefficients of the quartic
$\widetilde{F}(x,y,z) 
= F ( x \x + y (\nabla F \times \nabla H) + z H \e )$
are invariants. Using linear algebra to rewrite these invariants 
as polynomials in $F, H, c_4$ and $c_6$ we find
\small
\begin{align*}
  \widetilde{F}&(x,y,z) = F x^4 + 12 H^3 x^3 z 
   + (108 H^3 - 6 c_4 F ) x^2 y^2 
   - 8 c_6 F x y^3 
\\ &~+ 3 c_4 H^3 x^2 z^2 
     + (72 c_4 H^3  + 4128 F^2 H^4 + 48 c_4 F^3 H - 768 F^5 H^2) x y^2 z 
\\ &~+ (-108 c_4 H^3 - 3 c_4^2 F - 11376 F^2 H^4 + 32 c_4 F^3 H 
     + 3392 F^5 H^2 - 256 F^8) y^4 
\\ &~- 12 c_6 H^3 x y z^2 
   + (84 c_6 H^3 - 16 c_6 F^3 H) y^3 z + (c_4^2 H^3 + 688 F H^7 \\ 
&~+ 8 c_4 F^2 H^4  - 128 F^4 H^5) x z^3 
   + (-15 c_4^2 H^3 - 10512 F H^7  - 384 c_4 F^2 H^4  
\\ &~+ 6144 F^4 H^5 + 96 c_4 F^5 H^2 - 768 F^7 H^3) y^2 z^2 
  + (c_4 c_6 H^3 - 8 c_6 F^2 H^4) y z^3 
\\ &~+ (768 H^{10} - 36 c_4 F H^7 - c_4^2 F^2 H^4 + 176 F^3 H^8 
     + 16 c_4 F^4 H^5 - 64 F^6 H^6) z^4.
\end{align*}
\normalsize
Likewise $\widetilde{G}(x,y,z) 
= F ( x \nabla F + y(\x \times \e) + z H^2 \nabla H )$ becomes
\small
\begin{align*}
 \widetilde{G}&(x,y,z) = (3 H^2 + 28 F^3) x^4 
    + (c_4 H^2 + 168 F^2 H^3) x^3 z + (-18 c_4 H^2 \\
&~- 816 F^2 H^3 - 24 c_4 F^3 + 192 F^5 H) x^2 y^2 
    - 3 c_6 H^2 x^2 y z + 24 c_6 H^2 x y^3 
\\ &~+ (222 F H^6 + 24 F^4 H^4) x^2 z^2 
     + (3 c_4^2 H^2 + 3744 F H^6 - 576 F^4 H^4) x y^2 z 
\\ &~+ (-9 c_4^2 H^2 - 5184 F H^6 
    - 240 c_4 F^2 H^3 - 4 c_4^2 F^3 + 2240 F^4 H^4 + 64 c_4 F^5 H 
\\ &~- 256 F^7 H^2) y^4 
 + (-c_4 c_6 H^2 + 8 c_6 F^2 H^3) y^3 z + (168 H^9 + 3 c_4 F H^6  
  \\ &~+ 24 F^3 H^7) x z^3 + (1728 H^9 - 78 c_4 F H^6  
 + 816 F^3 H^7 + 24 c_4 F^4 H^4 \\ &~- 192 F^6 H^5) y^2 z^2 
   + c_6 F H^6 y z^3 + (5 c_4 H^9 + 35 F^2 H^{10} - 4 F^5 H^8) z^4.
\end{align*}
\normalsize

Let $(a:b:c)$ be the $\Kbar$-point on $X(7)$ corresponding to $(E,\phi)$ 
for some choice of symplectic isomorphism 
$\phi : E[7] \isom \mu_7 \times \Z/7\Z$.
By Lemma~\ref{lem:j7} 
we may scale $(a,b,c)$ to satisfy~(\ref{c4c6}). Then 
by Lemmas~\ref{TwistLem} and~\ref{Lem7} 
a formula for $X_E(7)$, respectively $X_E^-(7)$, is given by specialising 
the coefficients of $\widetilde{F}$, respectively $\widetilde{G}$, 
to this choice of $(a,b,c)$. Explicitly we set $F=0$,
divide through by $H^3$, respectively $H^2$, and replace $H^7$ by $\Delta$.
\end{Proof}

\begin{Remark}
\label{tidyup7}
The equations for $X_E(7)$ and $X_E^-(7)$ given in 
Theorems~\ref{MainThm7} and~\ref{Thm7} are related by
$\FF(x,y,z) = \tfrac{1}{4}
\BF( 6 c_4 z - \tfrac{1}{3} y, x, -18 z)$ and 
$\GG(x,y,z)  = \BG( 9 c_4 y + z, 3 x, 108 y)$
where $a = -27 c_4$ and $b= -54c_6$.
\end{Remark}

\subsection{Formulae in the case $n=9$}
\label{sec:invthy9}
We recall that $X(9) =\{ F_1=F_2=0 \} \subset \PP^3$ 
where $F_1 = a^2 b + b^2 c + c^2 a$ and
$F_2 = a b^2 + b c^2 + c a^2 - d^3$.
Let $G \isom \SL_2(\Z/9\Z)$ be the image of $\rho : \SL_2(\Z/9\Z)
\to \SL_4(\Kbar)$. It is generated by 
\[ 
\frac{1}{3} \begin{pmatrix} 
\zz_9     - \zz_9^{8} & \zz_9^{7} - \zz_9^{2} &  \zz_9^{4} - \zz_9^{5} &  \zz_9^{3} - \zz_9^{6} \\
\zz_9^{7} - \zz_9^{2} & \zz_9^{4} - \zz_9^{5} &  \zz_9 - \zz_9^{8} &  \zz_9^{3} - \zz_9^{6} \\
\zz_9^{4} - \zz_9^{5} & \zz_9 - \zz_9^{8} & \zz_9^{7} - \zz_9^{2} &  \zz_9^{3} - \zz_9^{6} \\
\zz_9^{3} - \zz_9^{6} & \zz_9^{3} - \zz_9^6 & \zz_9^{3} - \zz_9^{6} &  0
\end{pmatrix} 
\quad \text{ and } \quad
 \begin{pmatrix}
\zz_9 & 0 & 0 & 0 \\
0 & \zz_9^4 & 0 & 0 \\
0 & 0 & \zz_9^7 & 0 \\
0 & 0 & 0 & \zz_9^6 
\end{pmatrix} \]
Let $\chi: G \to \mu_3$ be the character given by 
$\Diag(\zz_9,\zz_9^4,\zz_9^7,\zz_9^6) \mapsto \zz_9^3$.

\begin{Definition}
An invariant of type $r \in \{0,1,2\}$ and degree $m$ is a 
homogeneous polynomial $I = I(a,b,c,d)$ of degree $m$ such that 
$I \circ g = \chi(g)^r I$ for all $g \in G$. 
\end{Definition}

The type $0$ invariants of smallest degree are 
$D = -(a^3 + b^3 + c^3 - 3 a b c) d$ and
$  I_6 = 2 (a^5 b + b^5 c + c^5 a) + 5 (a^4 c^2 + b^4 a^2 + c^4 b^2) 
       + 20 a b c (a^2 b + b^2 c + c^2 a) - 3 d^6.$
Writing $A$ and $B$ for their matrices of second partial derivatives we find
\[ \det (A + t B) = 81 D^2 - 1620 D I_6 t + 2700 I_{12} t^2 + 90000 D^4 t^4. \]
where $I_{12}$ is a type $0$ invariant of degree $12$.

\begin{Definition}
\label{defc9}
A covariant column, respectively contravariant column, of 
type $r \in \{0,1,2\}$ and
degree $m$ is a column vector $\vv=(v_1,v_2,v_3,v_4)^T$ of 
homogeneous polynomials
of degree $m$ in variables $a,b,c,d$ such that $\vv \circ g = \chi(g)^r g \vv$, 
respectively $\vv \circ g = \chi(g)^r g^{-T} \vv$, for all $g \in G$.
\end{Definition}

We note that $\x_1 = (a,b,c,d)^T$ is a covariant column 
of type $0$ and degree $1$,
whereas if $I$ is an invariant of type $r$ and degree $m$ then 
$\nabla I = (\tfrac{\partial I}{\partial a},\tfrac{\partial I}{\partial b},
\tfrac{\partial I}{\partial c},\tfrac{\partial I}{\partial d})^T$ 
is a contravariant column of type $r$ and degree $m-1$. To construct some 
further examples we 
define $4 \times 4$ alternating matrices
\[ \Lambda_i = \begin{pmatrix} 
0 & \qq_i(a,b,c) d^2 & -\qq_i(c,a,b) d^2 & \rr_i(a,b,c) \\
-\qq_i(a,b,c) d^2 & 0 & \qq_i(b,c,a) d^2 & \rr_i(b,c,a) \\
\qq_i(c,a,b) d^2 & -\qq_i(b,c,a) d^2 & 0 & \rr_i(c,a,b) \\
-\rr_i(a,b,c) & -\rr_i(b,c,a) & -\rr_i(c,a,b) & 0
\end{pmatrix} \]
for $i=0,1,2$ where
\begin{align*}
\qq_0(a,b,c) &= 3 (a^2 + 2 b c) &
\rr_0(a,b,c) &= 2 (a^3 b + c^3 b) - 3 a^2 c^2 - b^4 \\
\qq_1(a,b,c) &= 3 (c^2 + 2 a b) &
\rr_1(a,b,c) &= 2 (b^3 c + a^3 c) - 3 b^2 a^2 - c^4 \\
\qq_2(a,b,c) &= 3 (b^2 + 2 a c) &
\rr_2(a,b,c) &= 2 (c^3 a + b^3 a) - 3 c^2 b^2 - a^4.
\end{align*}
For $M$ a $4 \times 4$ alternating matrix and $\vv$ a column vector
we define a column vector $M \star \vv$ by the rule
$(M \star \vv)_l = \sum M_{ij} v_k$ where the sum is over all $(i,j,k)$ 
for which $(i,j,k,l)$ is an even permutation of $(1,2,3,4)$.
We define covariant and contravariant columns
\begin{align*}
  \x_7 &=  \Lambda_0 \nabla D &  \e_9 &= \Lambda_1 \nabla I_6 & \uu_5 &= \Lambda_1 \star \x_1 \\
  \e_7 &= \Lambda_1 \nabla D  & \f_7 &= \Lambda_2 \nabla D & \vv_{11} &= \Lambda_2 \star \x_7.
\end{align*}
Then $c_4 = \x_7 \cdot \uu_5$ is an invariant of type $1$ and degree $12$,
and $c_6 = \tfrac{1}{2} (\x_7 \cdot \nabla I_{12})$ is an 
invariant of type $0$
and degree $18$. Routine Gr\"obner basis calculations show that 
$c_4,c_6 \notin (F_1,F_2)$ yet
\[ c_4^3 - c_6^2 \equiv 1728 D^9 \mod{(F_1,F_2)}. \]
\begin{Lemma}
\label{lem:j9}
The $j$-invariant $X(9) \to \PP^1$ is given by $j = c_4^3/D^9$.
\end{Lemma}
\begin{Proof}
Both $j$ and $j_0 = c_4^3/D^9$ define maps $X(9) \to \PP^1$ 
that quotient out by the action of $G/\{ \pm I\} \isom \PSL_2(\Z/9\Z)$. 
So they can differ by at most a M\"obius map. We recall that $j$
is ramified above $0$, $1728$ and $\infty$ with ramification
indices $3$, $2$ and $9$. Since
\begin{align*}
\# \{ F_1 = F_2 = c_4 = 0 \} &\le 9 \deg (c_4) = \tfrac{1}{3} |G/\{ \pm I\}| \\
\# \{ F_1 = F_2 = c_6 = 0 \} &\le 9 \deg (c_6) = \tfrac{1}{2} |G/\{ \pm I\}| \\
\# \{ F_1 = F_2 = D = 0 \} &\le 9 \deg (D) = \tfrac{1}{9} |G/\{ \pm I\}| 
\end{align*}
and $j_0 - 1728 = c_6^2/D^9$ it follows that $j=j_0$ as required.
\end{Proof}
\begin{Lemma}
\label{Lem9}
Let $E/K$ be an elliptic curve and $\phi : E[9] \isom \mu_9 \times \Z/9\Z$
a symplectic isomorphism over $\Kbar$.
Let $(a:b:c:d)$ be the
corresponding $\Kbar$-point on $X(9) \subset \PP^3$ with co-ordinates
$(a,b,c,d)$ scaled so that 
\begin{equation}
\label{c4c6-nine}
c_4(a,b,c,d) = c_4(E) \quad \text{ and } \quad
c_6(a,b,c,d) = c_6(E) 
\end{equation}
where $E$ has Weierstrass equation $y^2 = x^3 - 27 c_4(E) x - 54 c_6(E)$.
If $j(E) \not= 0,1728$ and $h \in \GL_4(\Kbar)$ is a matrix whose columns 
are covariant columns, respectively contravariant columns, 
of types $r_1, \ldots, r_4$ and degrees $m_1, \ldots, m_4$ 
with $m_i + 6 r_i$ constant mod $18$, evaluated at $(a,b,c,d)$ then 
\[ \sigma(h) h^{-1} \propto \rho ( \sigma(\phi) \phi^{-1} ), \]
respectively
\[ \sigma(h) h^{-1} \propto \rho ( \sigma(\phi) \phi^{-1} )^{-T}, \]
for all $\sigma \in \Gal(\Kbar/K)$. 
\end{Lemma}
\begin{Proof}
Let $\xi_\sigma = \sigma(\phi) \phi^{-1} \in \SL_2(\Z/9\Z)$. 
Recalling that $\rho$ describes the action of $\SL_2(\Z/9\Z)$ 
on $X(9) \subset \PP^3$
we have
\begin{equation}
\label{ddagger-nine}
\sigma ((a,b,c,d)^T) = \lambda_\sigma \rho( \xi_\sigma) (a,b,c,d)^T
\end{equation}
for some $\lambda_\sigma \in \Kbar^\times$. Since $c_4$ and $c_6$ are invariants
of types $1$ and $0$ and degrees $12$ and $18$ we deduce
\[ \sigma (c_4(a,b,c,d)) = \lambda_\sigma^{12} \chi_\sigma c_4(a,b,c,d) 
\quad \text{ and } \quad
 \sigma (c_6(a,b,c,d)) = \lambda_\sigma^{18} c_6(a,b,c,d) \]
where $\chi_\sigma = \chi(\rho(\xi_\sigma))$. 
Since $c_4(E), c_6(E) \in K$ it follows by~(\ref{c4c6-nine}) and our assumption 
$j(E) \not= 0,1728$ that $\lambda_\sigma^{12} \chi_\sigma = \lambda_\sigma^{18}= 1$.
Hence $\lambda_\sigma$ is an $18$th root of unity and $\chi_\sigma = \lambda_\sigma^6$. 
Now suppose the columns of $h$ are obtained by specialising covariant columns,
respectively contravariant columns, of types $r_1, \ldots, r_4$ and 
degrees $m_1, \ldots, m_4$ with $m_i + 6 r_i \equiv r \pmod{18}$.
Then by~(\ref{ddagger-nine}) and Definition~\ref{defc9} we have
\[ \sigma(h) = h \circ (\lambda_\sigma \rho(\xi_\sigma)) 
= \lambda_\sigma^r \rho(\xi_\sigma) h,\]
respectively
\[ \sigma(h) = h \circ (\lambda_\sigma \rho(\xi_\sigma)) 
= \lambda_\sigma^r \rho(\xi_\sigma)^{-T} h.\]
\end{Proof}

We use Lemmas~\ref{TwistLem} and~\ref{Lem9} to compute equations
for $X_E(9)$ and $X^-_E(9)$. First we construct some more
covariant and contravariant columns. 
Let $\vv_9$ be the contravariant column of
type $2$ and degree $9$ given by
\[\vv_9 = (f_1 d, f_2 d, f_3 d, -a f_1 - b f_2 - c f_3)^T \]
where
\begin{align*}
  f_1(&a,b,c,d) = a^8 + 35 a^6 b c + 42 a^5 b^3 + 7 a^5 c^3 + 105 a^4 b^2 c^2 + 28 a^4 b d^3 
\\ & + 231 a^3 b^4 c - 196 a^3 b c^4 + 14 a^3 c^2 d^3 + 14 a^2 b^6 - 70 a^2 b^3 c^3 - 
      84 a^2 b^2 c d^3 \\ & - 21 a^2 c^6 + 28 a b^4 d^3 - 105 a b^2 c^5 - 14 a b c^3 d^3 - 
      27 a c d^6 + 19 b^7 c \\ & - 35 b^4 c^4 + 14 b^3 c^2 d^3 + 27 b^2 d^6 - 27 b c^7 + 
      14 c^5 d^3,
\end{align*}
$f_2(a,b,c,d) = f_1(b,c,a,d)$ and $f_3(a,b,c,d) = f_1(c,a,b,d)$. 
We further put
\begin{align*}
  \x_{13} &= \Lambda_1 \vv_9 &
  \x_{15} &= \Lambda_1 \vv_{11} &
  \f_{13} &= \Lambda_0 \vv_9. 
\end{align*}

In the following lemma the entries of the covariant and contravariant 
columns are viewed as elements of the co-ordinate ring 
$K[a,b,c,d]/(F_1,F_2)$.
The lemma is included to show that we have been systematic,
rather than because it is needed in what follows. We therefore 
omit the proof.
\begin{Lemma}
The covariant columns, respectively contravariant columns, of type
$r \in \{0,1,2\}$, mod $(F_1,F_2)$, form a free $K[D,c_6]$-module of
rank $4$ with basis as indicated in the following table. 
\[ \begin{array}{ccccccc}
\multicolumn{3}{c}{\text{Covariants}} & & \multicolumn{3}{c}{\text{Contravariants}} \\
\text{Type 0} & \text{Type 1} & \text{Type 2} & \qquad & 
\text{Type 0} & \text{Type 1} & \text{Type 2} \\
\x_1 & \e_7 & \f_7 & & \nabla D & \uu_5 & \vv_9 \\
\x_7 & \e_9 & \f_{13} & & \nabla I_6 & \nabla c_4 & \vv_{11} \\
\x_{13} & c_4 \x_1 & c_4 \e_7 & & \nabla I_{12} & c_4 \nabla D & c_4 \uu_5 \\
\x_{15} & c_4 \x_7 & c_4 \e_9 & & \nabla c_6 & c_4 \nabla I_6 & c_4 \nabla c_4 \\
\end{array} \]
\end{Lemma}

\begin{Theorem} 
\label{Thm9}
Let $E/K$ be an elliptic curve with Weierstrass equation 
$y^2 = x^3 - 27 c_4 x - 54 c_6$.
If $j(E) \not=0,1728$ then $X_E(9) \subset \PP^3$ has equations
$\BF_1 = \BF_2 = 0$ where  
\begin{align*}
  \BF_1 &= 24 x^2 t - 96 x y z + 64 y^3 + 24 c_4 y^2 t + 48 c_4 y z^2 
     - 48 c_6 y z t + 12 c_4^2 y t^2 \\ & \qquad - 8 c_6 z^3 + 12 c_4^2 z^2 t 
     - 6 c_4 c_6 z t^2 - (c_4^3 - 2 c_6^2) t^3, \\
  \BF_2 &= 16 x^2 z - 64 x y^2 - 16 c_4 x y t - 16 c_4 x z^2 + 16 c_6 x z t 
     - 4 c_4^2 x t^2 + 80 c_4 y^2 z \\ & \qquad - 32 c_6 y^2 t - 32 c_6 y z^2 
     + 32 c_4^2 y z t - 8 c_4 c_6 y t^2 + 8 c_4^2 z^3 - 12 c_4 c_6 z^2 t 
     \\ & \qquad + (2 c_4^3 + 4 c_6^2) z t^2 - c_4^2 c_6 t^3.
\end{align*}
and $X^-_E(9) \subset \PP^3$ has equations
$\BG_1 = \BG_2 = 0$ where  
\begin{align*}
  \BG_1 &= -72 x^2 y - 144 x^2 z + 24 c_4 x y t - 48 c_4 x z t - 8 c_6 x t^2 
     - c_4 y^3 + 18 c_4 y^2 z \\ & \qquad + 3 c_6 y^2 t 
     + 180 c_4 y z^2 + 12 c_6 y z t 
     - 3 c_4^2 y t^2 + 792 c_4 z^3 + 12 c_6 z^2 t  \\ & \qquad
     + 2 c_4^2 z t^2 + c_4 c_6 t^3, \\
  \BG_2 &= -864 x^3 + 216 c_4 x^2 t + 2592 c_4 x y z + 216 c_6 x y t 
     + 15552 c_4 x z^2 \\ & \qquad + 432 c_6 x z t - 72 c_4^2 x t^2 - 9 c_6 y^3 
     + 162 c_6 y^2 z + 27 c_4^2 y^2 t + 4212 c_6 y z^2 \\ & \qquad 
     + 108 c_4^2 y z t - 27 c_4 c_6 y t^2 + 26136 c_6 z^3 
     + 972 c_4^2 z^2 t + 18 c_4 c_6 z t^2   \\ 
     & \qquad + (5 c_4^3 + 4 c_6^2) t^3.
\end{align*}

\end{Theorem}
\begin{Proof}
The covariant columns $\x_1,\f_7,D \e_9,D \x_{15}$ have types $0,2,1,0$ and
degrees $1,7,13,19$. The contravariant columns 
$\uu_5,D^2 \nabla D,\nabla I_{12}, D^2 \vv_9$ have types $1,0,0,2$ 
and degrees $5,11,11,17$.
The determinants of the matrices formed from these columns
satisfy
\begin{align*}
 \det(\x_1,\f_7,D \e_9,D \x_{15}) & \equiv 3456 D^{10} \mod{(F_1,F_2)} \\
 \det(\uu_5,D^2 \nabla D,\nabla I_{12}, D^2 \vv_9) & \equiv 1152 D^{11} 
\mod{(F_1,F_2)}. 
\end{align*}

For certain $2 \times 2$ matrices $\alpha = (\alpha_{ij})$ and 
$\beta = (\beta_{ij})$,
specified in Remark~\ref{2by2} below, we put
\begin{align*}
 {\BF}_i(x,y,z,t) & = (\alpha_{i1} F_1 + \alpha_{i2} F_2) 
( x \x_1 + y \f_7 + z D \e_9 + t D \x_{15}) \\
 {\BG}_i(x,y,z,t) & = (\beta_{i1} F_1 + \beta_{i2} F_2) 
( x \uu_5 + y D^2 \nabla D + z \nabla I_{12} + t D^2 \vv_9) 
\end{align*}
for $i=1,2$. Using the Gr\"obner basis machinery in Magma to 
write each coefficient mod $(F_1,F_2)$ as a polynomial
in $c_4$ and $c_6$ we obtained the equations in the statement 
of the theorem.
By Lemmas~\ref{TwistLem}, \ref{lem:j9} and~\ref{Lem9} these are equations 
for $X_E(9)$ and $X_E^-(9)$.
\end{Proof}

\begin{Remark}
\label{2by2}
(i) The $2 \times 2$ matrices used in the proof of Theorem~\ref{Thm9} were
\[ \alpha = \frac{1}{9D^6} 
\begin{pmatrix}  3 B & 9 A \\ -(108 A^3+B^3) & 9 A B^2 
\end{pmatrix} \,\,\, \text{and} \,\,\,
 \beta = \frac{1}{D^6} 
\begin{pmatrix} 6 A & -B \\ 162 A B^2 & 9 (108 A^3+B^3) 
\end{pmatrix} \]
where $A = d^3$ and $B = a^3 + b^3 + c^3 + 6 abc$. \\
(ii) The equations for $X_E(9)$ and $X_E^-(9)$ in 
Theorems~\ref{MainThm9} and~\ref{Thm9} are related by
\begin{align*}
 \FF_i(x,y,z,t) &\propto \BF_i(x - 27 c_4 z - 162 c_6 t,
  9 y + 81 c_4 t,-54 z,324 t)  \\
 \GG_i(x,y,z,t) &\propto \BG_i(-x + 9 c_4 t, 36 y + 36 z,6 z, 108 t) 
\end{align*}
for $i=1,2$ where $a = -27c_4$ and $b = -54 c_6$.
\end{Remark}

\subsection{Formulae in the case $n=11$}
\label{sec:formulae11}

We recall that $X(11)$ is the singular locus of the Hessian of 
the cubic threefold $\{ F = 0 \} \subset \PP^4$ where
\[ F = a^2 b + b^2 c + c^2 d + d^2 e + e^2 a. \]
Let $G \isom \PSL_2(\Z/11\Z)$ be the image of $\rho : \SL_2(\Z/11\Z)
\to \GL_5(\Kbar)$. It is generated by
\[ \frac{1}{g_{11}} \begin{pmatrix} 
\zz_{11} - \zz_{11}^{-1} &  \zz_{11}^{3} - \zz_{11}^{-3} &  \zz_{11}^{9} - \zz_{11}^{-9} &  \zz_{11}^{5} - \zz_{11}^{-5} &
\zz_{11}^{4} - \zz_{11}^{-4} \\
\zz_{11}^{3} - \zz_{11}^{-3} &  \zz_{11}^{9} - \zz_{11}^{-9} &  \zz_{11}^{5} - \zz_{11}^{-5} &  \zz_{11}^{4} - \zz_{11}^{-4} &
\zz_{11} - \zz_{11}^{-1} \\
\zz_{11}^{9} - \zz_{11}^{-9} &  \zz_{11}^{5} - \zz_{11}^{-5} &  \zz_{11}^{4} - \zz_{11}^{-4} &  \zz_{11} - \zz_{11}^{-1} &
\zz_{11}^{3} - \zz_{11}^{-3} \\
\zz_{11}^{5} - \zz_{11}^{-5} &  \zz_{11}^{4} - \zz_{11}^{-4} &  \zz_{11} - \zz_{11}^{-1} &  \zz_{11}^3 - \zz_{11}^{-3} &
\zz_{11}^{9} - \zz_{11}^{-9} \\
\zz_{11}^{4} - \zz_{11}^{-4} &  \zz_{11} - \zz_{11}^{-1} &  \zz_{11}^3 - \zz_{11}^{-3} &  \zz_{11}^9 - \zz_{11}^{-9} &
\zz_{11}^{5} - \zz_{11}^{-5} 
\end{pmatrix} \]
and $\Diag(\zz_{11},\zz_{11}^9,\zz_{11}^4,\zz_{11}^3,\zz_{11}^5)$ where 
$g_{11} = 1 + 2 (\zz_{11} + \zz_{11}^3 + \zz_{11}^9 + \zz_{11}^5 + \zz_{11}^4)$. 
We define the invariants, covariant columns and contravariant columns
exactly as in Section~\ref{sec:formulae7}.
Let $\sum$ denote a sum over all cyclic permutations, so that for example
$F = \sum a^2 b$. Other examples of invariants of small degree include
\begin{align*}
  H &= 3 a b c d e + \textstyle\sum(a^3 c^2 - a^3 d e), \\
  I_7 &= \textstyle\sum(a^6 e + 3 a^5 d^2 - 15 a^4 b c e + 5 a^3 b^3 d + 15 a^3 b c d^2), \\
  I_8 &= \textstyle\sum(a^7 c - 7 a^4 b d^3 - 7 a^4 d e^3 
                    + 7 a^3 b^2 c^3 + 21 a^3 c^2 d^2 e). 
\end{align*}
Writing $A$ and $B$ for the matrices of second partial derivatives
of $F$ and $H$ we find
\[ \det (A + t B) = 32 H - 32 I_{7} t - 24 I_{9} t^2 - 8 c_4 t^3 + \ldots \]
where $I_{9}$ and $c_4$ are invariants of degrees 9 and 11.
We will not need a complete set of generators for the ring of invariants,
but note that this is given in~\cite{AdlerINV11} and may also 
be computed using Magma.
Let $\I$ be the homogeneous ideal of $X(11)$, i.e. the ideal generated 
by the $4 \times 4$ minors of the Hessian matrix of $F$. The degree $19$ 
polynomial
\begin{align*}
\widetilde{c}_6 &= a^9 b^{10} - 509 b^{18} d - 14107 b^{14} d^4 e 
   + 510 b^9 c^{10} + 
    42326 b^7 d^{12} + 20669 b^3 d^{15} e \\ & - 14107 b^2 d^2 e^{15} - 
    277419 b c^2 d^{10} e^6 - 248909 b c d^{16} e - 
    209926 b c d^5 e^{12} \\ & + 762409 b d^{11} e^7 + b e^{18} - 
    1018 c^{18} e - 14107 c^{16} d e^2 - 586835 c^{12} d^3 e^4 \\ & + 
    197780 c^{10} d^4 e^5 + 1019 c^9 d^{10} - 787130 c^8 d^5 e^6 + 
    15634 c^7 d^{11} e + 42326 c^7 e^{12} \\ & + 2007576 c^6 d^6 e^7 + 
    247382 c^5 d^{12} e^2 - 528424 c^5 d e^{13} - 616653 c^4 d^7 e^8 \\ &
    + 376744 c^3 d^{13} e^3 + 1067732 c^3 d^2 e^{14} - 
    225004 c^2 d^8 e^9 + 463659 c d^{14} e^4 \\ & - 582142 c d^3 e^{15} +
    70511 d^9 e^{10}
\end{align*}
is not an invariant but satisfies
\[ {\widetilde{c}_6}^{\,\,2} \equiv a b c d e (c_4^3 - 1728 F^{11}) 
\pmod{\I}. \]

\begin{Lemma}
\label{lem:j11}
The $j$-invariant $X(11) \to \PP^1$ is given by $j = c_4^3/F^{11}$.
\end{Lemma}
\begin{Proof}
Both $j$ and $j_0 = c_4^3/F^{11}$ define maps $X(11) \to \PP^1$ 
that quotient out by the action of $G \isom \PSL_2(\Z/11\Z)$. 
So they can differ by at most a M\"obius map. We recall that $j$
is ramified above $0$, $1728$ and $\infty$ with ramification
indices $3$, $2$ and $11$. It is shown in \cite[Corollary 23.28]{AR}
that $X(11) \subset \PP^4$ has degree 20. Since
\begin{align*}
\# X(11) \cap \{ c_4 = 0 \} &\le 20 \deg (c_4) = \tfrac{1}{3} |G|\\
\# X(11) \cap \{ \widetilde{c}_6 = 0 \} &\le 20 \deg (\widetilde{c}_6) < |G|\\
\# X(11) \cap \{ F = 0 \} &\le 20 \deg (F) = \tfrac{1}{11} |G| 
\end{align*}
and $j_0 - 1728 = \widetilde{c}_6^{\,\,2}/ ((abcde) F^{11})$ it follows 
that $j=j_0$ as required.
\end{Proof}

\begin{Lemma}
\label{Lem11}
Let $E/K$ be an elliptic curve and $\phi : E[11] \isom \mu_{11} 
\times \Z/11\Z$ a symplectic isomorphism over $\Kbar$. Let $(a:b:c:d:e)$ 
be the corresponding $\Kbar$-point on $X(11) \subset \PP^4$ with co-ordinates
$(a,b,c,d,e)$ scaled so that 
\begin{equation}
\label{c4:11}
c_4(a,b,c,d,e) = c_4(E) 
\end{equation}
where $E$ has Weierstrass equation $y^2 = x^3 - 27 c_4(E) x - 54 c_6(E)$.
If $j(E) \not= 0$ and $h \in \GL_5(\Kbar)$ is a matrix whose columns 
are covariant columns, respectively contravariant columns, 
of the same degree mod $11$ evaluated at $(a,b,c,d,e)$ then 
\[ \sigma(h) h^{-1} \propto \rho ( \sigma(\phi) \phi^{-1} ), \]
respectively
\[ \sigma(h) h^{-1} \propto \rho ( \sigma(\phi) \phi^{-1} )^{-T}, \]
for all $\sigma \in \Gal(\Kbar/K)$. 
\end{Lemma}
\begin{Proof}
The proof is similar to that of Lemma~\ref{Lem7}. Recall that
$c_4$ is a homogeneous polynomial of degree $11$ and so~(\ref{c4:11})
determines the scaling of $(a,b,c,d,e)$ up to an $11$th root of unity.
\end{Proof}

We use Lemmas~\ref{TwistLem} and~\ref{Lem11} to compute equations
for $X_E(11)$ and $X_E^-(11)$. First we construct some covariant columns.
Let $\x_1 = (a,b,c,d,e)^T$. If $\gamma \in \SL_2(\Z/11\Z)$ is 
diagonal then $\rho (\gamma)$ 
cyclically permutes the co-ordinates $a,b,c,d,e$. A covariant column
is therefore uniquely determined by its first entry. 
By averaging over the group
we found covariant columns $\x_4, \x_5, \x_9$ with first entries
\small
\begin{align*}
  f_4 &= 2 a^2 e^2 + 4 a b^2 c - 4 a c^2 d + 4 b c e^2 + d^4, \\
  f_5 &= -5 a^3 c e + 5 a^2 b^2 d + 5 a^2 c d^2 + 5 a b c^2 e - 10 a b d e^2 
      + b^5 - 5 b^3 c d + 5 b d^3 e + 5 c^2 e^3, \\ 
f_9 &= 
    -14 a^6 b d e - 8 a^5 b d^3 + 9 a^5 c^2 e^2 + 2 a^5 d e^3 + 8 a^4 b^4 e + 
        5 a^4 b^2 c^3 + 63 a^4 b^2 c d e \\ & + 6 a^4 c^4 d - 18 a^4 c^2 d^2 e + 
        8 a^4 d^3 e^2 + 31 a^3 b^4 d^2 - 21 a^3 b^3 e^3 + 47 a^3 b^2 c d^3 + 
        35 a^3 b c^3 e^2 \\ & + 14 a^3 b c d e^3 - 12 a^3 c^2 d^4 + 10 a^3 d^5 e + 
        3 a^2 b^5 c e - 26 a^2 b^3 c^4 - 42 a^2 b^3 c^2 d e \\ & - 75 a^2 b^3 d^2 e^2
        + 3 a^2 b^2 e^5 + 18 a^2 b c^5 d - 30 a^2 b c^3 d^2 e - 
        36 a^2 b c d^3 e^2 + 2 a^2 c^3 e^4 \\ & - 9 a^2 c d e^5 + a^2 d^7 - 2 a b^7 d
        - 6 a b^5 c d^2 + 50 a b^4 c e^3 - 7 a b^3 c^2 d^3 - 6 a b^3 d^4 e \\ & - 
        54 a b^2 c^4 e^2 - 3 a b^2 c^2 d e^3 - 9 a b^2 d^2 e^4 - 29 a b c^3 d^4 
        + 21 a b c d^5 e + a b e^7 + 9 a c^5 d e^2 \\ & + 25 a c^3 d^2 e^3 - 
        7 a c d^3 e^4 - 10 b^6 c^2 e - 2 b^6 d e^2 + 4 b^4 c^5 + 40 b^4 c^3 d e 
        - 6 b^4 c d^2 e^2 \\ & + 13 b^3 c e^5 - 3 b^3 d^6 - 15 b^2 c^4 d^2 e - 
        54 b^2 c^2 d^3 e^2 + 31 b^2 d^4 e^3 - 11 b c^4 e^4 + 3 b c^2 d e^5 \\ & - 
        2 b c d^7 - 7 b d^2 e^6 - c^7 d^2 + 5 c^5 d^3 e - 9 c^3 d^4 e^2 + 
        8 c d^5 e^3 - e^9.
\end{align*}
\normalsize
We temporarily write
$a_1, \ldots ,a_5$ for $a,b,c,d,e$ and let $\Xi$ be the $5 \times 5$
alternating matrix with entries 
\[ \Xi_{ij} = \frac{\partial F }{\partial a_r} \,
\frac{\partial I_7}{\partial a_s} - \frac{\partial F }{\partial a_s} 
\, \frac{\partial I_7}{\partial a_r} \]
where $r \equiv (i-j-2)^3 + i + 3 \pmod{5}$ and 
$s \equiv (j-i-2)^3 + j + 3 \pmod{5}$. Then $\x_{14} = \Xi \nabla I_7$
is a covariant column of degree $14$.

\begin{Theorem} 
\label{Thm11}
Let $E/K$ be an elliptic curve with Weierstrass equation 
$y^2 = x^3 - 27 c_4 x - 54 c_6$ and let $\Delta = (c_4^3- c_6^2)/1728$. 
If $j(E) \not=0,1728$ then $X_E(11) \subset \PP^4$ is the singular
locus of the Hessian of
\begin{align*}
 \BF &= v^3 + 3 v^2 w + c_4 v^2 y + 3 v w^2 + 2 c_4 v w y - c_4 \Delta v x^2 
    + 48 \Delta v x y \\ & + 9 w^3 + 5 c_4 w^2 y - c_4^2 w^2 z + c_4^2 w y^2 
    - 576 \Delta w y z + 72 c_4 \Delta w z^2 \\ &  - 4 \Delta^2 x^3 
    - 72 \Delta^2 x^2 z + 4 c_4 \Delta x y^2 + 2 c_4^2 \Delta x y z 
    - (c_4^3 \Delta - 1728 \Delta^2) x z^2 \\ & + 64 \Delta y^3 
    - 72 c_4 \Delta y^2 z + 12 c_4^2 \Delta y z^2 
    + (c_4^3 \Delta - 3456 \Delta^2) z^3,
\end{align*}
and $X_E^-(11) \subset \PP^4$ is the singular locus of the Hessian of
\begin{align*}
  \BG &= 5 v^3 - c_4 v^2 x - 60 v^2 y + 28 c_4 v^2 z - 2 c_4 \Delta v w^2 
    - 48 \Delta v w x \\ & - 240 \Delta v w z - 16 c_4 v x y + 1680 v y^2 
    - 872 c_4 v y z + 121 c_4^2 v z^2 \\ & 
     + 8 \Delta^2 w^3 + 44 c_4 \Delta w^2 y 
    - 11 c_4^2 \Delta w^2 z + c_4 \Delta w x^2 + 336 \Delta w x y \\ & 
    - 122 c_4 \Delta w x z + 25 c_4^2 w y^2 - 14160 \Delta w y z 
    + 817 c_4 \Delta w z^2 - 20 \Delta x^3 \\ & 
      + 5 c_4^2 x^2 y - 1884 \Delta x^2 z
    - 364 c_4 x y^2 + 160 c_4^2 x y z - 34764 \Delta x z^2 \\ & + 19840 y^3 
    - 10268 c_4 y^2 z + 1643 c_4^2 y z^2 - 129220 \Delta z^3.
\end{align*} 

\end{Theorem}
\begin{Proof}
The covariant columns $\x_1,\x_4,\x_5,\x_9,\x_{14}$ have degrees
$1,4,5,9,14$ and the contravariant columns 
$\nabla F,\nabla I_7, \nabla I_8, \nabla I_9, \nabla c_4$
have degrees $2,6,7,8,10$. The determinants of the matrices formed
from these columns satisfy
\begin{equation}
\label{det11}
\begin{aligned}
\det(\x_1,\x_4,\x_5,\x_9,\x_{14}) &=  c_4^3 - 1728 F^{11} \pmod{\I} \\
\det(\nabla F,\nabla I_7, \nabla I_8, \nabla I_9, \nabla c_4) 
&=  55(c_4^3 - 1728 F^{11}) \pmod{\I}.
\end{aligned}
\end{equation}

The coefficients of the cubic
$\widetilde{F}(v,w,x,y,z) 
     = F( v \x_1 + w \x_4 + x \x_5 + y \x_9 + z \x_{14})$
are invariants. Using the Gr\"obner basis machinery in Magma to rewrite
the coefficients mod $\I$ as polynomials in $c_4$ and $F$ we find
\begin{align*}
\widetilde{F} &=  F v^3 + 3 F^2 v^2 w + c_4 v^2 y + 3 F^3 v w^2 + 2 F c_4 v w y 
    - c_4 v x^2 + 48 F^5 v x y \\ & + 9 F^4 w^3 + 5 F^2 c_4 w^2 y - c_4^2 w^2 z 
    + c_4^2 w y^2 - 576 F^9 w y z + 72 F^7 c_4 w z^2 \\ & - 4 F^5 x^3 
    - 72 F^8 x^2 z + 4 F^4 c_4 x y^2 + 2 F^2 c_4^2 x y z 
    - (c_4^3 - 1728 F^{11}) x z^2 \\ & + 64 F^9 y^3 - 72 F^7 c_4 y^2 z 
    + 12 F^5 c_4^2 y z^2 + (F^3 c_4^3 - 3456 F^{14}) z^3.
\end{align*}
Likewise $\widetilde{G}(v,w,x,y,z) 
     = F( v \nabla F + w \nabla I_7 + x \nabla I_8 
+ y \nabla I_9 + z \nabla c_4)$ becomes
\begin{align*}
\widetilde{G} &= 5 F^2 v^3 - c_4 v^2 x - 60 F^4 v^2 y + 28 F c_4 v^2 z 
    - 2 F c_4 v w^2 - 48 F^5 v w x \\ & - 240 F^6 v w z - 16 F^2 c_4 v x y 
    + 1680 F^6 v y^2 - 872 F^3 c_4 v y z + 121 c_4^2 v z^2 \\ & + 8 F^6 w^3 
    + 44 F^3 c_4 w^2 y - 11 c_4^2 w^2 z + F^3 c_4 w x^2 + 336 F^7 w x y \\ & 
    - 122 F^4 c_4 w x z + 25 c_4^2 w y^2 - 14160 F^8 w y z + 817 F^5 c_4 w z^2 
    - 20 F^7 x^3 \\ & + 5 c_4^2 x^2 y - 1884 F^8 x^2 z - 364 F^4 c_4 x y^2 
    + 160 F c_4^2 x y z - 34764 F^9 x z^2 \\ & + 19840 F^8 y^3 
    - 10268 F^5 c_4 y^2 z + 1643 F^2 c_4^2 y z^2 - 129220 F^{10} z^3.
\end{align*}

Let $(a:b:c:d:e)$ be the $\Kbar$-point on $X(11)$ 
corresponding to $(E,\phi)$ for some choice of symplectic isomorphism 
$\phi : E[11] \isom \mu_{11} \times \Z/11\Z$. By Lemma~\ref{lem:j11} 
we may scale $(a,b,c,d,e)$ to satisfy~(\ref{c4:11}) and 
$F(a,b,c,d,e)^{11} = \Delta$. 
Moreover the determinants~(\ref{det11}) are non-zero by our assumption
$j(E) \not= 1728$. Then by Lemmas~\ref{TwistLem} and~\ref{Lem11} 
cubics describing $X_E(11)$ and $X_E^-(11)$ are obtained by putting
\begin{equation}
\label{FG11}
\begin{aligned}
 \BF(v,w,x,y,z) & = \frac{1}{F^{4}} \widetilde{F}(F v,w,F^7 x,F^2 y,F^4 z) \\
 \BG(v,w,x,y,z) & =\frac{1}{F^{8}}\widetilde{G}(F^2 v,F^8 w,F^4 x,y,F^3 z) 
\end{aligned}
\end{equation}
and replacing $F^{11}$ by $\Delta$.
\end{Proof}

\begin{Remark}
The cubic forms describing $X_E(11)$ and $X_E^-(11)$ in 
Theorems~\ref{MainThm11} and~\ref{Thm11} are related by
\begin{align*}
\FF(v,w,x,y,z) &= \frac{1}{2^3 c_6^3} 
\BF(-v',w',   -864 x, -36 c_4 x - 108 c_6 z, 72 y) \\
\GG(v,w,x,y,z) &= \frac{1}{2^5 3^6 (55 c_6)^3} 
\BG(v'',-427680 y ,   x'' , -y'',-z'') 
\end{align*}
where $a = -27c_4$, $b = -54 c_6$ and
\small
\begin{align*}
v' &= c_6 v + 2 c_6 w - 6 c_4^2 x + 3 c_4^2 y - 9 c_4 c_6 z &
w' &= c_6 v + 6 c_4^2 x + 3 c_4^2 y + 9 c_4 c_6 z \\
v'' &= 44 (2 c_4 v - 6 c_6 w + 33 c_6 x + 135 c_4^2 y + 810 c_4 c_6 z) &
x'' &= 60 (5 v + 729 c_4 y + 2187 c_6 z) \\
y'' &= 11 (c_4 v - 3 c_6 w - 6 c_6 x) &
z'' &= 60 (v + 27 c_4 y + 81 c_6 z).
\end{align*}
\normalsize
\end{Remark}

\section{Diagonal twists}
\label{sec:diag}

We give an alternative construction of $X_E(n)$
and $X_E^-(n)$ in the case $E$ is an elliptic curve whose $n$-torsion
contains a copy of the Galois module $\mu_n$.

Let $C \to D$ be an isogeny of elliptic curves with kernel 
a labelled copy of $\mu_n$. Then the dual isogeny has kernel a
labelled copy of $\Z/n\Z$. The pairs of such curves are parametrised
by the modular curve $Y_1(n)$. In the cases $n=7,9$ we choose a coordinate
$\la$ on $X_1(n) \isom \PP^1$. In the case $n=11$ 
we recall that $X_1(11)$ is the elliptic curve $\nu^2 + \nu = \la^3 - \la^2$.
We write $\lala$ to indicate $\la$ in the cases $n=7,9$ and the
pair $\la,\nu$ in the case $n=11$. Let $C_\lala$ and $D_\lala$
be the corresponding pairs of $n$-isogenous curves.
By \cite[Exercise 8.13]{Silverman} $D_\lala$ 
has Weierstrass equation
\begin{align*}
n & =7  &  y^2 - (\la^2-\la-1)xy - (\la^3-\la^2)y &= x^3 - (\la^3-\la^2)x^2, \\
n &= 9  &  y^2 + (\la^3 - 3 \la^2 + 4 \la - 1) xy & +
    \la(\la -  1)^4 (\la^2 - \la + 1)  y \\ &&& = x^3 
+ \la(\la-1)(\la^2 - \la + 1) x^2,  \\
n &= 11 & y^2 + (\la \nu+2 \la - (\nu+1)^2) xy 
   & -  \la^2 \nu (\nu+1) (\la-\nu-1) y  \\ &&&
= x^3 -\la \nu (\nu+1) (\la-\nu-1)x^2.
\end{align*}
On each of these curves $P=(0,0)$ is a point of order $n$.
If we write the Weierstrass equation for $D_\lala$ as
$y^2 + a_1 xy + a_3 y = x^3+ a_2 x^2$ then by V\'elu's formulae 
\cite{V1} the $n$-isogenous curve $C_\lala$ has Weierstrass equation 
\begin{equation}
\label{Clam}
 y^2 + a_1 xy + a_3 y = x^3 + a_2 x^2 -5t x - (a_1^2 + 4 a_2)t - 7 w
\end{equation}
where
$t = 6 s_2 + (a_1^2 + 4 a_2) s_1 + a_1 a_3 s_0$,
$w = 10 s_3 + 2(a_1^2 + 4 a_2) s_2 + 3 a_1 a_3 s_1 + a_3^2 s_0$ 
and $s_k = \sum_{j=1}^{(n-1)/2} x(j P)^k$.
The Weierstrass equations~(\ref{Clam}) have discriminant
\begin{equation}
\label{deltacla}
\begin{aligned}
n & = 7 & \Delta(C_\la) & = \la (\la-1) (\la^3 - 8 \la^2 + 5 \la + 1)^7 \\
n & = 9 & \Delta(C_\la) & = \la (\la-1) (\la^2-\la+1)^3 
(\la^3-6 \la^2+3 \la+1)^9 \\ 
n & = 11 & \Delta(C_{\la,\nu}) & = 
\la (\la-1)  (\la \nu + 2 \la^2 - 2 \la + 1) (\nu+1)^6 f(\la,\nu)^{11} 
\end{aligned}
\end{equation}
where $f(\la,\nu) = (-3 \la\nu + 2 \nu - \la^3 + 5 \la^2 - 5 \la + 1)/(\la-1)$.

\begin{Lemma}
\label{lem:defq}
Let $E/K$ be an elliptic curve and $\iota : \mu_n \inj E[n]$ an inclusion
of Galois modules. Then there exists $Q \in E(\Kbar)[n]$ and 
$q \in K^\times/(K^\times)^n$ such that
\begin{enumerate}
\item $e_n(\iota(\zeta),Q) = \zeta^{-1}$ for all $\zeta \in \mu_n$,
\item $\sigma(Q) - Q = \iota (  \frac{\sigma \sqrt[n]{q} }{ \sqrt[n]{q}} )$ 
for all $\sigma \in \Gal(\Kbar/K)$, and
\item $K(E[n]) = K(\mu_n, \sqrt[n]{q})$.
\end{enumerate}
\end{Lemma}
\begin{Proof}
This follows by standard properties of the Weil pairing,
together with Hilbert's Theorem 90.
\end{Proof}

Taking $E = C_{\lala}$ we compute $q=q(\lala)$ as described in 
\cite[Section 1.2]{JEMS}.
\begin{equation}
\label{qla} 
q(\lala) = \left\{ \begin{array}{ll}
\la^4 (\la-1) & \text{ if } n = 7 \\
\la (\la-1)^7 (\la^2 - \la + 1)^3 & \text{ if } n = 9 \\
\la \nu^2 (\la-1) (\la-\nu-1)^3
& \text{ if } n = 11. 
\end{array} \right.
\end{equation}

We consider the following diagonal twists of $X(n)$ for $n=7,9,11$.
Recall in the case $n=11$ we take the singular locus of the Hessian
of the cubic form.
\begin{align*}
X[\xi_1,\xi_2,\xi_3] 
            & = \{\xi_1 x^3 y + \xi_2 y^3 z + \xi_3 
z^3 x = 0\} \subset \PP^2\\
X[\xi_1,\xi_2,\xi_3;\eta]  &= \left\{ 
\begin{aligned}
\xi_1 x^2 y + \xi_2 y^2 z + \xi_3 z^2 x & =  0 \\
\xi_1 \xi_2 x y^2 + \xi_2 \xi_3 y z^2 
+ \xi_1 \xi_3 z x^2 & =  \eta t^3 
\end{aligned} \right\} \subset \PP^3 \\
X[\xi_1, \xi_2, \xi_3, \xi_4,\xi_5] 
&\sim \{ \xi_1 v^2 w + \xi_2 w^2 x + \xi_3 x^2 y 
    + \xi_4 y^2 z + \xi_5 z^2 v = 0 \} \subset \PP^4.  
\end{align*}

\begin{Lemma}
\label{lem:deftheta}
The diagonal twists are determined up to 
$K$-isomorphism by $\theta \in K^\times/(K^\times)^n$ where
\[ \theta = \left\{ \begin{array}{ll} 
\xi_1 \xi_2^2 \xi_3^4 & \text{ if } n = 7 \\
\xi_1 \xi_2^7 \xi_3^4 \eta^3 & \text{ if } n = 9 \\
(\xi_1 \xi_2^5 \xi_3^3  \xi_4^4  \xi_5^9)^2 & \text{ if } n= 11.
\end{array} \right. \]
Moreover if we replace $\theta$ by $\theta^k$ where $k$ is a 
square in $(\Z/n\Z)^\times$ then we obtain isomorphic twists.
\end{Lemma}
\begin{Proof}
The first part is proved by rescaling the 
co-ordinates and the second part by cyclically permuting them.
(The expression for $\theta$ in the case $n=11$ has been squared to
simplify the statement of the next theorem.)
\end{Proof}

\begin{Theorem}
\label{thm:diagtwist}
Let $E/K$ be an elliptic curve and $\iota : \mu_n \inj E[n]$ an inclusion
of Galois modules. If $q \in K^\times/(K^\times)^n$ is as specified in 
Lemma~\ref{lem:defq} then
\begin{enumerate}
\item $X_E(n)$ is the diagonal twist with $\theta = q$, and
\item $X^-_E(n)$ is the diagonal twist with $\theta = q^{-1}$.
\end{enumerate}
\end{Theorem}
\begin{Proof}
Let $\chi: \Gal(\Kbar/K) \to \Z/n\Z$ be defined by $\sigma (\sqrt[n]{q})/
\sqrt[n]{q} = \zeta_n^{\chi(\sigma)}$. Then the 
symplectic isomorphism $\phi : E[n] \isom
\mu_n \times \Z/n\Z; \,\, \iota(\zeta_n^x) + y Q \mapsto (\zeta_n^{-x},y)$ 
satisfies
\[ \sigma(\phi) \phi^{-1} : ( \zeta_n^x,y) \mapsto 
(\zeta_n^{x + \chi(\sigma)y},y). \]
Under our identification $\Aut(\mu_n \times \Z/n\Z) \isom \SL_2(\Z/n\Z)$
this corresponds to $T^{\chi(\sigma)}$. If $\rho(T) \propto
\Diag(\zeta_n^{r_1}, \ldots, \zeta_n^{r_m})$ for some integers 
$r_1, \ldots, r_m$ then 
\[   h_1 =  \Diag( (\sqrt[n]{q})^{r_1}, \ldots, (\sqrt[n]{q})^{r_m}) \]
and $h_2 = h_1^{-1}$ satisfy the hypotheses of Lemma~\ref{TwistLem}.

If $n=7$ then $(r_1,r_2,r_3) = (4,2,1)$ and we deduce
\[ X_E(7) \isom X[q,1,1] \qquad X_E^{-}(7) \isom X[q^{-1},1,1] \]
If $n=9$ then $(r_1,r_2,r_3,r_4) = (5,2,8,0)$ and we deduce
\[ X_E(9) \isom X[1,1,q;q^{-1}] \qquad 
  X_E^{-}(9) \isom X[1,1,q^{-1};q] \]
If $n=11$ then $(r_1,r_2,r_3,r_4,r_5) = (6,10,2,7,8)$ and we deduce
\[ X_E(11) \isom X[1,1,q^{-1},1,1] \qquad 
  X_E^{-}(11) \isom X[1,1,q,1,1] \]
The descriptions of $X_E(n)$ and $X_E^-(n)$ now follow by 
Lemma~\ref{lem:deftheta}.
\end{Proof}

It is sometimes convenient to rewrite 
$X[\frac{1}{\xi_2},\tfrac{1}{\xi_3}, \tfrac{1}{\xi_1}; 
\tfrac{\eta}{\xi_1 \xi_2 \xi_3} ]$ as
\[ X^-[\xi_1,\xi_2,\xi_3;\eta] = \left\{ 
\begin{aligned}
\xi_1 \xi_3 x^2 y + \xi_1 \xi_2 y^2 z 
+ \xi_2 \xi_3 z^2 x & =  0 \\
\xi_1  x y^2 + \xi_2 y z^2 + \xi_3 z x^2 & =  \eta t^3 
\end{aligned} \right\} \subset \PP^3 
\]
For this twist we have $\theta = \xi_1^2 \xi_2^5 \xi_3^8 
\eta^3 \mod{(K^\times)^9}$.

\begin{Corollary} 
\label{cor:diag}
Let $E/K$ be an elliptic curve whose $n$-torsion contains
a copy of the Galois module $\mu_n$. Then writing $E = C_\lala$ we have
\begin{align*}
X_E(7) & = X[\la,\la-1,1] & X_E(9) & = X[\la,\la-1,1; \la^2 - \la +1] \\
X^-_E(7) & = X[\la-1,\la,\la(\la-1)] & X^-_E(9) 
& = X^-[\la,\la-1,1; 1/(\la^2 - \la+1)] 
\end{align*} 
and 
\begin{align*}
X_E(11) & = X[\la^2(\la-1)^2,1,(\la-\nu-1)^2,\nu,1] \\ 
X^-_E(11) & = X[\la(\la-1),1,\la-\nu-1,\nu,\nu] 
\end{align*} 
\end{Corollary}
\begin{Proof} This follows from Theorem~\ref{thm:diagtwist} 
and the formula~(\ref{qla}) for $q(\lala)$.
\end{Proof}

\begin{Remark} The formula for $X_E(7)$ in Corollary~\ref{cor:diag}
was found by Halberstadt and Kraus \cite[Theorem 7.1]{HK} by specialising
the formula in Theorem~\ref{MainThm7} and then making a (rather
complicated) change of co-ordinates. We worked out the analogue of
this in the case $n=9$ before discovering the simpler proof 
presented here.
\end{Remark}

\section{Minimisation and reduction}
\label{sec:minred}

Let $n=7,9,11$ and $m = (n-1)/2$. Given an elliptic curve
$E/\Q$ the formulae in Section~\ref{sec:statres} give equations
for $X_E(n) \subset \PP^{m-1}$ and $X_E^-(n) \subset \PP^{m-1}$. 
We search for elliptic curves $n$-congruent to $E$ by searching
for $\Q$-rational points on these curves. It helps
with this search if we first make a change of co-ordinates over $\Q$
to simplify the equations, i.e. so that they have small
integer coefficients. Following \cite{CFS} this task naturally falls
into two parts, called minimisation and reduction. In minimisation one
seeks to remove primes from a suitably defined invariant. Then reduction,
which may be thought of as the analogue of minimisation at the infinite 
place, makes a final $\GL_m(\Z)$-transformation.

\subsection{The invariant}

\begin{Definition}
\label{def:inv}
We split into the cases $n=7,9,11$. \\
Case $n=7$. The invariant of a twisted form of $F= x^3 y + y^3 z + z^3 x$ is 
\[ \Psi( \mu (F \circ M)) = \mu^3 (\det M)^4. \]
Case $n=9$. The invariant of a twisted form of $(F_1,F_2) = 
(x^2 y + y^2 z + z^2 x, x y^2 + y z^2 + z x^2 - t^3)$ is
\[ \Psi( (\alpha F_1 + \beta F_2) \circ M, 
    (\gamma F_1 + \delta F_2) \circ M) = (\alpha \delta - \beta \gamma)^6 
(\det M)^9. \]
Case $n=11$. The invariant of a twisted form of 
$F= v^2 w + w^2 x + x^2 y + y^2 z + z^2 v$ is 
\[ \Psi( \mu (F \circ M)) = \mu^5 (\det M)^3. \]
\end{Definition}

\begin{Lemma} Let $\FF$ be one of the twisted forms in 
Definition~\ref{def:inv}. Then
\begin{enumerate}
\item $\Psi(\FF)$ is well-defined, i.e. it is independent of the choice
of $M \in \GL_m(\Kbar)$.
\item If $\FF$ has coefficients in $K$ then $\Psi(\FF) \in K$.
\end{enumerate}
\end{Lemma}
\begin{Proof} (i) This is easy to check for $M$ a scalar matrix.
In general we use that $\Aut(X(n)) \isom \PSL_2(\Z/n\Z)$.
This is proved in \cite[Lemma 20.4]{AR} for $n \ge 7$ prime, 
and the same proof works for $n=9$.
We are reduced to considering $M = \rho(\gamma)$ for some
$\gamma \in \SL_2(\Z/n\Z)$. It now suffices to recall (see the proof 
of Proposition~\ref{rho-lifts}) that the only 
$1$-dimensional characters of  $\SL_2(\Z/n\Z)$ are those of
order $3$ in the case $n=9$. This explains why in 
the case $n=9$ 
we defined the invariant as $(\alpha \delta - \beta \gamma)^6 
(\det M)^9$ and not just 
$(\alpha \delta - \beta \gamma)^2 (\det M)^3$. \\
(ii) This follows from (i) by Galois theory.
\end{Proof}

\begin{Remark}
\label{rem:inv}
(i) In the case $n=7$ it is shown in \cite[Section 7.1]{PSS} 
that $\Psi(\FF)$ is an 
integer coefficient polynomial in the coefficients of $\FF$. We expect 
that similar formulae exist in the cases $n=9$ and $n=11$. \\
(ii) The twisted forms in Theorems~\ref{MainThm7}, \ref{MainThm9}
and \ref{MainThm11} have the following invariants. 
These were computed by following the proofs
in Section~\ref{invthy}.
\[ \begin{array}{lcc}
& X_E(n) & X^-_E(n) \\
n = 7 & -4(4a^3 + 27b^2) & 16(4 a^3 + 27b^2)^2 \\
n = 9 & -2^{10}(4a^3 + 27b^2)^4 & 2^{8}(4 a^3 + 27b^2)^5 \\
n = 11 & -4(4a^3 + 27b^2)^2 & 8(4 a^3 + 27b^2). 
\end{array} \] 
(iii) The diagonal twists studied in Section~\ref{sec:diag} have invariants
\[ \begin{array}{lll}
n = 7 & X[ \xi_1,\xi_2,\xi_3] & \Psi = \xi_1 \xi_2 \xi_3 \\
n = 9 & X[ \xi_1,\xi_2,\xi_3;\eta] & \Psi = (\xi_1 \xi_2 \xi_3)^4 \eta^3 \\
    & X^-[ \xi_1,\xi_2,\xi_3;\eta] & \Psi = (\xi_1 \xi_2 \xi_3)^5 \eta^3 \\
n = 11 & X[ \xi_1,\xi_2,\xi_3,\xi_4,\xi_5] & \Psi = \xi_1 \xi_2 \xi_3 
\xi_4 \xi_5. 
\end{array} \]
\end{Remark}

\subsection{Minimisation}

The level of a model $\FF$ at a prime $p$ is the $p$-adic valuation 
of the invariant, i.e. $v_p(\Psi(\FF))$.
We seek to make a change of co-ordinates that
minimises the level. 
This is a local problem.

\begin{Theorem}
\label{thm:min}
Let $E$ be an elliptic curve over $\Q_p$ with $p \not=2,3$.
If $n = 7,9,11$ then $X_E(n)$ and $X_E^{-}(n)$ admit 
models (with coefficients in $\Z_p$) 
with the following levels. 
\[ \begin{array}{l|cc|cc} 
\multicolumn{1}{c|}{{\text{Kodaira Symbol}}} &  X_E(7) & X_E^{-}(7)  
&  X_E(11) & X_E^{-}(11)\\ \hline
{\rm I}_m,{\rm I}^*_m \quad  m \equiv 0 \pmod{n} & 0 & 0 & 0 & 0 \\
{\rm I}_m,{\rm I}^*_m  \quad (m/n) = + 1 & 1 & 2 & 2 & 1 \\
{\rm I}_m,{\rm I}^*_m \quad (m/n) = - 1 & 2 & 1 & 1 & 2 \\ 
{\rm II}, \, {\rm II}^*, \,  {\rm III}, \, {\rm III}^*, \, 
{\rm IV}, \, {\rm IV}^* & 2 & 2 & 2 & 2 
\end{array} \]
\[ \begin{array}{l|cc||c|cc} 
\multicolumn{1}{c|}{{\text{Kodaira Symbol}}} &  X_E(9) & X_E^{-}(9) &
{\text{Kodaira Symbol}} &  X_E(9) & X_E^{-}(9)\\ \hline
{\rm I}_m,{\rm I}^*_m  
 \quad m \equiv 0 \pmod{9} & 0 & 0 & 
{\rm II}, {\rm IV}^* & 5 & 7 \\
{\rm I}_m,{\rm I}^*_m  \quad m \equiv 3,6 \pmod{9} & 3 & 3 & 
{\rm III}, {\rm III}^* & 6 & 6 \\
{\rm I}_m,{\rm I}^*_m \quad m \equiv 1 \pmod{3} & 4 & 5 & 
{\rm IV}, {\rm II}^* & 7 & 5  \\
{\rm I}_m,{\rm I}^*_m  \quad m \equiv 2 \pmod{3} & 5 & 4 \\
\end{array} \]
\end{Theorem}

\begin{Proof}
Replacing $E$ by a quadratic twist does not change the curves 
$X_E(n)$ and $X_E^-(n)$. So we may assume that either (i) $E$ has
good reduction, or (ii) $E$ has split multiplicative reduction, or 
(iii) $E$ has additive reduction of type ${\rm II}$, ${\rm III}$
or ${\rm IV}$. We split into these three cases. \\
(i) If $E$ has good reduction then by Remark~\ref{rem:inv}(ii) 
the formulae in Theorems~\ref{MainThm7}, \ref{MainThm9} 
and~\ref{MainThm11} 
give models for $X_E(n)$ and $X^{-}_E(n)$ of level $0$. \\
(ii) If $E$ has split multiplicative reduction then by the Tate
parametrisation  $E[n]$ contains a copy of the Galois module $\mu_n$.
So we may apply the results of Section~\ref{sec:diag}.
More precisely
if $E$ has Kodaira symbol $I_m$ then $E(\Qbar_p) \isom 
\Qbar_p^\times/q^\Z$ for some $q \in \Q_p$ with $v_p(q)=m \ge 1$. 
It may be verified by computing the Weil pairing that the conditions
of Lemma~\ref{lem:defq} are satisfied. So by Theorem~\ref{thm:diagtwist}
the curves $X_E(n)$ and $X^{-}_E(n)$ are the diagonal twists 
with $\theta = q$ or $q^{-1}$. 
If $n=7$ or $11$ then by Lemma~\ref{lem:deftheta} they
admit models of the form $X[\xi_1,\xi_2,\ldots]$ where at most
two of the $\xi_i$ have $p$-adic valuation $1$ and the rest are units.
We read off the level from Remark~\ref{rem:inv}(iii). 
The case $n=9$ is similar. \\
(iii) Finally suppose $E$ has Kodaira symbol ${\rm II}$, 
${\rm III}$ or ${\rm IV}$. In terms of a minimal Weierstrass 
equation for $E$, say $y^2 = x^3 + a x + b$, these are the cases
\[ \begin{array}{cccc}
{\rm II} && v_p(a) \ge 1 & v_p(b) = 1, \\
{\rm III} && v_p(a) = 1 & v_p(b) \ge 2, \\
{\rm IV} && v_p(a) \ge 2 & v_p(b) = 2.
\end{array} \]
Integer-coefficient models of the required level are obtained by modifying
the formulae in Theorems~\ref{MainThm7}, \ref{MainThm9} 
and~\ref{MainThm11} as follows.
\[ \hspace{-0.7em} \begin{array}{c|cc|cc}
& X_E(7) &  X^-_E(7) & X_E(11) &  X^-_E(11) \\ \hline
{\rm II} & \FF(x,y,z) & \frac{1}{p^2} \GG(x,y,pz) 
& \frac{1}{p} \FF(pv+w,w,x,y,z) & \GG(v,w,x,y,z) \\
{\rm III} & \frac{1}{p^3} \FF(px,py,z) & \GG(x,\frac{1}{p} y,z) 
& \frac{1}{p^2} \FF(pv,pw,x,y,z) & \frac{1}{p^2}  \GG(pv,pw,px,y,z) \\
{\rm IV} & \frac{1}{p^2} \FF(x,py,z) & \frac{1}{p^2} \GG(x,\frac{1}{p} y,pz) 
& \frac{1}{p^3} \FF(pv,pw,x,py-x,z) 
& \frac{1}{p} \GG(pv,w,px,y, \frac{1}{p} z) 
\end{array} \hspace{-0.4em} \]
\[ \begin{array}{c|cc}
& X_E(9) &  X^-_E(9) \\ \hline
{\rm II} & 
\begin{array}{c}
\frac{1}{p^2} \FF_1(p x,p y,p z,t) \\
\frac{1}{p^3} \FF_2(p x,p y,p z,t) 
\end{array}
& 
\begin{array}{c}
\frac{1}{p} \GG_1(p x,y,z,t) \\ \medskip
\frac{1}{p} \GG_2(p x,y,z,t) 
\end{array}
\\ 
{\rm III} & 
\begin{array}{c}
\frac{1}{p^2} \FF_1(p x,p y,z,t) \\
\frac{1}{p^2} \FF_2(p x,p y,z,t) 
\end{array}
& 
\begin{array}{c}
\frac{1}{p} \GG_1(p x,y,z,t) \\ \medskip
\frac{1}{p^2} \GG_2(p x,y,z,t) 
\end{array}
\\
{\rm IV} & 
\begin{array}{c}
\frac{1}{p^4} \FF_1(p^2 x,p^2 y,p z,t) \\
\frac{1}{p^5} \FF_2(p^2 x,p^2 y,p z,t) 
\end{array}
& 
\begin{array}{c}
\frac{1}{p^2} \GG_1(p x,y,z,t) \\
\frac{1}{p^2} \GG_2(p x,y,z,t) 
\end{array}
\end{array} \]

\end{Proof}

To compute integer-coefficient models with level as specified
in Theorem~\ref{thm:min}
we could in principle follow the proof of the theorem.
In practice however it is simpler to use a range of ad hoc tricks.
We have not proved that these tricks always get down to the 
level specified in the theorem, but this does at least happen in
all numerical examples we have tried. 

Again on the basis of some numerical experimentation, 
we conjecture 
that the levels in Theorem~\ref{thm:min} are the minimal levels. 
At the primes $2$ and $3$ it seems the minimal levels can be larger. 
See Section~\ref{sec:ex} for some examples.
 
\subsection{Reduction}
In Section~\ref{sec:gpaction} we saw that 
the action of $\SL_2(\Z/n\Z)$ on 
$X(n) \subset \PP^{m-1} = \PP(V)$ lifts to an irreducible 
representation $\rho : \SL_2(\Z/n\Z) \to \GL(V)$. So by the Weyl unitary
trick there is an $\SL_2(\Z/n\Z)$-invariant inner product on $V$, and
this is unique up to scalars.
To reduce our equations for $X_E(n)$ and $X_E^-(n)$ we run the
LLL algorithm on the Gram matrix for this inner product.
In the untwisted case it is clear that the $\SL_2(\Z/n\Z)$-invariant 
inner product is the standard one on $\C^m$, i.e. $\rho$ is a
unitary representation. So to compute the inner product it 
suffices to have numerical approximations to the matrices
$h_1$ and $h_2$ in Lemma~\ref{TwistLem}. Since our method 
for finding equations for $X_E(n)$ and $X_E^-(n)$ involved
explicitly computing $h_1$ and $h_2$ this is obviously something
we can do. To use these formulae we need numerical approximations
for $a,b,c\ldots$ in Lemmas~\ref{Lem7}, \ref{Lem9} and \ref{Lem11}.
These are computed by evaluating suitable $q$-expansions.
See \cite{E} or \cite{HK} for the relevant formulae in the case $n=7$.

\section{Modular interpretation}
\label{modint}

In Section~\ref{invthy} we gave explicit formulae for $X_E(n)$
and $X^-_E(n)$ for $n=7,9,11$. In this section we give equations
for the families of curves they parametrise. 

\subsection{Computing the $j$-invariant}
\label{sec:jinv}
We first give formulae for the $j$-map $X_E(n) \to \PP^1$
and $X_E^{-}(n) \to \PP^1$. This is sufficient for some applications:
see for example \cite{PSS}. In each case we 
adapt the formulae in Lemmas~\ref{lem:j7}, \ref{lem:j9} and~\ref{lem:j11} 
by writing them in an way that behaves well under all 
changes of co-ordinates.

Case $n=7$. Let $X = \{ \FF = 0 \} \subset \PP^2$ be a twist of $X(7)$. 
Starting with $\FF$ in place of the Klein quartic $F$ 
the formulae in Section~\ref{sec:formulae7} define polynomials
$H(\FF)$, $c_4(\FF)$ and $c_6(\FF)$.
If $\FF = \mu (F \circ M)$ then $\Psi(\FF) = \mu^3 (\det M)^4$ and
\begin{equation}
\label{covprop7}
\begin{aligned}
H(\FF) & = \mu^3 (\det M)^2 (H \circ M) \\
c_4(\FF) & = \mu^8 (\det M)^6 (c_4 \circ M) \\
c_6(\FF) & = \mu^{12} (\det M)^9 (c_6 \circ M) 
\end{aligned}
\end{equation}
As noted in \cite{PSS} the syzygy $c_4^3 - c_6^2 \equiv 1728 H^7 
\pmod{F}$ becomes
\[  c_4(\FF)^3 - c_6(\FF)^2 \equiv 1728 \, \Psi(\FF) \, H(\FF)^7 \pmod{\FF}. \] 
In particular the $j$-map
$X \to \PP^1$ is given by
\[ j = \frac{c_4(\FF)^3}{\Psi(\FF) H(\FF)^7}. \]

Case $n=11$. Let $X \subset \PP^4$ be a twist of 
$X(11)$ given as the singular locus of the Hessian of a cubic form
$\FF = \FF(v,w,x,y,z)$. Starting with $\FF$ in place of the cubic
form $F = v^2 w + w^2 x + x^2 y + y^2 z + z^2 v$ the formulae in 
Section~\ref{sec:formulae11} define polynomials $H(\FF)$ and $c_4(\FF)$.
If $\FF = \mu (F \circ M)$ then $\Psi(\FF) = \mu^5 (\det M)^3$ and
\begin{equation}
\label{covprop11}
\begin{aligned}
H(\FF) & = \mu^5 (\det M)^2 (H \circ M) \\
c_4(\FF) & = \mu^{17} (\det M)^8 (c_4 \circ M) 
\end{aligned}
\end{equation}
By Lemma~\ref{lem:j11} the $j$-map $X \to \PP^1$ is given by
\[ j = \frac{c_4(\FF)^3}{ \Psi(\FF)^8 \FF^{11}}. \]

Case $n=9$. Our construction of $c_4(a,b,c,d)$
and $c_6(a,b,c,d)$ in Section~\ref{sec:formulae11} does not
immediately generalise to twists of $X(9)$. 
Instead we exploit the fact that
the pencil of cubics defining $X(9) \subset \PP^3$ is naturally
a copy of $X(3) \isom \PP^1$. The notation $X_M$ for a twist of 
$X(n)$ was introduced in Section~\ref{sec:defcrvs}.

\begin{Theorem}
\label{thm:9->3}
Let $X_M = \{ \FF_1 = \FF_2 = 0 \} \subset \PP^3$ be a twist of $X(9)$
where $M$ is a symplectic Galois module with $M \isom (\Z/9\Z)^2$ 
as an abelian group. Then
\begin{enumerate}
\item Writing $\HH$ for the determinant of the 
matrix of second partial derivatives we have
\[ \HH(r \FF_1 + s \FF_2) = f(r,s) D(a,b,c,d) \]
where $f$ and $D$ are homogeneous polynomials of degree $4$.
\item We have $X_{M[3]} \isom \PP^1$ with cusps at the roots
of $f(r,s) = 0$. 
\item For $P \in X_M$ with tangent line $P + t Q$ write
$\FF_i(P+ t Q) = \gamma_i t^2 + \delta_i t^3$ for $i=1,2$. 
Then the forgetful map $X_M \to X_{M[3]}$ is 
$P \mapsto (-\gamma_2:\gamma_1)$. 
\end{enumerate}
\end{Theorem}

\begin{Proof}
We first prove the theorem in the case $M \isom \mu_9 \times \Z/9\Z$.\\
(i) Taking $F_1 = a^2b + b^2 c + c^2 a$ and 
$F_2 = a b^2 + b c^2 + c a^2 - d^3$ we compute
\[ \HH(r F_1 + s F_2) 
    = 48 (r^3 + s^3) s (a^3 + b^3 + c^3 - 3 a b c) d. \]
(ii) 
Since\footnote{The pairing $e_9$ on $M$ induces a 
pairing $e_3$ on $M[3]$ by the
rule $e_3(3S,T) = e_9(S,T)$ for all $S \in M$ and $T \in M[3]$.} 
$M[3] \isom \mu_3 \times \Z/3\Z$ we have $X_{M[3]} = X(3)$.
We recall from Section~\ref{sec:n=3} that $X(3) \isom \PP^1$ with 
cusps at the roots of $A(27 A^3 + B^3)=0$. Our two choices of co-ordinates
on $X(3)$ are now related by $(r:s) = (B:3A)$. \\
(iii) We temporarily write $a_1, a_2, a_3, a_4$ for $a,b,c,d$ and let 
$\Lambda_2$ be the $4 \times 4$ alternating matrix with $(i,j)$ entry
\[ \frac{\partial F_1}{\partial a_l} \frac{\partial F_2}{\partial a_k} 
-\frac{\partial F_1}{\partial a_k} \frac{\partial F_2}{\partial a_l} \]
where $(i,j,k,l)$ is an even permutation of $(1,2,3,4)$. This is 
the same as the matrix $\Lambda_2$ in Section~\ref{sec:invthy9}.
By direct calculation we find
\[ \Lambda_2 \begin{pmatrix} \frac{\partial^2 F_i}{\partial a^2} 
& \ldots &\frac{\partial^2 F_i}{\partial a \partial d} \\ 
\vdots & & \vdots \\
\frac{\partial^2 F_i}{\partial a \partial d}
 & \ldots & \frac{\partial^2 F_i}{\partial d^2} 
\end{pmatrix}
\Lambda_2
 \equiv 
\gamma_i  \, D 
\begin{pmatrix} a \\ \vdots \\ d \end{pmatrix}
\begin{pmatrix} a & \cdots & d \end{pmatrix}
\mod{(F_1,F_2)} 
\]
for $i=1,2$,
where $\gamma_1 = 18 d^3$, $\gamma_2 = -6 (a^3 + b^3 + c^3 + 6 a b c)$ 
and $D = -(a^3 + b^3 + c^3 - 3abc)d$.

In Corollary~\ref{cor:forget} below 
we show that the forgetful map $X(9) \to X(3)$ is 
\begin{equation}
\label{forget:9->3}
 (a:b:c:d) \mapsto (A:B) = (d^3:a^3 + b^3 + c^3 + 6 abc). 
\end{equation}
Thus $(B:3A) = (-\gamma_2:\gamma_1)$ as required.

Relaxing our restriction on $M$, the general case of the theorem
follows since the forgetful map $X_M \to X_{M[3]}$ may be characterised 
geometrically, i.e. it quotients out by the symplectic automorphisms 
of $M$ that act trivially on $M[3]$.
\end{Proof}

\begin{Remark}
\label{rem:9to3}
If we only use Theorem~\ref{thm:9->3} to compute the $j$-invariant, then
we can replace the forward reference~(\ref{forget:9->3}) in the above proof by
the observation that the polynomials $c_4(A,B)$ and $c_6(A,B)$ 
in Section~\ref{sec:n=3} are related to the polynomials $c_4(a,b,c,d)$ 
and $c_6(a,b,c,d)$ in Section~\ref{sec:invthy9} by
\begin{align*}
c_4(d^3, a^3 + b^3 + c^3 + 6 a b c) & \equiv c_4(a,b,c,d) \mod{(F_1,F_2)}, \\
c_6(d^3, a^3 + b^3 + c^3 + 6 a b c) & \equiv c_6(a,b,c,d) \mod{(F_1,F_2)}. 
\end{align*}
\end{Remark}

Since we already gave formulae for the universal families above
$Y_E(3)$ and $Y_E^{-}(3)$ in Theorem~\ref{MainThm3}, the following 
corollary gives formulae for the universal families above $Y_E(9)$ 
and $Y_E^{-}(9)$.

\begin{Corollary}
Let $E/K$ be the elliptic curve $y^2 = x^3 - 27 c_4 x - 54 c_6$. 
Let $X_E(3)$ and $X_E^{-}(3)$ be as given in Theorem~\ref{MainThm3},
and let $X_E(9)$ and $X_E^-(9)$ be as given in Theorem~\ref{MainThm9}
with $a = -27c_4$ and $b = -54 c_6$. Then the forgetful maps
$X_E(9) \to X_E(3)$ and $X_E^-(9) \to X_E^-(3)$ are given by
$(\lambda : \mu )= (\gamma_2: 3\gamma_1)$
where $\gamma_1, \gamma_2$ are 
computed by the tangent line construction in Theorem~\ref{thm:9->3}(iii).
\end{Corollary}
\begin{Proof}
(i) We have $X_E(9) = \{ \FF_1 = \FF_2 = 0 \} \subset 
\PP^3$ where $\FF_1$ and $\FF_2$ are as
given in Theorem~\ref{MainThm9}.
We compute
$ \HH(r \FF_1 + s \FF_2) = f(r,s) D(x,y,z,t) $
with 
\begin{equation}
\label{eqnD3a}
  f(r,s) = r^4 + 2 a r^2 s^2 - 4 b r s^3 - \tfrac{1}{3} a^2 s^4.
\end{equation}
On the other hand the family of curves parametrised by $X_E(3)$
in Theorem~\ref{MainThm3} has cusps at the roots of 
\begin{equation}
\label{eqnD3b}
 \DD(\la,\mu) = 
    \la^4 - 6 c_4 \la^2 \mu^2 - 8 c_6 \la \mu^3 - 3 c_4^2 \mu^4.
\end{equation}
Comparing~(\ref{eqnD3a}) and~(\ref{eqnD3b}) we see that our two 
choices of co-ordinates on $X_E(3) \isom \PP^1$ are related by 
$(\la:\mu) = (r:-3 s)$. Taking $(r:s) = (-\gamma_2:\gamma_1)$ gives
the result.

(ii) We have $X_E^-(9) = \{ \GG_1 = \GG_2 = 0 \} \subset 
\PP^3$ where $\GG_1$ and $\GG_2$ are as given in 
Theorem~\ref{MainThm9}. We compute 
$ \HH(r \GG_1 + s \GG_2) = f(r,s) D(x,y,z,t) $
with 
\begin{equation}
\label{eqnD3ma}
  f(r,s) = a r^4 - 6 b r^3 s - 2 a^2 r^2 s^2 + 2 a b r s^3 
                        - (\tfrac{1}{3} a^3 + 3 b^2) s^4.
\end{equation}
On the other hand the family of curves parametrised by $X_E^-(3)$
in Theorem~\ref{MainThm3} has cusps at the roots of 
\begin{equation}
\label{eqnD3mb}
\cc_4(\la,\mu) = c_4 \la^4 + 4 c_6 \la^3 \mu + 6 c_4^2 \la^2 \mu^2 
+ 4 c_4 c_6 \la \mu^3 - (3 c_4^3 - 4 c_6^2) \mu^4.
\end{equation}
Comparing~(\ref{eqnD3ma}) and~(\ref{eqnD3mb}) we see that our 
choices of co-ordinates on $X_E^-(3) \isom \PP^1$ are related by 
$(\la:\mu) = (r:-3s)$. 
Taking $(r:s) = (-\gamma_2:\gamma_1)$ gives the result.
\end{Proof}

\subsection{Modular interpretation of $X(n)$}
\label{sec:famXn}

In Section~\ref{klein-xn} we gave equations for $X(n)$. 
We recall from \cite{mythesis} that  
(analogous to Definition~\ref{def:zn}) the elliptic curve 
$E \subset \PP^{n-1}$ above $(0:a_1:a_2: \ldots : -a_2:-a_1) \in Y(n)$ 
has equations
\begin{equation*}
 \rank ( a_{i-j} x_{i+j} )_{i,j=0}^{n-1} \le 2.
\end{equation*}
Taking $n=7,9,11$ we now put this curve in Weierstrass form.
Notice that the Weierstrass equations have coefficients 
that are homogeneous polynomials of degree $4m$ and $6m$ 
for some integer $m$.

\begin{Theorem}
\label{weqns}
We split into the cases $n=7,9,11$.
\begin{enumerate}
\item
The family of curves parametrised by 
$X(7) = \{ a^3 b + b^3 c + c^3 a = 0 \} 
\subset \PP^2$ is
\begin{equation}
\label{dagger}
 y^2 = x^3 - 27 (abc)^2 c_4(a,b,c) x - 54 (abc)^3 c_6(a,b,c)
\end{equation}
where $c_4,c_6 \in K[a,b,c]$ are as defined 
in Section~\ref{sec:formulae7}.
\item 
The family of curves parametrised by 
$X(9) = \{ 
F_1 = F_2 = 0 \} \subset \PP^3$ is
\begin{equation*}
 y^2 = x^3 - 27 c_4(a,b,c,d) x - 54 c_6(a,b,c,d)
\end{equation*}
where $c_4,c_6 \in K[a,b,c,d]$ are as defined 
in Section~\ref{sec:invthy9}.
\item
The family of curves parametrised by $X(11) \subset \PP^4$ is 
\begin{equation}
\label{dagger11}
 y^2 = x^3 - 27 (abcde)c_4(a,b,c,d,e) x
            - 54 (abcde) \widetilde{c}_6(a,b,c,d,e).  
\end{equation}
where $c_4,\widetilde{c}_6 \in K[a,b,c,d,e]$ are as defined 
in Section~\ref{sec:formulae11}.
\end{enumerate}
\end{Theorem}
\begin{Proof}
In Section~\ref{sec:diag} we wrote $\lala$ for a co-ordinate
(or pair of co-ordinates) on $X_1(n)$ for $n=7,9,11$. With
$q(\lala)$ as defined in~(\ref{qla}), we see by Lemma~\ref{lem:defq}
that $X(n)$ is birational to $\{q(\lala)=\tau^n\} \subset X_1(n) 
\times \Gm$. In the case $n=7$ an explicit birational map is given
in \cite[Section 2.2]{JEMS}. Applying the same method 
for $n=9,11$ we obtain
\begin{align*}
n &= 7 & (a:b:c) &\mapsto (\la,\tau) = (-ac^2/b^3,ac/b^2), \\
n &= 9 & (a:b:c:d) &\mapsto (\la,\tau) = (-ac/b^2,a^2 d/ b^2 c ), \\
n &= 11 & (a:b:c:d:e) & \mapsto 
(\la,\nu,\tau) = (-a b d/ c^2 e, a b^3/ c^3 e, -a b/c^2). 
\end{align*}
We checked directly that these are birational maps, and 
that the cusps of $X(n)$, i.e. $(1:0: \ldots :0)$ and its translates
under the action of $\SL_2(\Z/n\Z)$, map to the cusps of $X_1(n)$,
i.e. the roots of~(\ref{deltacla}).

Let $c_4(\lala)$ and $c_6(\lala)$ be the invariants of the Weierstrass
equation for $C_\lala$ in Section~\ref{sec:diag}. Using Magma we compute
\[ n=7 \qquad   \begin{aligned}
c_4(-ac^2/b^2)& \equiv \xi_7^4 (abc)^2 c_4(a,b,c) \mod{(a^3b+b^3c+c^3a)}\\
c_6(-ac^2/b^2)& \equiv \xi_7^6 (abc)^3 c_6(a,b,c) \mod{(a^3b+b^3c+c^3a)} 
\end{aligned} \]
\[ n=9 \qquad  \begin{aligned}
c_4(-ac/b^2)& \equiv \xi_9^4  c_4(a,b,c,d) \mod{(F_1,F_2)} \\
c_6(-ac/b^2)& \equiv \xi_9^6  c_6(a,b,c,d) \mod{(F_1,F_2)} 
\end{aligned} \]
\[ n=11 \qquad   \begin{aligned}
c_4(-abd/c^2e,ab^3/c^3e)& \equiv \xi_{11}^4  (abcde)c_4(a,b,c,d,e) \mod{\I}\\
c_6(-abd/c^2e,ab^3/c^3e)& \equiv \xi_{11}^6  
(abcde) \widetilde{c}_6(a,b,c,d,e) \mod{\I}
\end{aligned} \]
where $\xi_7 = a/ b^5 c$, $\xi_9 = a^2/b^4 c$ and $\xi_{11} = a^3 b /c^6 e^2$.

See \cite[Section 3]{HK} for a sketch of an alternative proof
in the case $n=7$.
\end{Proof}

\begin{Corollary}
\label{cor:forget}
The forgetful map $X(9) \to X(3)$ is given by
\[ (a:b:c:d) \mapsto (d^3 : a^3 + b^3 + c^3 + 6 abc). \]
\end{Corollary}

\begin{Proof}
This follows from Theorem~\ref{weqns}(ii) and Remark~\ref{rem:9to3}.
\end{Proof}

\subsection{An alternative projective embedding}
We take $p \ge 5$ a prime and let $G= \PSL_2(\Z/p\Z)$ act 
on $X(p)$ in the usual way. 

\begin{Theorem}[Adler, Ramanan] 
\label{lem:AR}
The group of $G$-invariant divisor classes on $X(p)$ is free
of rank 1 generated by a divisor class $[\Lambda]$ of degree 
$(p^2 -1)/24$.
\end{Theorem}
\begin{Proof} See \cite[Theorem 24.1]{AR}.
\end{Proof}

Let $m = (p-1)/2$. Klein showed there are embeddings
$X(p) \subset \PP^{m-1}$ and $X(p) \subset \PP^{m}$ with linear $G$-action. 
The images are called the $z$-curve and the $A$-curve respectively.
The corresponding hyperplane sections are $(m-1) \Lambda$ and $m \Lambda$,
and indeed the divisor $\Lambda$ in Theorem~\ref{lem:AR} 
is constructed by taking the difference of these. 
It is conjectured that each of these embeddings is via a complete linear 
system (the WYSIWYG Hypothesis in \cite{AR}) 
and this is certainly known for $p=7$ and $p=11$.
The equations for $X(p)$ we have used so far (introduced in 
Section~\ref{klein-xn}) are for the $z$-curve. 
However in Sections~\ref{modint7} and~\ref{modint11} below
we also need the $A$-curve.

\begin{Remark}
\label{rem:noc6}
Let $\pi : {\mathfrak X} \to X$ be the universal family over $X = X(p)$.
It can be shown that $\pi_*(\Omega_{{\mathfrak X}/X}) \isom \OO (p \Lambda)$.
Hence the family of curves ${\mathfrak X} \to X$ has Weierstrass equation 
$y^2 = x^3 - 27 c_4 x - 54 c_6$ where
\[ c_k \in H^0(X(p), \OO(kp\Lambda)) \]
for $k=4,6$. If we realise $X(p) \subset \PP^{m-1}$ as the
$z$-curve then to write $c_k$ as a polynomial
(in $m$ variables) we need that $m-1 = (p-3)/2$ divides $kp$. 
Thus in the case $p=7$ we were
able to write $c_4$ and $c_6$ as polynomials in $K[a,b,c]$.
In the case $p=11$ we likewise constructed $c_4 \in K[a,b,c,d,e]$,
but there was no polynomial $c_6$.
\end{Remark}

Case $p=7$. The $z$-curve is the Klein quartic
\[ X(7) = \{ x^3 y + y^3 z + z^3 x = 0 \} \subset \PP^2. \]
The cusps of $X(7)$ are the $24$ points of inflection.
We recall from \cite{PSS} that the cusps are naturally 
partitioned into eight sets of three $\{P_1,P_2,P_3\}$ with
\[ P_1 + 3 P_2 \sim P_2 + 3 P_3 \sim P_3 + 3 P_1 \sim H \]
where $H \sim 2 \Lambda$ is the hyperplane section. We write 
$T_0, \ldots, T_7$ for the effective divisors of degree $3$ of
the form $P_1 + P_2 + P_3$ and note that one of these divisors,
say $T_0$, satisfies $X(7) \cap \{xyz = 0\} = 4 T_0$.
We also recall from \cite{PSS} that $2 T_i \sim 2 T_j$ for 
all $0 \le i,j \le 7$. 
It follows by  Theorem~\ref{lem:AR} that $2 T_0 \sim 3 \Lambda$.
Since $3 \Lambda \sim 3H - 2T_0$ and $\CL(3H - 2T_0)$ has basis
$x^2 y, y^2 x, z^2x, xyz$ the $A$-curve is the image of
\[ X(7) \to \PP^3 \, ; \quad (x:y:z) \mapsto  
(t_1:t_2:t_3:t_4) = (x^2 y:y^2 z:z^2 x:x y z) \]
with equations
\[ \rank
\begin{pmatrix}
 t_1 &  0  & t_4 & -t_2 \\
 t_2 & -t_3 &  0 &  t_4 \\
 t_3 &  t_4 & -t_1 &  0 
\end{pmatrix} \le 2. \]

Case $p=11$. The $z$-curve is the singular locus of the Hessian of
\[ \{ F = v^2 w + w^2 x + x^2 y + y^2 z + z^2 v = 0 \} \subset \PP^4. \]
We write $H \sim 4 \Lambda$ for the hyperplane section. 
The cusps are the $60$ points of intersection of $X(11)$ with
$\{F = 0\}$. They are naturally 
partitioned into twelve sets of five $\{P_1,\ldots ,P_5\}$ with
\[  P_1  + 6 P_3 + 3 P_4 + 10 P_5 \sim H \]
and likewise under all cyclic permutations of the $P_i$.
We write $T_0, \ldots, T_{11}$ for the effective divisors of degree $5$ of
the form $P_1 + \ldots + P_5$ and note that one of these divisors,
say $T_0$, satisfies $X(11) \cap \{vwxyz = 0\} = 20 T_0$.
It may be shown that $5 T_i \sim 5 T_j$ for all $0 \le i,j \le 11$ 
and hence $5 T_0 \sim 5 \Lambda$ by Theorem~\ref{lem:AR}.
Since $5 \Lambda \sim 5H - 15 T_0$ we find by computing a basis
for $\CL(5H - 15 T_0)$ that the $A$-curve is the image of
the morphism $X(11) \to \PP^5$ given by
\[  (t_1:t_2:t_3:t_4:t_5:t_6) = 
(v^2 w x z:v w^2 x y:w x^2 y z:v x y^2 z:v w y z^2:v w x y z). \]
It is shown in \cite[Theorem 51.1]{AR}, and we checked using Magma, that 
this is the singular locus of the quartic hypersurface
\[  t_6^4 - (t_1^2 t_2 + t_2^2 t_3 + t_3^2 t_4 + t_4^2 t_5 + t_5^2 t_1)t_6
  + t_1^2 t_3 t_5 + t_2^2 t_4 t_1 + t_3^2 t_5 t_2 + t_4^2 t_1 t_3 + t_5^2 t_2 t_4 = 0. \]

\subsection{Formulae in the case $n=7$}
\label{modint7}

\begin{Theorem}
\label{uptoquad}
Let $ {\mathcal X} = \{\FF = 0 \} \subset \PP^2$ be a twist of 
the Klein quartic, with hyperplane section $H$. Let $T = P_1 + P_2 + P_3$
where $P_1, P_2, P_3$ are points of inflection on ${\mathcal X}$ with
\[ P_1 + 3 P_2 \sim P_2 + 3 P_3 \sim P_3 + 3 P_1 \sim H. \]
Let $d \in K[x,y,z]$ be a cubic form with $\{d = 0\}$ meeting ${\mathcal X}$
in a divisor $2D$ with $D \sim 2T$. Then there is a $\Gal(\Kbar/K)$-module
$M$ such that for every field extension $L/K$ and rational point
$P=(x:y:z) \in  {\mathcal X}(L) \setminus \{ d= 0 \}$, 
not a point of inflection, the elliptic curve
\begin{equation}
\label{star}
  Y^2 = X^3 - 27 \frac{c_4(\FF)(x,y,z)}{d(x,y,z)^2} X - 54
  \frac{c_6(\FF)(x,y,z)}{d(x,y,z)^3} 
\end{equation}
has $7$-torsion isomorphic to $M$ as a $\Gal(\Lbar/L)$-module.
\end{Theorem}
\begin{Proof}
We first note that if $d_1,d_2 \in K[x,y,z]$ are cubic forms meeting
${\mathcal X}$ 
in divisors $2 D_1$ and $2 D_2$ with $D_1 \sim D_2$ then $d_1/d_2$
is the square of a rational function and hence the elliptic 
surfaces~(\ref{star}) with $d = d_1$ and $d = d_2$ are isomorphic 
over $\Kbar$. Since ${\mathcal X}$ is a twist of the Klein quartic it follows 
(by taking $D=2T_0$ as defined in the last section)
that the elliptic surfaces~(\ref{dagger}) and~(\ref{star}) are 
isomorphic over $\Kbar$.
(Notice it does not matter whether we write the terms $d(x,y,z)$ in the 
numerator or in the denominator.) 
We are done by \cite[Proposition 2.1]{RubinSilverberg}.
\end{Proof}

In Theorem~\ref{families7} below we determine rational functions $d$ satisfying
the hypothesis of Theorem~\ref{uptoquad} 
in the cases $ \mathcal{X} = X_E(7)$
and $\mathcal{X} = X^-_E(7)$. We also show how to scale these functions
to give the quadratic twist with $M \isom E[7]$.

\begin{Remark}
Recall that $X_E(7)$ has a trivial
$K$-rational point corresponding to $E$ itself. 
Following~\cite{RubinSilverberg} one method for finding the right
quadratic twist would be to specialise at this point.
However this approach fails when $d$ vanishes at the trivial point.
Neither does the method generalise to $X^-_E(7)$. 
\end{Remark}

\begin{Theorem}
\label{families7}
Let $E/K$ be an elliptic curve with Weierstrass equation 
$y^2 = x^3 - 27 c_4 x - 54 c_6$ and let $\Delta = (c_4^3- c_6^2)/1728$. 
If $j(E) \not= 0,1728$
then the families of elliptic curves parametrised by $Y_E(7)$ and $Y_E^-(7)$
are given by~(\ref{star}) with $({\mathcal F},d) =(\BF,d_1)$ and 
$(\BG,d_2)$ where $\BF$ and $\BG$ are the quartics in
Theorem~\ref{Thm7} and 
\begin{align*}
d_1(x,y,z) & = -6 ( 3 x^2 + c_4 x z - 3 c_4 y^2 + c_6 y z) z \\
d_2(x,y,z) & = 2 \Delta (4 x^3 + c_4 x^2 z - 12 c_4 x y^2 - 2 c_6 x y z 
+ 8 c_6 y^3 + c_4^2 y^2 z + 200 \Delta z^3).
\end{align*}
\end{Theorem}
\begin{Proof}
We fix a symplectic isomorphism 
$\phi : E[7] \isom \mu_7 \times \Z/7\Z$ and let $(a:b:c)$ be the
$\Kbar$-point on $X(7)$ corresponding to $(E,\phi)$. 
As in the proof of 
Theorem~\ref{Thm7} we scale $(a,b,c)$ 
so that $c_4(a,b,c) = c_4$ and $c_6(a,b,c) = c_6$. 

Consideration of the action of $\SL_2(\Z/7\Z)$ on both the $z$-curve and
the $A$-curve suggests we start with the forms
\begin{align*}
s_1(x,y,z) & = 
   (a^2 c^3 - 2 a b^3 c) x^2 y + (a^3 b^2 - 2 a b c^3) y^2 z \\
  & \text{~~\hspace{2em}} ~+ (b^3 c^2 -2 a^3 b c) x z^2 
       + (a^3 c^2 + a^2 b^3 + b^2 c^3) x y z, \\
s_2(x,y,z) & = a^2 b x^2 y + b^2 c y^2 z + c^2 a z^2 x + 2 a b c x y z.
\end{align*}
We then let $r_1$ and $r_2$ be the unique cubic forms satisfying
\begin{equation}
\label{rs-identity}
 r_i(x,y,z) xyz  \equiv s_i(x,y,z)^2  \mod{(x^3y + y^3z + z^3x)}
\end{equation}
for $i=1,2$. The coefficients of $r_1$ and $r_2$
are homogeneous polynomials in $a,b,c$ of degrees $10$ and $6$. 
Recall that in the proof of Theorem~\ref{Thm7} we put
\begin{align*}
\BF (x,y,z)& = \frac{1}{H^3} 
F ( x \x + y ( \nabla F \times \nabla H) + z H \e), \\
\BG (x,y,z)& = \frac{1}{H^2} F ( x \nabla F + y (\x \times \e)
 + z H^2 \nabla H ).
\end{align*}
The cubics $d_1$ and $d_2$ in the statement of the theorem are likewise
found by putting
\begin{equation}
\label{d1d2}
\begin{aligned}
d_1 (x,y,z)& = \frac{1}{2 abc H^4}  
r_1 ( x \x + y ( \nabla F \times \nabla H) + z H \e), \\
d_2 (x,y,z)& = \frac{2 H^5}{abc} r_2 (  x \nabla F + y (\x \times \e) + z H^2 \nabla H ).
\end{aligned}
\end{equation}

It is clear from these constructions that $\{ d_1 = 0\}$ and $\{d_2= 0\}$
meet the corresponding twists of the Klein quartic in divisors of the 
form specified in Theorem~\ref{uptoquad}.
Hence our formulae for the families of elliptic curves parametrised by $Y_E(7)$
and $Y_E^-(7)$ are correct up to quadratic twist, say by $\delta \in K^\times$. 
It remains to show that $\delta$ is a square. As noted in \cite[Section 7.1]{HK}
it suffices to check this in the case $\phi : E[7] \isom \mu_7 \times \Z/7\Z$
is defined over $K$. Then $(a:b:c)$ is a $K$-rational point on $X(7)$.
We write $(a,b,c) = (\la a_0, \la b_0, \la c_0)$ with $a_0,b_0,c_0 \in K$.
By our earlier choice of scaling for $a,b,c$ we have $\la^7 \in K$.
Comparing the Weierstrass equation~(\ref{dagger})
for $E$ with that in the statement of the theorem we deduce that
$\la^7 a_0b_0c_0 \in (K^\times)^2$. Hence $a^7, b^7, c^7 \in K$ and 
$(abc)^7 \in (K^\times)^2$. 
Using~(\ref{det7}) and~(\ref{covprop7}) we compute
\begin{align*}
c_k( \BF )(x,y,z) & = (2^93^6)^{k/2}
    c_k(x \x + y ( \nabla F \times \nabla H) + z H \e) \\
c_k( \BG )(x,y,z) & = (2^93^6H^7)^{k/2}
    c_k(  x \nabla F + y (\x \times \e) + z H^2 \nabla H ). 
\end{align*}
for $k=4,6$. It follows by~(\ref{d1d2}) that
\begin{align*}
\frac{c_k(\BF)(x,y,z)}{ \,\,\, d_1(x,y,z)^{k/2}} 
& = \xi^k \frac{ 
c_k(x H\x + y H ( \nabla F \times \nabla H) + z H^2 \e)}
{((abc)^{6} r_1(x H \x + y H( \nabla F \times \nabla H) + z H^2 \e ))^{k/2}} \\
\frac{ c_k(\BG)(x,y,z)}{ \,\,\, d_2(x,y,z)^{k/2}} 
& = \eta^k \frac{ 
c_k(x H^3\nabla F + y H^3 (\x \times \e) + z H^5 \nabla H )}
{((abc)^{6} H^{10} 
  r_2(x H^3 \nabla F + y H^3 (\x \times \e) + z H^5 \nabla H ))^{k/2}}
\end{align*}
for some $\xi, \eta \in K^\times$.
The covariant columns $H \x$, $H(\nabla F \times \nabla H)$, $H^2 \e$
have degrees $7,14,21$ and the contravariant columns $H^3 \nabla F$, 
$H^3 (\x \times \e)$, $H^5 \nabla H$ have degrees $21, 28, 35$.
Since each column has degree a multiple of $7$, its evaluation at
$(a,b,c)$ is $K$-rational. Thus the families of curves 
in the statement of the theorem are $K$-isomorphic to
\[  Y^2  = X^3 - 27 \frac{c_4(x,y,z)}{((abc)^6 r_1(x,y,z))^2} X - 54
  \frac{c_6(x,y,z)}{((abc)^6 r_1(x,y,z))^3}  \]
and
\[ Y^2 = X^3 - 27 \frac{c_4(x,y,z)}{( (abc)^6 H^{10} r_2(x,y,z))^2} X - 54
  \frac{c_6(x,y,z)}{((abc)^6 H^{10} r_2(x,y,z))^3}. \]
To identify these with~(\ref{dagger}) we note that the cubic forms 
$(abc)^3 s_1(x,y,z)$ and $(abc)^3 H^5 s_2(x,y,z)$ have coefficients
in $K$ (since the degree of each coefficient is a multiple of $7$)
and then use~(\ref{rs-identity}).
\end{Proof}

Making the change of co-ordinates in Remark~\ref{tidyup7} we can replace
$d_1$ and $d_2$ by cubic forms that satisfy the conditions of 
Theorem~\ref{uptoquad} for $X_E(7) = \{ \FF = 0 \} \subset \PP^2$
and $X_E^-(7) = \{ \GG = 0 \} \subset \PP^2$ where $\FF$ and $\GG$
are the quartics in Theorem~\ref{MainThm7}. Moreover having 
found one such form we can use the Riemann-Roch
machinery in Magma to find further such forms. 

In the case of $X_E(7)$ we obtain a cubic form $d_{11}$ with $X_E(7) \cap
\{d_{11} = 0\} = 2D_1$ for some divisor $D_1 \sim 2T$. 
Then $\CL(3H - D_1)$
has basis 
\begin{align*}
d_{11} &= -2 (a x^2 + 3 b x z + 3 y^2 + 2 a y z) z, \\
d_{12} &= 2 (a x^2 + 3 b x z + 3 y^2 + 2 a y z) x, \\
d_{13} &= 4 (3 b x^2 - 2 a x y - 2 a^2 x z - 3 b y z - 2 a b z^2) z, \\
d_{14} &= 4 (a^2 x^2 + 3 b x y + 4 a b x z + a y^2 + 3 b^2 z^2) z.
\end{align*}
More generally we compute
cubic forms $d_{ij}$ for $1 \le i,j \le 4$
such that the matrix $(d_{ij})$ is symmetric and each 
$2 \times 2$ minor vanishes mod $\FF$. 
The remaining $d_{ij}$ are computed using
$d_{11} d_{ij} \equiv d_{1i} d_{1j} \pmod{\FF}$.
Then $X_E(7) \cap \{ d_{ij} = 0 \} = D_i + D_j$ where $D_1, \ldots, D_4$ 
are divisors all linearly equivalent to $2 T$. 
The family of elliptic curves parametrised by $Y_E(7)$ is now given
by~(\ref{star}) with $(\FF,d) = (\FF,d_{ii})$ for any $1 \le i \le 4$. 
The $A$-curve is the image of $X_E(7) \to \PP^3;$ $(x:y:z) 
\mapsto (d_{11}: \ldots :d_{14})$ with equations
\[ \rank \begin{pmatrix}   0 & t_3 & -t_4 & 2 a t_1 + t_4 \\
  t_1 & 2 a t_1 + t_4 & 2 b t_1 + a t_2 + a t_3 & 2 a t_2 + a t_3 \\
  t_2 & 2 b t_1 + a t_3 & -a^2 t_1 + b t_3 - a t_4 & 2 b t_2 - b t_3 - a t_4 
    \end{pmatrix} \le 2. \]

Our formula for the elliptic curve corresponding to $P \in Y_E(7)$
fails at points $P$ with $d_{ii}(P)=0$. These are the points whose image
on the $A$-curve lies on the co-ordinate hyperplane $\{t_i=0\}$. 
Therefore for any given point $P$ 
we have $d_{ii}(P) \not= 0$ for some $i$.
So unlike the treatment in \cite[Theorem 5.2]{HK}, where only the cubic form 
$d_{11}$ was given, we have found formulae that cover all cases.

In the case of $X_E^-(7)$ we likewise 
find cubic forms $d'_{ij}$ for $1 \le i,j \le 4$
such that the matrix $(d'_{ij})$ is symmetric and each 
$2 \times 2$ minor vanishes mod $\GG$. Explicitly 
\begin{align*}
d'_{11} &= -7 a x^2 y + 6 x^2 z + 3 a^2 y^3 - 8 a y^2 z + 3 y z^2, \\
d'_{22} &= 2 a x^3 + 12 b x^2 y - 2 a x y z - 3 a b y^3 + 6 b y^2 z, \\
d'_{33} &= 2 a^2 x y^2 - 10 a x y z + 6 x z^2 + 5 a b y^3 - 12 b y^2 z, \\
d'_{44} &= 2 a^2 x^2 y - 3 a x^2 z + 5 a b x y^2 - 12 b x y z
 - 3 a^2 y^2 z + 8 a y z^2 - 3 z^3.
\end{align*}
The remaining $d'_{ij}$ are computed using 
$d'_{11} d'_{ij} \equiv d'_{1i} d'_{1j} \pmod{\GG}$.
The family of elliptic curves parametrised by $Y_E^-(7)$ is now given
by~(\ref{star}) with $(\FF,d) = (\GG,\Delta d'_{ii})$ for any $1 \le i \le 4$. 
Exactly as before these formulae cover all cases.

\subsection{Formulae in the case $n=11$}
\label{modint11}

Our approach is similar to that in the last section. As one would expect
the formulae in the case $n=11$ are more complicated than those in 
the case $n=7$. There are however two further complications. One
as noted in Remark~\ref{rem:noc6} 
is the absence of a polynomial $c_6$. The other is
that the form we are looking for is no longer uniquely determined by its
image in the co-ordinate ring. Indeed in the case $n=7$ we were looking
for a cubic form, and in the case $n=11$ we are looking for a quintic
form. But in both cases the homogeneous ideal is generated by quartics.

Consideration of the action of $\SL_2(\Z/11\Z)$ on both the $z$-curve and
the $A$-curve suggests we start with the forms 
\begin{align*}
s_1(v,w,x,y,z) &=  
   (a^3 b c^3 + b^4 c d^2 - a b^2 c^2 d e - 2 b c^2 d e^3) v^2 w x z \\
& \quad + (b^3 c d^3 + c^4 d e^2 - a b c^2 d^2 e - 2 a^3 c d^2 e) v w^2 x y \\
& \quad + (c^3 d e^3 + a^2 d^4 e - a b c d^2 e^2 - 2 a b^3 d e^2) w x^2 y z \\
& \quad + (a^3 d^3 e + a b^2 e^4 - a^2 b c d e^2 - 2 a^2 b c^3 e) v x y^2 z \\
& \quad + (a b^3 e^3 + a^4 b c^2 - a^2 b^2 c d e - 2 a b^2 c d^3) v w y z^2 \\
& \quad \hspace{-3em} 
   + 2 (a^2 b^2 c^2 e + a^2 b^2 d e^2 + a^2 c d^2 e^2 + a b^2 c^2 d^2 
      + b c^2 d^2 e^2) v w x y z, \\
s_2(v,w,x,y,z) &= 
    a^2 b c e v^2 w x z + a b^2 c d v w^2 x y + b c^2 d e w x^2 y z 
+ a c d^2 e v x y^2 z  
\\ & \quad \hspace{3em} 
+ a b d e^2 v w y z^2 + 2 a b c d e v w x y z.
\end{align*}
We then solve for $r_1$ and $r_2$ satisfying
\begin{equation}
\label{reln11}
r_i(v,w,x,y,z) (vwxyz)^3 \equiv s_i(v,w,x,y,z)^4 \pmod{\I,\I'} 
\end{equation}
where $\I$ and $\I'$ are the homogeneous ideals for $X(11) \subset \PP^4$
with respect to the two sets of variables $a,b,c,d,e$ and $v,w,x,y,z$.
The coefficients of $r_1$ and $r_2$ are homogeneous polynomials of degrees
$28$ and $20$ in $a,b,c,d,e$. It is important to note that $r_1$ and $r_2$
are not uniquely determined by~(\ref{reln11}). However by 
averaging over the group we were able to choose
$r_i= (abcde)^3 \widetilde{r}_i$ in such a way that the coefficients of 
\[ \widetilde{r}_1(v \x_1 + w \x_4 + x \x_5 + y \x_9 + z \x_{14}) \]
and
\[ \widetilde{r}_2 (v \nabla F + w \nabla I_7 + x \nabla I_8+ y \nabla I_9
+ z \nabla c_4)  \]
are congruent mod $\I$ to certain polynomials in $F$ and $c_4$.
The result is a pair of quintic forms $\widetilde{d}_1(v,w,x,y,z)$
and $\widetilde{d}_2(v,w,x,y,z)$ with coefficients in $\Q[F,c_4]$.
We then put 
\begin{align*}
 d_1(v,w,x,y,z) & = \widetilde{d}_1(F v,w,F^7 x,F^2 y,F^4 z) \\
 d_2(v,w,x,y,z) & =\frac{1}{F^{4}}\widetilde{d}_2(F^2 v,F^8 w,F^4 x,y,F^3 z) 
\end{align*}
and replace $F^{11}$ by $\Delta$ so that $d_1$ and $d_2$ have coefficients
in $\Q[c_4,\Delta]$. 

\begin{Remark}
Unfortunately the polynomials $\widetilde{r}_i$ and 
$d_i$ would take several pages to print out, so we must refer the 
reader to the accompanying Magma file \cite{mywebsite}
for further details. We should also remark that the computation 
of $d_1$ and $d_2$ took several hours of computer time (whereas 
no other calculation up to this point took more than a few seconds).
\end{Remark}

\begin{Theorem} 
\label{getc4-11}
Let $E/K$ be an elliptic curve with 
Weierstrass equation $y^2 = x^3 - 27 c_4 x - 54 c_6$ and let
$\Delta = (c_4^3- c_6^2)/1728$. Assume $j(E) \not= 0,1728$ and let 
$X = X_E(11)$, respectively $X_E^-(11)$, be as given in Theorem~\ref{Thm11}. 
If $(v:w:x:y:z) \in X(K) \setminus \{d_i= 0\}$, not a cusp, then 
the corresponding elliptic curve $E'/K$ satisfies
\[  c_4(E') \equiv d_1(v,w,x,y,z) \, c_4(\BF)(v,w,x,y,z) \mod{(K^\times)^4}, \]
respectively
\[  c_4(E') \equiv d_2(v,w,x,y,z) \, c_4(\BG)(v,w,x,y,z) \mod{(K^\times)^4}. \]
\end{Theorem}

\begin{Proof}
As noted in \cite[Section 7.1]{HK} we are free to extend our field $K$ so 
that $\phi : E[11] \isom \mu_{11} \times \Z/11\Z$ is defined over $K$.
Let $(a:b:c:d:e)$ be the corresponding $K$-point on $X(11)$. We scale
$a,b,c,d,e$ so that $c_4(a,b,c,d,e) = c_4$. Then 
$a^{11}, \ldots, e^{11} \in K$ 
and by comparing the Weierstrass equation
for $E$ in the statement of the theorem with~(\ref{dagger11}) we deduce
$(abcde)^{11} \in (K^\times)^4$. The polynomials $\BF$ and $\BG$ were
computed in Section~\ref{sec:formulae11} as twists of $F$. Putting 
\begin{align*}
(v',w',x',y',z')^T & = v F^7 \x_1 + w F^6 \x_4 
 + x F^{13} \x_5 + y F^8 \x_9 + z F^{10} \x_{14}, \\
(v'',w'',x'',y'',z'')^T & = v F^3 \nabla F + w F^9 \nabla I_7 + x F^5 \nabla I_8
+ y F \nabla I_9 + z F^4 \nabla c_4,
\end{align*}
it follows by (\ref{det11}), (\ref{FG11}) and~(\ref{covprop11}) that
\begin{align*} 
c_4(\BF)(v,w,x,y,z) 
&= \frac{(c_4^3 - 1728 F^{11})^8}{F^{22}} c_4(v',w',x',y',z'), \\
c_4(\BG)(v,w,x,y,z) &= \frac{(55(c_4^3 - 1728 F^{11}))^8}{F^{11}} 
c_4(v'',w'',x'',y'',z''). 
\end{align*}
By construction of $d_1$ and $d_2$ we have
\begin{align*} 
d_1(v,w,x,y,z) &= \frac{1}{(abcde)^3 F^{30}} r_1(v',w',x',y',z'), \\
d_2(v,w,x,y,z) &= \frac{1}{(abcde)^3 F^{9}} r_2(v'',w'',x'',y'',z''). 
\end{align*}
In view of Theorem~\ref{weqns} our aim is to show that
\small
\begin{align*}
d_1(v,w,x,y,z) \, c_4(\BF)(v,w,x,y,z) & \equiv v'w'x'y'z' c_4(v',w',x',y',z') 
\mod{(K^\times)^4}, \\
d_2(v,w,x,y,z) \, c_4(\BG)(v,w,x,y,z) 
& \equiv v''w''x''y''z'' c_4(v'',w'',x'',y'',z'') 
\mod{(K^\times)^4},
\end{align*}
\normalsize
equivalently
\begin{align*}
(abcde)^8 F^{36} r_1(v',w',x',y',z') & \equiv v'w'x'y'z' \mod{(K^\times)^4}, 
\\
(abcde)^8 F^{24} r_2(v'',w'',x'',y'',z'') & \equiv v''w''x''y''z'' 
\mod{(K^\times)^4}. 
\end{align*}
To finish the proof we note that the quintic forms 
\[ (abcde)^2 F^9 s_1(v',w',x',y',z') \quad \text{ and } \quad 
(abcde)^2 F^6 s_2(v'',w'',x'',y'',z'')\] 
have coefficients in $K$ (since the degree of each coefficient 
is a multiple of $11$) and then use~(\ref{reln11}). 
\end{Proof}

We already gave a formula for the $j$-invariant in Section~\ref{sec:jinv}.
So (assuming $j(E') \not=0$) Theorem~\ref{getc4-11} determines $E'$ 
up to quadratic twist by $-1$.
In the case $K = \Q$ it is easy to decide which of the remaining 
two possibilities is correct by looking at traces of Frobenius.

In principle it should be possible to find alternative quintic forms
to be used at points where $d_1$ or $d_2$ vanishes. (The quintic
forms in question are those meeting the $z$-curve in a divisor $4D$ 
where $D$ is a hyperplane section for the $A$-curve.) In the case $n=7$ we
managed to find the alternative forms using the Riemann-Roch machinery
in Magma. Unfortunately the analogue of this in the case $n=11$ does
not appear to be practical. In the case of $X_E^-(11)$ this is not
a problem, since the $25$ points with $d_2=0$ correspond to the elliptic 
curves $\ell$-isogenous to $E$ for $\ell = 2,7,13$. We can also account 
for $7$ of the points on $X_E(11)$ with $d_1=0$ as corresponding to
the elliptic curve $E$ itself and the elliptic curves $5$-isogenous to $E$.
We are yet to encounter an example (over $K= \Q$) where one of the
remaining points with $d_1=0$ is rational.

\section{Examples}
\label{sec:ex}

We use the formulae in Theorems~\ref{MainThm7}, \ref{MainThm9} 
and~\ref{MainThm11} to give examples 
of non-trivial $n$-congruences for $n=7,9,11$ over $\Q$ and $\Q(T)$. 
By ``non-trivial'' we mean that the elliptic curves are not isogenous. 
The examples over $\Q$ illustrate the value
of minimising and reducing as described in Section~\ref{sec:minred}.
The examples over $\Q(T)$ were found by
setting $a=b=-27j/(4(j-1728))$ 
to obtain a surface fibred over the $j$-line and then intersecting
with one of the co-ordinate hyperplanes in the hope of finding 
a rational curve.
We refer to elliptic curves over $\Q$ by their labels in Cremona's tables
\cite{Cr}. For elliptic curves beyond the range of 
Cremona's tables we simply write the conductor followed by a ${\tt *}$.

\subsection{Examples in the case $n=7$}

\begin{Example} Let $E$ be the elliptic curve $162c1$.
Let $\FF$ and $\GG$ be the equations for $X_E(7)$ and $X_E^-(7)$ as given in
Theorem~\ref{MainThm7} with $a = 3645$ and $b = -13122$. These have 
invariants $\Psi(\FF) = -2^{11} \cdot 3^{18}$ and $\Psi(\GG) = 2^{22} \cdot 3^{36}$.
Minimising and reducing suggests that 
we substitute
\begin{align*}
F(x,y,z) & = \frac{1}{2^{10} 3^{14}} 
  \FF (36 y - 9 z, 1944 x - 972 y - 1215 z, z )  \\
G(x,y,z) &= \frac{1}{2^{12} 3^{20}} 
\GG (18 x + 18 y + 9 z,z,-486 x + 1458 y + 1944 z)
\end{align*}
to give quartics 
\begin{align*}
F(x,y,z) &= 3 x^3 z + 3 x^2 y^2 - 6 x^2 y z + 3 x^2 z^2 - 3 x y^3 \\ 
&~\text{\hspace{8em}}~
  + 3 x z^3 + 2 y^4 - y^3 z - 9 y^2 z^2 + 4 y z^3 - 5 z^4 \\
G(x,y,z) &= -x^3 y - x^3 z - 6 x^2 z^2 + 6 x y^2 z - 6 x y z^2 \\ 
&~\text{\hspace{8em}}~
   + 6 x z^3 + 2 y^4 + 2 y^3 z - 6 y^2 z^2 - 38 y z^3 - 8 z^4
\end{align*}
with invariants $\Psi(F) = -2 \cdot 3^4$ and $\Psi(G) = 2^2 \cdot 3^4$.
We find rational points $P_1 = (1:0:0)$, $P_2 = (3:-2:-1)$ 
on $\{F = 0\} \subset \PP^2$, and rational points $P_3 = (1:0:0)$, 
$P_4 = (1:1:-1)$, $P_5 = (4:-1:1)$ on $\{G=0\} \subset \PP^2$. 
The corresponding elliptic curves $7$-congruent to $E$ are
\begin{align*}
& P_1 & & 162c1 & y^2 + x y & = x^3 - x^2 + 3 x - 1 \\
& P_2 & & 293706* & y^2 + x y & = x^3 - x^2 - 62930562 x - 192134303740  \\
& P_3 & & 162c2 & y^2 + x y & = x^3 - x^2 - 42 x - 100 \\
& P_4 & & 17334f1 & y^2 + x y & = x^3 - x^2 - 5473977 x - 4956193171 \\
& P_5 & & 624186* & y^2 + x y & = x^3 - x^2 - 11751402282 x + 360746315347508. 
\end{align*}
The fact $162c1$ and $162c2$ are reverse $7$-congruent is already clear
since they are $3$-isogenous and $(3/7)=-1$.
\end{Example}

It is shown in \cite{HK} that there are infinitely many $6$-tuples of 
directly $7$-congruent non-isogenous elliptic curves over $\Q$.
The following example shows that there are infinitely many 
pairs of reverse $7$-congruent non-isogenous elliptic curves over $\Q$.

\begin{Example}
Let $E/\Q(T)$ be the elliptic curve $y^2 = x^3 + a x + b$
where $a=b= -27 j/(4 (j-1728))$ and $j = 27 T^3 (5 T-56)/(T-1)$.
Then on $X_E^-(7)$, with equation as given in Theorem~\ref{MainThm7}, 
we find the rational point
\[ (x:y:z) = (0 : -4 (T^2 - 12 T + 8) (5 T^2 + 4 T + 8): 
  9 T^2 (T + 4) (5 T-56)). \]
Specialising $T$ (and taking quadratic twists by $d$ as indicated) 
we obtain the following pairs of reverse $7$-congruent elliptic 
curves $E_1$ and $E_2$.
\[ \begin{array}{ccll}
T & d & \,\,\,\,\, E_1 & \,\,\,\,\, E_2 \\ \hline
-16 & -38 & 361a1 & 361a2 \\
8 & -10 & 700g1 & 2100q1 \\
2 & -2  & 2116b1 & 10580h1 \\
16/5 & -42 & 24255r1 & 24255m2
\end{array} \]
The existence of specialisations $E_1$ and $E_2$ that are not 
isogenous is enough to show that there are
infinitely many such specialisations.
\end{Example}

\subsection{Examples in the case $n=9$}

\begin{Example}
\label{ex9-3}
Let $E$ be the elliptic curve $47775z1$. Let $\FF_1 = \FF_2 = 0$ be the 
equations for $X_E(9) \subset \PP^3$ as given in Theorem~\ref{MainThm9} with 
$a = -41489280$ and $b = 102867483600$. The invariant is 
$\Psi(\FF_1,\FF_2) = -2^{42} \cdot 3^{60} \cdot 5^{12} \cdot 7^{32} \cdot 13^{4}$.
Minimising and reducing suggests that we substitute
\small
\[ \begin{pmatrix} x \\ y \\ z \\ t \end{pmatrix} \leftarrow
\begin{pmatrix}
 2520473760 & 937149484320 & -1998984627360 & -152410870080 \\
 0 & 79644600 & -185343480 & -3827880 \\ 
 0 & -22932 & 47040 & 6468 \\ 0 & -6 & 13 & 1 
\end{pmatrix} \begin{pmatrix} x \\ y \\ z \\ t \end{pmatrix} \]
\normalsize
so that $X_E(9)$ is defined by
\begin{align*}
-x^2 z + x^2 t + 4 x y z + 2 x y t - 3 x z^2 + 2 x z t - 3 x t^2  + 6 y^3 + 
14 y^2 z \hspace{6em} 
 & \\ +~y^2 t + 6 y z^2 - 4 y z t + 9 y t^2 - 6 z^3 + 27 z^2 t
 - 13 z t^2 - t^3 & = 0 \\
-3 x^2 y + 4 x^2 z + 3 x^2 t + 3 x y^2 + 20 x y z - 12 x y t - 3 x z^2  - 32 x z t
+ 25 x t^2  + 21 y^3 & \\ + 16 y^2 z - 24 y^2 t - 12 y z^2 + 100 y z t + 34 y t^2 + 
39 z^3 - 21 z^2 t - 56 z t^2 - 11 t^3 & = 0
\end{align*}
with invariant $-3^3 \cdot 5^6 \cdot 7^5 \cdot 13^4$.
We find rational points $P_1 = (1:0:0:0)$, $P_2 = (4:-1:-1:0)$ and 
$P_3 = (1:2:-1:0)$. The corresponding elliptic curves directly $9$-congruent 
to $E$ are
\begin{align*}
& P_1 & & 47775z1 & y^2 + y &= x^3 - x^2 - 32013 x + 2215478 \\
& P_2 & & 429975*  & y^2 + y &= x^3 - 314688780 x - 2148671872069 \\ 
& P_3 & & 494901225*  & y ^2 + y &= x^3 - 23634650164230 x - 21037908383222056594 
\end{align*}
Since $X_E^-(9)$ is not locally soluble at $p=7$ there
are no elliptic curves reverse $9$-congruent to $E$. 
\end{Example}

In addition to Example~\ref{ex9-3} we have found two further triples of 
directly $9$-congruent non-isogenous elliptic curves over $\Q$. These are
\[ \begin{array}{ll}
4650j1  & y^2 + x y = x^3 + x^2 - 2700 x + 54000  \\
553350*   & y^2 + x y = x^3 + x^2 - 10472207700 x - 
    455228489646000  \\
1966950*   & y^2 + x y = x^3 - x^2 - 20654522386242 x - 
    36130051534030639084  \bigskip \\
27606c1  & y^2 + x y = x^3 - 10289707 x + 12703497719  \\
358878*   & y^2 + x y = x^3 + 2940333 x - 1416695391  \\
1242270*   & y^2 + x y + y = x^3 - x^2 - 359912 x - 322105301 
\end{array} \]
The elliptic curves $1701a1$, $1701g1$ and $22113c1$ are also 
$9$-congruent but only the last two of these
are directly $9$-congruent.

\begin{Example}
Let $E$ be the elliptic curve $201c1$. Let $\GG_1 = \GG_2 = 0$ be the 
equations for $X_E^-(9) \subset \PP^3$ as given in Theorem~\ref{MainThm9} with 
$a = -1029699$ and $b = 402173694$. The invariant is 
$\Psi(\GG_1,\GG_2) = 2^{48} \cdot 3^{85} \cdot 67^5$.
Minimising and reducing suggests that we substitute
\small
\[ \begin{pmatrix} x \\ y \\ z \\ t \end{pmatrix} \leftarrow
\begin{pmatrix}
 -26471709 & -23136696 & 20106774 & -20376135 \\ 
  -45147 & -39828 & 33990 & -34509 \\ 
  90294 & 79332 & -68304 & 69342 \\ 77 & 68 & -58 & 59 
\end{pmatrix} \begin{pmatrix} x \\ y \\ z \\ t \end{pmatrix} \]
\normalsize
so that $X_E^-(9)$ is defined by
\begin{align*}
-x^3 + 4 x^2 y + 3 x^2 z - x^2 t + 6 x y^2 + 2 x y z - 2 x y t - 6 x z^2 + 
4 x z t \hspace{2em} 
 & \\ - 11 x t^2 + y^3 + 7 y^2 t - 2 y z^2 + 4 y z t - 4 y t^2 + 6 z^3 - 
7 z^2 t + 4 z t^2 + t^3 & = 0 \\
2 x^3 - x^2 y + 5 x^2 t - 10 x y^2 - 2 x y z + 16 x y t - 3 x z^2 + 4 x z t + 
8 x t^2 \hspace{2em} 
 & \\ - 5 y^3 - y^2 z - 3 y^2 t - y z^2 - 2 y z t + 12 y t^2 + 3 z^3 - 4 z^2 t 
+ 2 z t^2 - 3 t^3 & = 0 
\end{align*}
with invariant $-3^4 \cdot 67^5$. On this curve we find the rational point
$(1:-2:-1:0)$. The corresponding elliptic curve reverse 
$9$-congruent to $E$ is 
\begin{align*}
& 374865* & y^2 + x y = x^3 + x^2 - 60068738107 x + 4858035498982726 
\end{align*}
\end{Example}

The following example shows that there are infinitely many 
pairs of directly $9$-congruent non-isogenous 
elliptic curves over $\Q$.

\begin{Example} Let $E/\Q(T)$ be the elliptic curve 
$y^2 = x^3 + a(T) x + b(T)$ where
\begin{align*}
a(T) & = \tfrac{1}{2} (39 T^4 - 60 T^3 - 162 T^2 + 60 T + 39), \\
b(T) & = 47 T^6 + 120 T^5 + 21 T^4 + 21 T^2 - 120 T + 47.
\end{align*}
Then on $X_E(9)$, with equations as given in Theorem~\ref{MainThm9}, 
we find the rational point
\[ (x:y:z:t) = (\tfrac{15}{2}(3 T^4 + 8 T^3 - 2 T^2 - 8 T + 3)
: T^2 +1 : 1 :0).\]
The corresponding curve directly $9$-congruent to $E$ is the 
curve directly $3$-congruent to $E$ constructed in Theorem~\ref{MainThm3} with
$c_4 = -a(T)/27, c_6 = -b(T)/54$ and 
\small
\begin{align*}
(\la:\mu) &= \big(47 T^6 - 78 T^5 - 153 T^4 + 244 T^3 + 153 T^2 - 78 T - 47  
 \\ & \hspace{18em} : 18 (T^2 + 1) (T^2 + 6 T - 1)\big).
\end{align*}
\normalsize
Specialising to $T=0$ gives a pair of curves with conductors
$80640$ and $5886720$. In particular these curves are not isogenous.
\end{Example}

Next we give an example to show there are infinitely many non-trivial
pairs of reverse $9$-congruent elliptic curves over $\Q$.

\begin{Example}
Let $E/\Q(T)$ be the elliptic curve $y^2 = x^3 + a(T) x + b(T)$ 
where
\begin{align*}
a(T) & = 3 (3 T + 1) (6 T^3 - 3 T - 1) (9 T^3 - 9 T - 4)^2, \\
b(T) & = 2 (3 T^3 + 27 T^2 + 21 T + 4) (6 T^3 - 3 T - 1)^2 (9 T^3 - 9 T - 4)^2. 
\end{align*}
Then $X^-_E(9)$, with equations 
as given in Theorem~\ref{MainThm9}, has rational point
\[ (x:y:z:t) = (-(6 T^3-3 T-1) (9 T^3 - 9 T - 4):T:1:0). \]
The corresponding curve reverse $9$-congruent to $E$ is the curve
reverse $3$-congruent to $E$ constructed in Theorem~\ref{MainThm3} with
$c_4 = -a(T)/27, c_6 = -b(T)/54$ and 
\small
\begin{align*}
(\la:\mu) &= \big( (3 T+1) (9 T^3-9 T-4) (6 T^3-3 T-1)
               (180 T^4 + 321 T^3 + 216 T^2 + 66 T + 8) \\
& \hspace{5em} 
: 3 (369 T^6 + 1107 T^5 + 1431 T^4 + 1017 T^3 + 414 T^2 + 90 T + 8) \big). 
\end{align*}
\normalsize
Specialising $T$ (and taking quadratic twists by $d$ as indicated) 
we obtain the following pairs of reverse $9$-congruent elliptic 
curves $E_1$ and $E_2$.
\[ \begin{array}{ccll}
T & d & \,\,\,\,\, E_1 & \,\,\,\,\, E_2 \\ \hline
-1 & 6 & 24a4 & 24a5 \\ 
-1/3 & 6 & 243a1 & 243b2 \\ 
-1/4 & 3 & 768d1 & 114432o1 \\ 
-1/2 & -6 & 6400u1 & 6400u2 \\ 
-2/3 & 6 & 23814v1 & 23814i1 
\end{array} \]
\end{Example}

\subsection{Examples in the case $n=11$}

\begin{Example}
Let $E$ be the elliptic curve 1782b1. Let $\FF$ be the cubic form 
describing $X_E(11) \subset \PP^4$ as given in Theorem~\ref{MainThm11} with 
$a = 765$ and $b = 15102$. The invariant is 
$\Psi(\FF) = -2^{28} \cdot 3^{12} \cdot 11^6$. 
Minimising and reducing suggests that 
we substitute
\small
\[ \begin{pmatrix} v \\ w \\ x \\ y \\ z \end{pmatrix} \leftarrow
\begin{pmatrix}
  984 & 12900 & -9093 & -34056 & 13689 \\
  -2040 & -24252 & -3315 & 0 & -16857 \\ 
  328 & 164 & -435 & 0 & -57 \\ 
  -352 & 88 & -264 & 264 & -1056 \\ 
  -8 & -4 & -13 & 0 & 25
\end{pmatrix} \begin{pmatrix} v \\ w \\ x \\ y \\ z \end{pmatrix} \]
\normalsize
so that $X_E(11) \subset \PP^4$ is the singular locus of the Hessian of
\begin{align*} -v^2 w + v^2 x - v^2 y + 2 v^2 z - v w^2 + 4 v w z 
  - 4 v x^2 - 8 v x y + 2 v x z + 6 v y z \hspace{2em} 
\\ + 3 v z^2 + 2 w^3  - 3 w^2 x
  - 2 w^2 y + 8 w^2 z + 6 w x^2 + 2 w x y + 2 w x z + 6 w y^2 - 6 w y z \\
  + 9 w z^2 - x^3 - x^2 z - 3 x y^2 - 6 x y z - 9 x z^2
  - 6 y^3 + 9 y^2 z + 3 y z^2 - 7 z^3 & = 0 
\end{align*}
with invariant $2^2 \cdot 3^4 \cdot 11^2$. 
We find rational points $P_1 = (-1:5:1:2:1)$, $P_2 = (0:0:0:1:0)$
and $P_3 = (1:1:-1:0:-4)$. The corresponding elliptic curves directly 
$11$-congruent to $E$ are
\begin{align*}
& P_1 & & 1782b1 &  y^2 + x y &= x^3 - x^2 + 48 x + 224 \\
& P_2 & & 1782b2  & y^2 + x y &= x^3 - x^2 - 447 x - 7795 \\
& P_3 & & 447282*  & y^2 + x y &= x^3 - x^2 - 17552171922 x - 227953575178678
\end{align*} 
The fact $1782b1$ and $1782b2$ are directly $11$-congruent is already
clear since they are $3$-isogenous and $(3/11)=1$.
\end{Example}

\begin{Example}
Let $E$ be the elliptic curve 4466c1. Let $\GG$ be the cubic form 
describing $X_E^-(11) \subset \PP^4$ as given in Theorem~\ref{MainThm11} with 
$a = 85$ and $b = -83162$. The invariant is 
$\Psi(\FF) = 2^{21} \cdot 7 \cdot 11^2 \cdot 29^2$.
Minimising and reducing suggests that 
we substitute
\small
\[ \begin{pmatrix} v \\ w \\ x \\ y \\ z \end{pmatrix} \leftarrow
\begin{pmatrix}
4096 & -1408 & 128 & -1312 & 45088 \\ 
0 & 128 & 128 & 32 & 110 \\ 0 & 0 & -256 & -96 & -103 \\
0 & 0 & 0 & -32 & -11 \\ 0 & 0 & 0 & 0 & -1 
\end{pmatrix} \begin{pmatrix} v \\ w \\ x \\ y \\ z \end{pmatrix} \]
\normalsize
so that $X_E^-(11) \subset \PP^4$ is the singular locus of the Hessian of
\begin{align*}
  - 2 v^2 z - 4 v w y + 12 v x y + 4 v x z + 5 v y^2 + 6 v y z - 43 v z^2 
  - w^2 x + w^2 y \hspace{2em} \\ - 4 w x y - 2 w x z - 3 w y^2 
  + 196 w y z + 83 w z^2 - 11 x^3 - 12 x^2 y - 9 x^2 z \hspace{2em} \\ 
  - 11 x y^2 + 366 x y z + 125 x z^2 + 322 y^3 
  + 447 y^2 z + 275 y z^2 + 632 z^3 &= 0
\end{align*}
with invariant $-2^2 \cdot 7 \cdot 11^2 \cdot 29^2$.
We find rational points $P_1 = ( -7 : 11 : 3 : 1 : 1)$
and $P_2 = ( 7830 : -3553 : 510 : -281 : 71 )$.
The corresponding elliptic curves reverse $11$-congruent to $E$ are
\begin{align*}
& P_1 & & 4466c2 &  y^2 + x y + y &= x^3 - x^2 - 1755 x - 27349 \\
& P_2 & & 1174558* & y^2 + x y + y &= x^3 - x^2 + 117885809240 x + 
16240157710556505 
\end{align*} 
The fact $4466c1$ and $4466c2$ are reverse $11$-congruent is already
clear since they are $2$-isogenous and $(2/11)=-1$.
\end{Example}

We did not find any triples of $11$-congruent non-isogenous
elliptic curves over~$\Q$. The following example 
shows that there are infinitely many pairs of 
directly $11$-congruent non-isogenous elliptic curves over $\Q$.
\begin{Example}
Let $E/\Q(T)$ be the elliptic curve $y^2 = x^3 + a(T) x + b(T)$ 
where
\begin{align*}
a(T) &= -3 (T - 3) (T^4 - 5 T^2 - 24 T - 92)/(T^3 - T^2 + 4 T + 24)  \\
b(T) &= -2 (T - 3) (T^5 - T^4 - 11 T^3 - 43 T^2 - 62 T - 316)
                 /(T^3 - T^2 + 4 T + 24).
\end{align*}
Then $X_E(11)$, with equations 
as given in Theorem~\ref{MainThm11}, has rational point
\[ \begin{pmatrix} v \\ w \\ x \\ y \\ z \end{pmatrix}
= \begin{pmatrix}
T^6 + T^5 + 31 T^4 + 259 T^3 + 520 T^2 + 676 T + 1248 \\
-(T - 3) (T^5 + 4 T^4 + 43 T^3 + 100 T^2 - 44 T - 320) \\
-(T^2 + 3 T + 14) (T^3 - T^2 + 4 T + 24) \\
0 \\
(T + 4) (T^3 - T^2 + 4 T + 24)
\end{pmatrix} \]
Specialising $T$ (and taking quadratic twists by $d$ as indicated) 
we obtain the following pairs of directly $11$-congruent elliptic 
curves $E_1$ and $E_2$.
\[ \begin{array}{ccll}
T & d & \,\,\,\,\, E_1 & \,\,\,\,\, E_2 \\ \hline
2 & -6 & 11a3 & 11a2 \\ 
1 & 42 & 49a1 & 49a4 \\ 
-3 & -2 & 216b1 & 1512c1 \\ 
11 & -426 & 10082c1 & 70574h1 
\end{array} \]
The elliptic curve $11$-congruent to $E$ is 
$y^2 = x^3 + A(T) x + B(T)$ where
\small
\begin{align*}
A(T) &= -3 (T - 3) (T^2 - 8 T - 17) (T^3 - T^2 + 4 T + 24) 
  (T^{12} - 250 T^{11} + 3473 T^{10} \\ & \quad
   - 23824 T^9 + 106654 T^8 - 354556 T^7 + 890186 T^6 
        - 1710568 T^5 \\ & \quad  + 2386357 T^4 - 2054170 T^3 + 1799781 T^2 + 956680 T + 
       3570796), \\
B(T) &= -2 (T - 3) (T^3 - T^2 + 4 T + 24)^2 
  (T^{20} + 476 T^{19} - 27815 T^{18} + 556718 T^{17} \\ & \quad
    - 6046664 T^{16} + 42450848 T^{15} - 
        213832636 T^{14} + 823702888 T^{13} \\ & \quad 
- 2497998850 T^{12} + 5954643736 T^{11} - 
        10798748818 T^{10}  + 13644339892 T^9 \\ & \quad - 7927895108 T^8 - 10398245632 T^7 + 
        25581636532 T^6 - 10366268760 T^5 \\ & \quad - 60876061719 T^4 
    + 164062110060 T^3 - 
      98120800447 T^2 + 262948421518 T \\ & \quad + 141270230564).
\end{align*}
\normalsize
These elliptic curves have discriminants
\small
\begin{align*}
   & 2^{12} 3^6 (T - 5) (T - 3)^2 (T + 1)^5 
     (T^2 + 7)/ (T^3 - T^2 + 4 T + 24)^3, \\
   -& 2^{12} 3^6 (T - 5)^4 (T - 3)^2 (T + 1)^3 (T^2 + 7)
   (T^3 - T^2 + 4 T + 24)^3 (T^3 - T^2 + 15 T - 31)^{11}.
\end{align*}
\normalsize
\end{Example}

We did not find any pairs of reverse 11-congruent non-isogenous 
elliptic curves over $\Q(T)$. We note that 
according to \cite[Theorem 4]{KS} the modular diagonal surface 
in this case is of general type.

\subsection{Tables}
We have written a program in Magma that given an elliptic curve
$E/\Q$ and $n \in \{7,9,11\}$ searches for elliptic curves $n$-congruent
to $E$. For $n \in \{9,11\}$ we have run these programs for all elliptic
curves in the Cremona database (up to conductor $130000$ at the time of
our calculation). 
The resulting list of pairs of $n$-congruent elliptic curves
is available from the author's website \cite{mywebsite}. 
We have been careful to remove all pairs that could be deduced from 
earlier entries by any combination of the following observations. 
\begin{itemize}
\item $n$-congruence of elliptic curves is an equivalence relation. 
\item If $\phi : E \to E'$ is a rational isogeny of
degree coprime to $n$ then $E$ and $E'$ are $n$-congruent. 
\item If $E$ and $E'$ are $n$-congruent then so are their 
twists by the same quadratic character. 
\end{itemize}

There is no guarantee that our tables are complete, since 
we have only searched for points of small height on the corresponding curves of
genus $10$ and $26$. However the points we found (on minimised and reduced
models) were generally much smaller than the search bound used. 
So in any given case it 
seems likely  that we have found all the rational points, but
there is little prospect of proving this for curves of such large genus. 

\tiny
\[ \hspace{-0.1em} \begin{array}{lll|lll|lll|lll}
\multicolumn{12}{c}{\text{\small{Table 1 : Pairs of $9$-congruent elliptic curves}}} \\
17a^{ 2 } & + & 493a & 430b & + & 20210b & 968a & + & 132616* & 1950h & - & 122850l \\ 
33a^{ 2 } & - & 297c & 434c^{ 2 } & + & 62062h & 1050j2 & - & 11550bg2 & 1950m2 & + & 17550ba2 \\ 
35a2 & - & 77b2 & 446b & - & 150302* & 1066a & - & 7462e & 1952c & + & 44896e \\ 
35a3 & - & 1015b1 & 456c & - & 20520e & 1106c & + & 2281678* & 1963a & + & 374933* \\ 
66a1^{ 2 } & + & 3102c1 & 459c1 & + & 2295a1 & 1110n2 & + & 13098q2 & 1998d & - & 3774c \\ 
84a3^{ 2 } & - & 1932b2 & 506b & + & 2530b & 1111a & + & 121215655* & 2040h & + & 4073880* \\ 
91a & + & 5005c & 525c^{ 2 } & + & 4725r & 1134d1 & + & 5670h1 & 2072c^{ 2 } & - & 184408* \\ 
106d & + & 2438a & 537b & + & 148749* & 1155j1 & + & 47355s1 & 2093c & + & 2576483* \\ 
110c2 & + & 60610m2 & 570k3^{ 2 } & - & 1254j2 & 1176a & + & 10584n & 2118b1 & - & 7660806* \\ 
115a & + & 366505* & 573c & - & 21491511* & 1176b & + & 34104b & 2135d & + & 207095* \\ 
118c & + & 27494c & 600c^{ 2 } & + & 5400r & 1190c2 & - & 393890* & 2190b & + & 116070i \\ 
123b & + & 3813a & 606f^{ 5 } & + & 155742* & 1209a & + & 1357707* & 2190e & + & 15330k \\ 
131a & - & 136633* & 627b2 & + & 5643d2 & 1215b & - & 1215g & 2190l & + & 24090k \\ 
140a2 & - & 3220c2 & 648b & + & 26568a & 1215c & + & 23085b & 2198b1 & + & 964922* \\ 
142c^{ 2 } & + & 232454* & 696a & + & 355656* & 1218h^{ 2 } & + & 1218i & 2209a & - & 11045b \\ 
153a & - & 117351a & 702f & + & 30186h & 1275d & + & 263925* & 2211c & + & 134871* \\ 
162a1 & + & 810c1 & 710a & + & 137030* & 1281d & + & 4209c & 2314b & - & 363298* \\ 
166a & + & 848426* & 714c & + & 6426n & 1288e & + & 240856* & 2373f & - & 21357c \\ 
174a1 & + & 17574f1 & 715b & + & 1080365* & 1330d1 & + & 860510* & 2418a & + & 2884674* \\ 
174c & + & 1914g & 741b & + & 3335241* & 1470e2 & + & 4998j2 & 2443c & + & 925897* \\ 
195c & + & 1755a & 741d & + & 65949e & 1470n^{ 2 } & + & 13230dv & 2450l^{ 13 } & + & 129850r \\ 
200b^{ 2 } & - & 20600v & 768b^{ 2, 5 } & + & 114432j & 1482e1 & + & 16302l1 & 2451b^{ 2 } & + & 2451c \\ 
201a & + & 82209a & 781a & - & 1349b & 1482e2 & + & 117078p2 & 2530e & + & 908270* \\ 
201c & - & 374865* & 805b^{ 2 } & + & 10465a & 1482f2 & + & 8128770* & 2537b & + & 210571* \\ 
229a & + & 2357555* & 805c^{ 2 } & + & 24955b & 1518r & + & 13662m & 2541a & - & 22869o \\ 
235b & + & 329a & 808a & + & 891224* & 1617b & + & 289443* & 2568e & + & 7408680* \\ 
238d^{ 2 } & - & 91154b & 810a2 & + & 208170* & 1653a & - & 520695* & 2571a & + & 12749589* \\ 
243a1 & - & 243b2 & 810d2 & + & 23490n2 & 1701a & - & 1701g & 2616c & + & 6574008* \\ 
246c^{ 2 } & + & 2214a & 843a & + & 1600857* & 1701a & - & 22113c & 2646a1 & + & 60858h1 \\ 
249a & + & 42579e & 850c & + & 16150q & 1725i^{ 2 } & - & 343275* & 2646e1 & + & 13230e1 \\ 
258c & + & 1290h & 858d2 & + & 4290l2 & 1757a & + & 335587* & 2674e & + & 259378* \\ 
270a1 & + & 2970i1 & 861b & + & 192003* & 1771c1 & + & 2807035* & 2678a & - & 20806a \\ 
302b^{ 2 } & + & 226802* & 861d & + & 469245* & 1782c & + & 290466* & 2699b & + & 5249555* \\ 
338a^{ 7 } & + & 22646d & 867b^{ 2 } & + & 141321* & 1806a & + & 2775822* & 2706i^{ 2 } & + & 65855922* \\ 
350b2 & - & 6650bh2 & 897d^{ 2 } & + & 8073a & 1830h & + & 16470w & 2728d & + & 832040* \\ 
354a^{ 2 } & + & 55578e & 906b & + & 610222710* & 1848g & + & 9097704* & 2758b & + & 1069527578* \\ 
364a & + & 26572d & 906e & + & 28086d & 1862a2 & + & 102410n2 & 2787d & + & 393187173* \\ 
390d3^{ 2 } & + & 3510f2 & 930d^{ 2 } & + & 564510* & 1870e & + & 136510* & 2793b & + & 790419* \\ 
406b2 & - & 62930k2 & 930i2 & + & 160890* & 1887b & - & 1054833* & 2835a1 & + & 65205g1 \\ 
430a & + & 7310e & 930l & + & 8370n & 1911a & - & 17199o & 2835d & + & 53865d 
\end{array} \hspace{-0.1em} \]
\normalsize

\tiny
\[ \hspace{-0.3em} \begin{array}{lll|lll|lll|lll}
\multicolumn{12}{c}{\text{\small{Table 2 : Pairs of $11$-congruent elliptic curves}}} \\
138b^{ 2, 3 } & + & 17526b & 3910l & + & 43010o & 16354h & + & 703222* & 55470m^{ 2 } & + & 1275810* \\ 
190b & + & 2470a & 3990z & + & 3990ba & 16390m & + & 23585210* & 55594a & - & 10062514* \\ 
216b & + & 1512c & 4046n & + & 416738* & 17550a & + & 298350* & 56100j^{ 2 } & - & 15539700* \\ 
238b^{ 2 } & + & 4522b & 4158f & + & 54054m & 17550z & + & 122850r & 56154t & - & 27346998* \\ 
258b^{ 2 } & + & 461562* & 4200g^{ 2 } & + & 79800h & 17556d & + & 9778692* & 57354j & + & 3498594* \\ 
294a^{ 7 } & + & 4998ba & 4275a^{ 2 } & - & 175275* & 17710h & + & 726110* & 58968k & + & 58968l \\ 
325a^{ 3 } & + & 23075e & 4466c^{ 2 } & - & 1174558* & 18850i & + & 131950* & 60088b & + & 53658584* \\ 
329a & + & 59549a & 4598b & + & 32186g & 19950bl & + & 618450* & 64596b & + & 53162508* \\ 
426b^{ 2 } & + & 77106c & 4704a & - & 192864* & 20838f & - & 10273134* & 65220c & + & 534086580* \\ 
497a & + & 148603* & 4760b & + & 775880* & 21175b & + & 21175g & 66300q & + & 1259700* \\ 
513a & + & 77463a & 5070d & - & 55770j & 21760a & + & 805120* & 66930e & + & 101265090* \\ 
600b & + & 4200c & 5265b & + & 194805* & 22950bi & + & 849150* & 67158c & + & 3559374* \\ 
645f & - & 70305g & 5454b & + & 59994e & 23595g & - & 165165* & 67650v & + & 16168350* \\ 
648b & - & 4536c & 5577b & - & 39039b & 24906a & + & 70558698* & 70890w & + & 2197590* \\ 
700e^{ 2 } & - & 11900d & 5880b & - & 417480* & 25110j & - & 426870* & 71610bj & + & 109634910* \\ 
1080d & + & 39960h & 6710h & + & 4381630* & 26026d & - & 2836834* & 73255b & - & 89151335* \\ 
1115a & + & 125995a & 7315e & + & 8214745* & 27610f^{ 2 } & + & 5273510* & 73326t & + & 3886278* \\ 
1134b & - & 12474f & 7350n & + & 95550ch & 27650j & + & 1023050* & 73370h & + & 130790* \\ 
1155b & - & 42735a & 7830f^{ 3 } & + & 477630* & 27885p & - & 2704845* & 78650w^{ 2 } & - & 78650bc \\ 
1176a & - & 15288a & 8096a & - & 13984a & 28420d & + & 5087180* & 79800s & + & 16837800* \\ 
1210d & + & 1210f & 8120d & + & 332920* & 28830a & - & 2219910* & 80325u & + & 80325v \\ 
1254b & + & 264594* & 8410c & - & 58870e & 28840c & + & 30253160* & 81675q & + & 3348675* \\ 
1470j & - & 16170bl & 8670o & - & 8670p & 30774h & + & 11417154* & 83148g & - & 582036* \\ 
1666a & + & 18326b & 9438d & + & 1009866* & 30855a & + & 215985* & 84270j & + & 589890* \\ 
1782b^{ 3 } & + & 447282* & 9450l & + & 9450m & 30888g & + & 216216* & 84930m & + & 59196210* \\ 
1848f & + & 97944h & 9450o & + & 292950* & 31746q & + & 222222* & 84930p & + & 31509030* \\ 
1870c & + & 24310c & 9450bo & + & 160650* & 32802d^{ 2 } & - & 117660774* & 85050bf & + & 2466450* \\ 
1925g & + & 140525* & 9555a & - & 181545* & 32856f & + & 229992* & 87210l & + & 3226770* \\ 
2093b & + & 26611d & 9966h & + & 827178* & 34385a & + & 378235* & 87990d & + & 5895330* \\ 
2184c & + & 2184d & 10010t & + & 670670* & 35574b & - & 35574e & 88842l & + & 120198k \\ 
2755b & + & 476615* & 10082b & + & 70574c & 36498s & + & 56406i & 90986e & + & 423175886* \\ 
3234e & + & 947562* & 11774a & + & 24172022* & 37446i & - & 4456074* & 92950s & + & 9016150* \\ 
3322c & + & 235862* & 12376j & + & 12376k & 37830o & + & 31285410* & 93795r & + & 2157285* \\ 
3325a & + & 1173725* & 12696b & + & 12696c & 38070k & - & 5291730* & 95370ch & + & 667590*^{ 2 } \\ 
3610d & + & 25270j & 13650ba & + & 2716350* & 39732a & + & 24514644* & 95550cb & + & 95550cc \\ 
3718d & + & 107822e & 14196i & - & 269724* & 45738b & + & 869022* & 95700b & + & 2201100* \\ 
3822u & - & 42042ce & 14844a & + & 40182708* & 47974d^{ 2 } & - & 48885506* & 98070z & + & 270771270* \\ 
3850a & + & 65450b & 14910ba & + & 40450830* & 53290a^{ 2 } & + & 53290b & 98982c & - & 19697418* \\ 
3850i^{ 2 } & + & 65450i & 14950l & + & 1599650* & 53865e & - & 10395945* & 104910d & + & 3252210* \\ 
3900e & - & 66300t & 16150p & + & 4731950* & 54813a & - & 1589577* & 105522e & + & 139183518* 
\end{array}  \hspace{-0.3em} \]
\normalsize

The entries are given as $E_1 \pm E_2$ where the $\pm$ indicates
whether the $n$-congruence respects the Weil pairing.
In all cases except when $n=9$ and the curves admit a rational $3$-isogeny 
we have omitted the final number from the Cremona reference (which may be 
taken to be $1$). A superscript $\ell$ indicates an elliptic curve that admits a
rational $\ell$-isogeny where $\ell$ is a prime not dividing $n$.
If $\ell$ is not a square mod $n$ then we may change the sign of the
congruence by passing to an isogenous curve.
For elliptic curves beyond the range of Cremona's tables we have again
written the conductor followed by a $*$. The extended version of 
our tables~\cite{mywebsite} also gives the Weierstrass equations.

\end{document}